\documentclass[preprint,12pt]{elsarticle}
\usepackage[margin=2.5cm]{geometry}

\usepackage{graphicx}
\usepackage{amssymb}

\usepackage[toc,page]{appendix}

\usepackage[centertags]{amsmath}
\usepackage{amsfonts}
\usepackage{amssymb}
\usepackage{amsthm}
\usepackage{mnsymbol}

\usepackage{epstopdf}
\usepackage{subfig}
\usepackage{pgf,tikz}
\usetikzlibrary{calc,patterns,angles,quotes,intersections,patterns}
\usepackage{tkz-euclide}
\usetkzobj{all} 

\definecolor{green(ryb)}{rgb}{0.173, 0.627, 0.173}
\definecolor{azure}{rgb}{0.0, 0.5, 1.0}

\usepackage{verbatim}
\usepackage[norelsize,ruled,vlined,linesnumbered,titlenumbered]{algorithm2e}
\makeatletter
\newcommand{\nosemic}{\renewcommand{\@endalgocfline}{\relax}}
\newcommand{\dosemic}{\renewcommand{\@endalgocfline}{\algocf@endline}}
\let\oldnl\nl
\newcommand{\nonl}{\renewcommand{\nl}{\let\nl\oldnl}}
\makeatother

\newtheorem{thm}{Theorem}[section]
\newtheorem{prop}[thm]{Proposition}
\newtheorem{dfn}[thm]{Definition}
\newtheorem{lem}[thm]{Lemma}
\newtheorem{oss}[thm]{Remark}
\newtheorem{ex}[thm]{Example}
\newtheorem{cor}[thm]{Corollary}

\def\RR{{\mathbb{R}}}

\newcommand{\eps}{\varepsilon}
\newcommand{\R}{\Rightarrow}
\newcommand{\cM}{{\mathcal{M}}}
\newcommand{\cB}{{\mathcal{B}}}
\newcommand{\SSS}{{\mathbb{S}}}

\newcommand{\cE}{{\mathcal{E}}}
\newcommand{\cF}{{\mathcal{F}}}
\newcommand{\sse}{\Leftrightarrow}

\journal{Computer Aided Geometric Design}

\begin{document}

\begin{frontmatter}


\title{Linear dependence of bivariate Minimal Support and\\ Locally Refined B-splines over LR-meshes\tnoteref{funding}}
\tnotetext[funding]{This work has received funding from the European Union's Horizon 2020 research
and innovation programme under grant agreement No 675789.}


\author[Si,uio]{Francesco Patrizi\corref{cor1}}
\ead{francesco.patrizi@sintef.no}
\author[Si]{Tor Dokken}
\ead{tor.dokken@sintef.no}
\address[Si]{Department of Applied Mathematics and Cybernetics, SINTEF, Forskningsveien 1, 0373 Oslo, Norway}
\address[uio]{Department of Mathematics, University of Oslo, Moltke Moes vei 35, 0851 Oslo, Norway}
\cortext[cor1]{Corresponding Author}

\begin{abstract}
The focus on locally refined spline spaces has grown rapidly in recent years due to the need in Isogeoemtric analysis (IgA) of spline spaces with local adaptivity: a property not offered by the strict regular structure of tensor product B-spline spaces. However, this flexibility sometimes results in collections of B-splines spanning the space that are not linearly independent. 
In this paper we address the minimal number of B-splines that can form a linear dependence relation for Minimal Support B-splines (MS B-splines) and for Locally Refined B-splines (LR B-splines) on LR-meshes.  We show that the minimal number is six for MS B-splines and eight for LR B-splines. 
Further  results are established to help detecting collections of B-splines that are linearly independent.
\end{abstract}


\end{frontmatter}


\section{Introduction}

In 2005 Thomas J.R. Hughes et al. \cite{hughes2005isogeometric} proposed to reconstitute finite element analysis (FEA) within the geometric framework of CAD technologies. This gave rise to Isogeometric Analysis (IgA). It unifies the fields of CAD and FEA by extending the isoparametric concept of the standard finite elements to other shape functions, such as B-splines and non-uniform rational B-splines (NURBS), used in CAD. This does not only allow for an accurate geometrical description, but it also improves smoothness properties. As a consequence, IgA methods often reach a required accuracy using a much smaller number of degrees of freedom \cite{grossmann2012isogeometric}. Moreover, in some situations, the increased smoothness also improves the stability of the approximations resulting in fewer nonphysical oscillations \cite{hughes2005isogeometric, manni2011isogeometric}.

However, in numerical simulations, local (adaptive) refinements are frequently used for balancing accuracy and computational costs. 
Traditional B-splines and NURBS spaces are formulated as tensor products of univariate B-spline spaces. This means that refining in one of the univariate B-spline spaces will cause the insertion of an entire new row or column of knots in the bivariate spline space, resulting in a global refinement. In order to break the tensor product structure of the underlying mesh,  new formulations of multivariate B-splines have been introduced addressing local refineability. 
\subsection{Overview of locally refined spline methods}\label{overview-local}
The first locally refined method introduced were the Hierarchical B-splines, or HB-splines \cite{forsey1988hierarchical}, whose properties were further analyzed in \cite{kraft1997adaptive}. The HB-splines are linearly independent and non-negative. However, partition of unity, which is a necessary for the convex hull property (essential for interpreting the B-spline coefficients as control points), was still missing. To rectify this, Truncated Hierarchical B-splines, or THB-splines, were proposed in \cite{giannelli2012thb} and further analyzed in \cite{giannelli2014strongly}. Here they show how the construction of HB-splines can be modified while preserving the properties of HB-splines, gaining the partition of unity and smaller support of the basis functions.

A different approach, for local refinements, was introduced in \cite{sederberg2003t} with the T-splines. These are defined over \emph{T-meshes}, where T-junctions between axis aligned segments are allowed. T-splines have been used efficiently in CAD applications, being able to produce watertight and locally refined models. However, the use of the most general T-spline concept in IgA is limited by the risk of linear dependence of the resulting splines  \cite{buffa2010linear}. 


It is desirable in numerical simulations to use linear independent basis functions to ensure that the resulting mass and stiffness matrices have full rank and avoid the algorithmic complexity posed by singular matrices.

Analysis-Suitable T-splines, or AST-splines, were therefore introduced in \cite{da2012analysis}. As T-splines, AST-splines provide watertight models, obey the convex hull property, and moreover are linearly independent. 

There are many other definitions of B-splines over meshes with local refinements, such as PHT-splines \cite{deng2008polynomial}, PB-splines \cite{engleitner2017patchwork} and LR B-splines \cite{dokken2013polynomial}. A discusssion of the differences and similarities of  HB-splines, THB-splines, T-splines, AST-splines and LR B-splines can be found in \cite{dokken2018trivariate}.

\subsection{LR B-splines and MS B-splines}
In this paper we look at Locally Refined B-splines, or LR B-splines, introduced in \cite{dokken2013polynomial}. The idea is to  extend the knot insertion refinement of univariate B-splines to insertion of local line segments in tensor meshes. The process starts by considering the tensor product B-spline space over a coarse tensor mesh. Then, when a new inserted local line segment divides the support of one or more LR B-splines in two parts, we perform knot insertion to split such B-splines into two (or more) new ones. The resulting LR B-splines are a collection of scaled B-splines forming a non-negative partition of unity on the final mesh. 



The LR B-splines are a subset of the Minimal Support B-splines, or MS B-splines. As one can guess from their name, MS B-splines are the tensor product B-splines with minimal support, i.e., without superfluous line segments crossing their support, identifiable on the LR-mesh. The main difference between LR and MS B-splines is that the former ones are defined algorithmically, while the latter are defined by the topology of the mesh.
\subsection{Content of the paper}
The freedom in the refinement process can result in undesirable collections of LR B-splines. Namely, the refinement may create linear dependence relations among the B-splines. Assumptions on the refinement process have to be established in order to ensure linear independence. We start such analysis by looking at conditions on the mesh necessary for linear dependence. 
We say that functions $\phi_1, \ldots, \phi_n:\RR^d\to \RR$ are \textbf{actively linearly dependent} on $\RR^d$ if there exist $\alpha_i \in \RR,$ $\alpha_i \neq 0$ for all $i =1, \ldots, n$, such that $$
\sum_{i=1}^n \alpha_i\phi_i(\pmb{x}) = 0,\quad \forall~\pmb{x} \in \RR^d. 
$$
Note that we consequentially look at the minimal set of linearly dependent functions by forbidding zero coefficients in the linear combination.

In this work we show that:\begin{itemize}
\item For any bidegree $\pmb{p}$, the minimal number of active MS B-splines in a linear dependence relation on an LR-mesh is six, while for LR B-splines it is eight.
\item These numbers are sharp for any bidegree $\pmb{p} = (p_1,p_2)$ with $p_k \geq 1$ for the MS B-splines and $p_k\geq 2$ for LR B-splines for some $k \in \{1,2\}$.
\item If $(p_1,p_2) = (0,0)$, both MS and LR B-splines are linearly independent on any LR-mesh.
\end{itemize}
We look at the minimal configurations of linear dependence because we believe that any linear dependence relation is a refinement of one of these minimal cases. In other words, they are the roots for the linear dependence. Avoiding the minimal cases, the MS B-splines and LR B-splines are  linearly independent and form a basis. Furthermore, to get such lower bounds, we prove results that can be used to understand if the set of B-splines considered, on a given LR-mesh, is linearly independent or not. In particular, they can be used to improve the Peeling Algorithm \cite[ Algorithm 6.3]{dokken2013polynomial} to verify if the LR B-splines are linearly independent.
\subsection{Structure of the paper} 
 
In Section \ref{LRm} an introduction to the concepts of box-partitions, box-meshes and LR-meshes in 1D and 2D is provided.

In Section \ref{splspa} we define the spline spaces on such meshes and recall the dimension formula presented in \cite{pettersen2013dimension}. We then discuss conditions on the mesh for ensuring that the dimension formula depends only on the topology of the mesh and not on the position of the local line segments.

In Section \ref{Bspl} we recall univariate B-splines and tensor product B-splines, their basic properties and knot insertion.

In Section \ref{MSLRB} we define the MS B-splines and the LR B-splines on an LR-mesh. We show when these two sets are different.
 
In Section \ref{HIHP} we study the spanning properties of the LR and MS B-splines. In particular we  state necessary and sufficient conditions for spanning the full spline space. Knowing the dimension of the spline space, we can check linear dependencies just by counting the elements in the LR (or MS) B-spline set.

In Section \ref{LDMS}, we identify features needed for a linear dependence relation and we derive the minimal number of active MS B-splines needed in a linear dependence.

In Section \ref{LDLR}, we derive the minimal number of active LR B-splines in a linear dependence relation.

In Section \ref{peel}  we recall briefly the Peeling algorithm for checking linear independence and we show how to improve it using the results of Section \ref{LDMS}.

Finally, we summarize the main results and discuss future work in Section \ref{conc}.
\section{LR-meshes}\label{LRm}
The purpose of this section is to describe box-partitions in 1D and 2D and define bivariate LR-meshes. For our scope, and sake of simplicity, we decided to restrict general definitions, valid in any dimension, to the 2D case; we refer to  \cite{dokken2013polynomial} for the general case.
\subsection{Box-partitions in 1D}
 Consider a closed interval $\gamma = [a, b] \subseteq \RR$. We say that $\gamma$ is \textbf{trivial} if $a=b$ and \textbf{nontrivial} otherwise.
 \begin{dfn} Let $\gamma$ be nontrivial. A \textbf{box-partition} $\mathcal{E}$ of $[a, b]$ is a set of the form $$\mathcal{E} =\{ [\tau_i, \tau_{i+1}] \subseteq [a,b] : a = \tau_1 < \ldots < \tau_n = b\}.$$
The points $\tau_1, \ldots, \tau_n$ are called \textbf{knots} of the box-partition.\\
We can associate an integer $\mu(\tau_i) \geq 1$ to every knot and define a \textbf{knot vector on} $\gamma$ as the increasing sequence $\pmb{\tau}_\gamma^\mu= (\tau_1, \ldots, \tau_n)$ with $\mu(\tau_i)$ as assigned multiplicities to the knots.
\end{dfn}
 In some definitions it will be more convenient to write this knot vector on $\gamma$ as the non-decreasing sequence $\pmb{t}_\gamma = (t_1, \ldots, t_\ell)$, i.e, with $t_i \leq t_{i+1}$, where $\ell = \sum_{i=1}^n \mu(\tau_i)$ and $$
\underbrace{t_1 = \ldots = t_{\mu(\tau_1)}}_{=\, \tau_1} < \underbrace{t_{\mu(\tau_1)+1} = \ldots = t_{\mu(\tau_1) + \mu(\tau_2)}}_{=\,\tau_2} < \ldots
$$
Although it is an abuse of notation, the two sequences $\pmb{\tau}_\gamma^\mu$ and $\pmb{t}_\gamma$ encode the same knot vector on $\gamma$. 
We use bold Greek letters with the multiplicity function in superscript in the first way of expression for the knot vector and bold Latin letters for the second one. However, we will sometimes omit the subscript if the interval where these knots are defined is not relevant.
\subsection{LR-meshes}
\begin{dfn}
Let $\Omega \subseteq \RR^2$  be an axis-aligned rectangle. A \textbf{box-partition} of $\Omega$ is a finite collection $\mathcal{E}$ of (axis-aligned) rectangles, or \textbf{elements}, in $\Omega$ such that:
\begin{enumerate}
\item $\mathring{\beta}_1 \cap \mathring{\beta}_2 = \emptyset$ for any $\beta_1, \beta_2 \in \cE$, with $\beta_1 \neq \beta_2$.
\item $\cup_{\beta \in \cE} \beta= \Omega$.
\end{enumerate}\vspace{-.2cm}
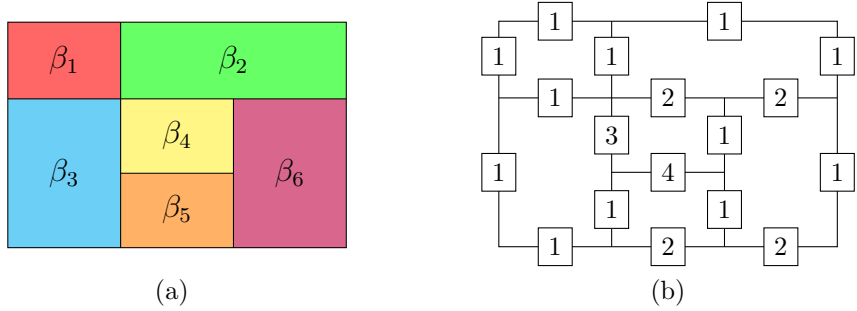
\begin{figure}[ht!]\centering
\subfloat[]{
\begin{tikzpicture}[scale=3]
\fill[red!60] (0,.66) -- (.5,.66) -- (.5, 1) -- (0,1) -- cycle;
\fill[cyan!50] (0,0) -- (.5,0) -- (.5, .66) -- (0,.66) -- cycle;
\fill[orange!60] (.5, 0) -- (1,0) -- (1,.33) -- (.5, .33) --cycle;
\fill[yellow!60] (.5, .33) -- (1,.33) -- (1,.66) -- (.5, .66) --cycle;
\fill[green!60] (.5, .66) -- (1.5,.66) -- (1.5,1) -- (.5, 1) --cycle;
\fill[purple!60] (1, 0) -- (1.5,0) -- (1.5,.66) -- (1, .66) --cycle;
\draw (0,0) -- (1.5,0) -- (1.5, 1) -- (0,1) --cycle;
\draw (0,.66) --(1.5, .66);
\draw (.5, .33) -- (1, .33);
\draw (.5, 0) -- (.5, 1);
\draw (1, 0) -- (1, .66);
\draw (.25, .83) node{$\beta_1$};
\draw (1, .83) node{$\beta_2$};
\draw (.25, .33) node{$\beta_3$};
\draw (.75, .17) node{$\beta_5$};
\draw (.75, .5) node{$\beta_4$};
\draw (1.25, .33) node{$\beta_6$};
\draw[white] (0,-.08);
\end{tikzpicture}}\qquad\qquad
\subfloat[]{
\begin{tikzpicture}[scale=3]
\tikzset{%
  block/.style    = {draw=black, fill=white, rectangle, minimum height = .05cm,
    minimum width = .05cm}
    }
\draw (0,0) -- (1.5,0) -- (1.5, 1) -- (0,1) --cycle;
\draw (0,.66) --(1.5, .66);
\draw (.5, .33) -- (1, .33);
\draw (.5, 0) -- (.5, 1);
\draw (1, 0) -- (1, .66);
\draw (0,.33) node[block]{\footnotesize{$1$}};
\draw (0,.845) node[block]{\footnotesize{$1$}};
\draw (.25,0) node[block]{\footnotesize{$1$}};
\draw (.25,.66) node[block]{\footnotesize{$1$}};
\draw (.25,1) node[block]{\footnotesize{$1$}};
\draw (.5,.165) node[block]{\footnotesize{$1$}};
\draw (.5, .495) node[block]{\footnotesize{$3$}};
\draw (.5,.845) node[block]{\footnotesize{$1$}};
\draw (.75,.66) node[block]{\footnotesize{$2$}};
\draw (.75,.33) node[block]{\footnotesize{$4$}};
\draw (.75,0) node[block]{\footnotesize{$2$}};
\draw (1.25,.66) node[block]{\footnotesize{$2$}};
\draw (1,1) node[block]{\footnotesize{$1$}};
\draw (1,.165) node[block]{\footnotesize{$1$}};
\draw (1,.495) node[block]{\footnotesize{$1$}};
\draw (1.5,.33) node[block]{\footnotesize{$1$}};
\draw (1.5,.845) node[block]{\footnotesize{$1$}};
\draw (1.25,0) node[block]{\footnotesize{$2$}};
\end{tikzpicture}
}
\caption{(a) Example of box-partition of a 2D element $\Omega$. (b) Example of $\mu$-extended box-mesh. Each meshline $\gamma$ has its own multiplicity $\mu(\gamma)$.}\label{fig2}
\end{figure}
We call \textbf{vertices of} $\mathcal{E}$ the set of vertices of its elements. A \textbf{meshline} in $\mathcal{E}$ is a segment connecting two and only two vertices of $\mathcal{E}$. The collection $\cM$ of meshlines is called \textbf{box-mesh} on $\Omega$. 
We can associate an integer $\mu(\gamma)\geq 1$ to every meshline $\gamma \in \cM$, called \textbf{multiplicity} of $\gamma$.  Then the pair $(\cM, \mu)$ is called $\pmb{\mu}$\textbf{-extended box-mesh}. See Figure \ref{fig2} for an example.
\end{dfn}
Note that each interior vertex in a box-mesh belongs either to 4 or 3 rectangles of the box-partition, corresponding to a \textbf{cross-vertex} ($+$) or a \textbf{T-vertex} ( \rotatebox[origin=c]{-180}{$\perp$}), respectively.

It is convenient for the rest of the paper to distinguish between \textbf{vertical and horizontal meshlines}, also called \textbf{1-meshlines and 2-meshlines} respectively.
\begin{dfn}
Given a box partition $\mathcal{E}$, corresponding to a box-mesh $\cM$, and an axis-aligned segment $\gamma$, we say that $\gamma$ \textbf{traverses} $\beta \in \mathcal{E}$ if $\beta\backslash \gamma$ is not connected, i.e., the interior of $\beta$ is divided into two pieces by $\gamma$. If $\gamma$ traverses $\beta$ then $\gamma \notin \cM$ and $\beta \backslash \gamma$ has two connected components $\beta_1, \beta_2$. We define $X_{\beta, \gamma} := \{\bar{\beta}_1, \bar{\beta}_2\}$ where $\bar{\beta}_j$ is the closure of $\beta_j$, for $j =1,2$.

A \textbf{split} is a finite union of contiguous and colinear segments $\gamma = \cup_i \gamma_i$   such that for every $i$ either $\gamma_i$ traverses a rectangle in $\mathcal{E}$ or it is a meshline already in $\cM$. Vertical splits and horizontal splits will be sometimes called \textbf{1-splits} and \textbf{2-splits} respectively.

A \textbf{maximal split $\gamma$ of }$\cM$ is a union of meshlines in $\cM$ such that any meshline in $\cM$ contiguous and colinear with $\gamma$ is contained in $\gamma$.
\end{dfn}
\begin{oss}No split can end in the middle of a rectangle of the box-partition $\mathcal{E}$.
\end{oss}
In order to define LR-meshes, we have to explain how the box-partition and associated box-mesh change when a new split is inserted.
\begin{dfn}
Let $\cE_\gamma$ be the set of all elements in $\cE$ whose interior is traversed by $\gamma$. We define $$
\cE + \gamma = \left(\cE\backslash\cE_\gamma\right) \cup \left(\bigcup_{\beta \in \cE_\gamma} X_{\beta, \gamma}\right),
$$
and we indicate the mesh related to $\cE+\gamma$ as $\cM + \gamma$. If we assign an integer $\mu_\gamma\geq 1$ to $\gamma$, 
the multiplicity of a meshline $\hat{\gamma} \in \cM+\gamma$ is \begin{itemize}
\item unchanged, $\mu(\hat{\gamma}) \to \mu(\hat{\gamma})$, if $\hat{\gamma} \notin \gamma$,
\item increased by $\mu_\gamma$, $\mu(\hat{\gamma}) \to \mu(\hat{\gamma})+\mu_\gamma$, if $\hat{\gamma} \in \gamma$ and it was meshline of $\cM$,
\item assigned to be $\mu_\gamma$, $\mu(\hat{\gamma}) = \mu_\gamma$, if $\hat{\gamma} \in \gamma$ and it traverses a rectangle of $\cE$.
\end{itemize}
\end{dfn}

We further need to define a box-mesh with constant splits and describe a tensor product mesh in our setting.
\begin{dfn}
A $\mu$-extended box-mesh $(\cM, \mu)$ has \textbf{constant (maximal) splits} if any maximal split $\gamma$ in it is made of meshlines of the same multiplicity, called \textbf{multiplicity of the split, $\mu(\gamma)$}.
\end{dfn}
\begin{dfn}
Given two increasing sequences $\pmb{\tau}_1$, $\pmb{\tau}_2$ with $\pmb{\tau}_k = (\tau_{k,1}, \ldots, \tau_{k,n_k})$, the associated \textbf{tensor product mesh} $\cM(\pmb{\tau}_1, \pmb{\tau}_2)$ is the box-mesh corresponding to the box-partition $\cE = \{[\tau_{1, i_1}, \tau_{1, i_1+1}]\times [\tau_{2,i_2}, \tau_{2, i_2+1}]: 1 \leq i_k \leq n_k -1, k=1, 2\}$ of $\Omega = [\tau_{1,1}, \tau_{1, n_1}]\times [\tau_{2,1}, \tau_{2,n_2}]$.
\end{dfn}
As for general box-partitions, we can associate a multiplicity $\mu(\gamma)$ to every meshline $\gamma \in \cM$ and define a $\mu$-extended tensor product mesh.


We can now define the LR-meshes.
\begin{dfn}
An \textbf{LR-mesh} is a $\mu$-extended box-mesh $\cM$ obtained through a sequence of split insertions:\begin{equation*}
\begin{array}{l}
\cM_1 \mbox{ is a $\mu$-extended tensor product mesh with constant splits,}\\\\\cM_{i+1} = \cM_i +\gamma_i \mbox{ has constant splits, for }i =1, \ldots, N-1
\end{array}
\end{equation*}
and $\cM = \cM_N$, for some $N$.   
\end{dfn}
Finally, in order to simplify the statements in what follows, we conclude this section by defining the length of a split and the knot vector on it. 
\begin{dfn} 
We say that a split $\gamma$ in an box-mesh $\cM$ \textbf{has length} $d$, or \textbf{is $d$ long}, if it intersects $d$ orthogonal meshlines on the mesh counting their multiplicities. 
\end{dfn}
\begin{dfn}
Given a box-mesh $\cM$ of box-partition $\cE$, for any vertex $\pmb{v}$ of $\cE$ we define \begin{equation*}
\begin{split}
\mu_1(\pmb{v}) &= \max\{\mu(\gamma): \pmb{v} \in \gamma \mbox{ and } \gamma \mbox{ 1-meshline of } \cM\}\\
\mu_2(\pmb{v}) &=\max\{\mu(\gamma): \pmb{v} \in \gamma \mbox{ and } \gamma \mbox{ 2-meshline of } \cM\}
\end{split}
\end{equation*}
$\mu_1(\pmb{v})$ is called \textbf{vertical multiplicity} and $\mu_2(\pmb{v})$ \textbf{horizontal multiplicity} of vertex $\pmb{v}$. See Figure \ref{5} for an example of computation of horizontal and vertical multiplicities.
\end{dfn}
\begin{figure}[h!]
\centering
\begin{tikzpicture}[scale=3]
\tikzset{%
  block/.style    = {draw=black, fill=white, rectangle, minimum height = .05cm,
    minimum width = .05cm},
    }
\draw (0,0) -- (1.5, 0) -- (1.5, 1) -- (0,1) -- cycle;
\draw (0,.5) -- (1.5,.5);
\draw (.5, 0) -- (.5, 1);
\draw (1, .5) -- (1, 1);
\draw (0,.25) node[block]{\footnotesize{$1$}};
\draw (0,.75) node[block]{\footnotesize{$1$}};
\draw (.25,0) node[block]{\footnotesize{$1$}};
\draw (.25,.5) node[block]{\footnotesize{$1$}};
\draw (.25,1) node[block]{\footnotesize{$1$}};
\draw (.5,.25) node[block]{\footnotesize{$1$}};
\draw (.5,.75) node[block]{\footnotesize{$1$}};
\draw (.75,.5) node[block]{\footnotesize{$2$}};
\draw (1.25,.5) node[block]{\footnotesize{$2$}};
\draw (.75,1) node[block]{\footnotesize{$1$}};
\draw (1.25,1) node[block]{\footnotesize{$1$}};
\draw (1,0) node[block]{\footnotesize{$1$}};
\draw (1,.75) node[block]{\footnotesize{$1$}};
\draw (1.5,.25) node[block]{\footnotesize{$1$}};
\draw (1.5,.75) node[block]{\footnotesize{$1$}};
\fill (.5, .5) node[below]{\footnotesize{$\quad~\pmb{v}_1$}} circle (.02);
\fill (1, .5) node[below]{\footnotesize{$\pmb{v}_2$}} circle (.02);
\end{tikzpicture}
\caption{Example of computation of vertical and horizontal multiplicities. The meshlines on the left and right hand-side of $\pmb{v}_1$ have multiplicity 1 and 2 respectively. So $\mu_2(\pmb{v}_1) = \max\{1, 2\} = 2$. The meshlines above and below $\pmb{v}_1$ have both multiplicity 1, so that $\mu_1(\pmb{v}_1) = 1.$ Concerning $\pmb{v}_2$, we have $\mu_2(\pmb{v}_2) = 2$, whereas $\mu_1(\pmb{v}_2) = \max\{1\} = 1$ since there is no meshline below $\pmb{v}_2$.}\label{5}
\end{figure}
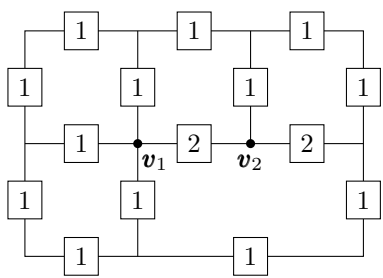
\begin{dfn}\label{knotvec} Given a $k$-split $\gamma$, for $k \in \{1,2\}$, of length $d$ in an box-mesh $\cM$, we can define the  \textbf{knot vector on} $\gamma$, $\pmb{\tau}_\gamma^{\mu_{3-k}}$, as the points in $\RR$ corresponding to the $(3-k)$th-components of the vertices  of $\cM$  where $\gamma$ intersects the orthogonal splits, with assigned multiplicity equal to the multiplicity $\mu_{3-k}$ of these vertices. 
\end{dfn}
\begin{ex}
Assume $\gamma$ is a 2-split of length $d$ and intersecting the vertices $\pmb{v}_1, \ldots, \pmb{v}_\ell$ with $\pmb{v}_i = (v_i, v_\gamma)$ for $i=1,\ldots,\ell$. Then the knot vector on $\gamma$ is $\pmb{\tau}_\gamma^{\mu_1} = (\tau_1, \ldots, \tau_\ell)$ with $\tau_i = v_i$ and knot multiplicities $\mu_1(\tau_i) = \mu_1(\pmb{v}_i)$ where $\sum_i \mu_1(\tau_i) = d$.  
\end{ex}
\section{Spline space over a box-mesh in 1D and 2D}\label{splspa}
In this section we define the spline spaces over box-meshes in 1D and 2D and present the dimension formulas of them. In 2D, the general formula \cite{pettersen2013dimension} presents terms depending on the paramenters of the meshlines. We recall sufficient conditions for nullifying such terms, making the formula dependent only on the topology of the mesh and, therefore, stable.
\subsection{Spline space in 1D}
\begin{dfn}
Given a degree $p$, we define $\Pi_p \subset \RR[t]$ to be the vector space spanned by the monomials $t^j$ such that $0\leq j\leq p$.

Given a knot vector $\pmb{\tau}_\gamma^\mu$ from a box-partition $\mathcal{E}$ of a nontrivial interval $\gamma =[a, b]$, with $\mu(\tau_i) \leq p+1$ for every knot $\tau_i \in \pmb{\tau}_\gamma^\mu$, 
the \textbf{univariate spline space of degree} $p$ on $\pmb{\tau}_\gamma^\mu$, denoted by $\SSS_p(\pmb{\tau}_\gamma^\mu)$ (or $\SSS_p(\pmb{t}_\gamma)$), is the set of all functions $f: \RR \to \RR$ such that \begin{itemize}
\item $f$ is zero outside $\gamma$,
\item for every element $\gamma_i= [\tau_i, \tau_{i+1}]\in \mathcal{E}$, $f$ restricted to $[\tau_i, \tau_{i+1})$ is a polynomial in $\Pi_p$,
\item for each knot $\tau_i\in \pmb{\tau}_\gamma^\mu$, $f$ is $C^{p-\mu(\tau_i)}$-continuous across $\tau_i$. 
\end{itemize}
\end{dfn}
The following result is the dimension formula for the spline space in the 1D case. The general statement for any dimension and its proof can be found in \cite{pettersen2013dimension}. 
\begin{thm}\label{dim1d}
Let $\gamma \subset \RR$ be a nontrivial interval. Given a degree $p$ and a knot vector $\pmb{\tau}_\gamma^\mu= (\tau_1, \ldots, \tau_n)$ on $\gamma$, then  
\begin{equation}\label{dim1d2}
\dim\,\SSS_p(\pmb{\tau}_\gamma^\mu) = \max\left\{\sum_{i=1}^n \mu(\tau_i) - (p+1),0\right\}.
\end{equation}
\end{thm}
Therefore, 
if $\pmb{t}_\gamma$ has cardinality $p+r+1$  for some $r\geq 1$, then $\dim\SSS_p(\pmb{t}_\gamma) = r.$

There are many possible bases for $\SSS_p(\pmb{t}_\gamma)$. One possibility comes out from a classical result in spline theory, called Curry-Schoenberg Theorem. It ensures that the so called B-spline functions of degree $p$, defined on the knot vector $\pmb{t}_\gamma$, can be used as a possible basis: $$
\SSS_p(\pmb{t}_\gamma) = \mbox{span}\,\{B[\pmb{t}_{\gamma}^i]\}_{i=1}^{r}\quad \mbox{with }\pmb{t}_{\gamma}^i = (t_i, \ldots, t_{i+p+1})\subseteq \pmb{t}_\gamma.
$$ 
For a brief introduction to B-splines we refer to Section \ref{Bspl}.
\begin{oss}
Definition \ref{knotvec} tells us that there is a knot vector on every split of a box-mesh. We can therefore define a univariate spline space on it where appropriate.
\end{oss}
\subsection{Spline space in 2D}
\begin{dfn}
Given a bidegree $\pmb{p} = (p_1, p_2)$, we define $\Pi_{\pmb{p}} \subset \RR[x, y]$ to be the vector space spanned by the monomials $x^{i_1}, y^{i_2}$ such that $0\leq i_k \leq p_k$ for $k = 1, 2$.

Given a $\mu$-extended box mesh $\cM$ from a box partition $\cE$ of $\Omega = [a_1, b_1]\times [a_2, b_2] \subset \RR^2$, suppose $\mu(\gamma) \leq p_k +1$ for every $k$-meshline $\gamma \in \cM$. For any element $\beta \in \cE$, $
\beta = J_1\times J_2 \mbox{ with }J_k = [a_{\beta, k}, b_{\beta, k}]
$,
we set $$
\tilde{\beta} = \tilde{J}_1 \times \tilde{J}_2\mbox{ with }\tilde{J}_k = \left\{\begin{array}{ll} [a_{\beta, k}, b_{\beta, k}) & \mbox{if } b_{\beta, k} < b_k\\
\\

[a_{\beta, k}, b_{\beta, k}] & \mbox{if } b_{\beta, k} = b_k
\end{array}\right.
$$
The \textbf{spline space of degree p on }$\pmb{\cM}$, denoted by $\SSS_{\pmb{p}}(\cM)$ is the set of all function $f:\RR^2\to \RR$ such that \begin{enumerate}
\item $f$ is zero outside $\Omega$,
\item for each element $\beta \in \cE$, the restriction $f$ to $\tilde{\beta}$ is a bivariate polynomial function in $\Pi_{\pmb{p}}$,
\item for each $k$-meshline $\gamma \in \cM$, $f \in C^{p_k-\mu(\gamma)}(\gamma)$.
\end{enumerate}
\end{dfn}

The general dimension formula of the spline space over box-partitions is presented in \cite{pettersen2013dimension} and has terms depending on the parametrization of the splits. This makes the dimension of the spline space unstable \cite{li2011instability}. However, if we consider the spline space over an LR-mesh built so that \begin{itemize}
\item[\fbox{\textbf{LR-rule 1}}] the starting tensor product mesh $\cM_1$ has at least $p_1+2$ vertical splits and $p_2+2$ horizontal splits counting their multiplicities,
\item[\fbox{\textbf{LR-rule 2}}] For $k\in \{1,2\}$, any maximal $k$-split at any step in the construction of the LR-mesh has length at least $p_{3-k}+2$.
\end{itemize}
then, called $\mathcal{F}_0(\cM)$ the set of vertices of $\cE$, $V$ the set of 1-meshlines, $H$ the set of 2-meshlines and $|\mathcal{E}|$ the cardinality of $\mathcal{E}$, we have \begin{equation}\begin{split}\label{dimS}
\dim \SSS_{\pmb{p}}(\cM) &= \sum_{\pmb{v}\in \cF_0(\cM)} [(p_1-\mu_1(\pmb{v})+1)(p_2-\mu_2(\pmb{v})+1)]\\ &\quad- (p_2+1)\sum_{\beta \in V} [(p_1-\mu(\beta) +1)] - (p_1+1)\sum_{\beta \in H} [(p_2-\mu(\beta) +1)]\\ &\quad+ |\mathcal{E}|(p_1+1)(p_2+1),
\end{split}
\end{equation}
and it depends only on the topology of the mesh and not on its parametrization. In this paper we will always assume the LR-rules for constructing LR-meshes.
\begin{oss}
In the LR-mesh building process, 
any extension of an older split is allowed being LR-rule 2 satisfied on the new mesh. 
\end{oss}
From \eqref{dimS}, it is possible to prove the dimension increasing formula \cite[Theorem 5.5]{dokken2013polynomial}. Knowing $\dim\mathbb{S}_{\pmb{p}}(\cM)$, through this formula, one can easily compute the dimension of the spline space on a refined mesh $\cM + \gamma$. First, we need to introduce the concept of expanded split.
\begin{dfn}
Given a degree $p_{3-k}$, an LR-mesh $\cM$, and a $k$-split $\gamma$, we define the \textbf{expanded $k$-split} corresponding to $\gamma$, as the $k$-split $\bar{\gamma} \in \cM + \gamma$ which is equal to $\gamma$ if $\gamma$ is not an extension of a split in $\cM$; while if  $\gamma$ is an extension of a split $\tilde{\gamma}$ in $\cM$, $\bar{\gamma}$ is $\gamma$ plus $p_{3-k}+1$ meshlines of $\tilde{\gamma}$ contiguous to $\gamma$. 
\end{dfn}
In particular, if $\gamma$ is an extesion of two $k$-splits $\gamma_1, \gamma_2$ in $\cM$, i.e., $\gamma$ is the link between $\gamma_1, \gamma_2$, we have to consider $p_{3-k}+1$ meshlines both in $\gamma_1$ and $\gamma_2$. 

We can now given the dimension increasing formula.
\begin{thm}\label{diminc}Given an LR-mesh $\cM$ and a $k$-split $\gamma$ such that the corresponding expanded $k$-split $\bar{\gamma}$ has length $p_{3-k}+r+1$ on $\cM + \gamma$ with $r \geq 1$, then
$$ 
\dim \SSS_{\pmb{p}}(\cM+\gamma) = \dim \SSS_{\pmb{p}}(\cM) + r = \dim \SSS_{\pmb{p}}(\cM) + \dim \SSS_{p_{3-k}} (\pmb{t}_{\bar{\gamma}}).
$$
\end{thm}
\section{Univariate B-splines and tensor product B-splines in 2D}\label{Bspl}
Here we recall the definition of B-splines and their main properties. In particular, we state the knot insertion algorithm, which is used for the definition of LR B-splines.
For a complete overview on B-splines we refer to \cite{de1978practical} and \cite{schumaker2007spline}.
\subsection{Univariate B-splines}\label{unispl}
\begin{dfn}
For a non-decreasing sequence $\pmb{t} = (t_1, t_2, \ldots, t_{p+2})$ we define a \textbf{B-spline} $B[\pmb{t}]:\RR \to \RR$ \textbf{of degree} $p\geq 0$ recursively by 
\begin{equation}
B[\pmb{t}](t) = \frac{t-t_1}{t_{p+1}-t_1} B[t_1, \ldots, t_{p+1}](t) + \frac{t_{p+2}-t}{t_{p+2}-t_2}B[t_2, \ldots, t_{p+2}](t),
\end{equation}
where each time a fraction with zero denominator appears, it is taken as zero. The initial B-splines of degree 0 on $\pmb{t}$ are defined as \begin{equation}
B[t_i, t_{i+1}](t) := \left\{\begin{array}{ll} 1&\mbox{if }t_i \leq t<t_{i+1};\\\\0 &\mbox{otherwise};\end{array}\right. \quad \mbox{for }i=1, \ldots, p+1.
\end{equation}
The sequence $\pmb{t}$ is called \textbf{knot vector} of $B[\pmb{t}]$ and $t_j$ are called \textbf{knots}. A knot $t_j$ has multiplicity $\mu(t_j)$ if it appears $\mu(t_j)$ times in $\pmb{t}$. 
\end{dfn}
\begin{prop}[Properties]
Given a degree $p\geq0$ and a knot vector $\pmb{t} = (t_1, \ldots, t_{p+2})$, \begin{itemize}
\item $\mbox{supp}\,B[\pmb{t}] = [t_1, t_{p+2}]$,
\item $B[\pmb{t}]$ restricted to every nontrivial half-open element $[t_i, t_{i+1})$ is in $\Pi_p$, 
\item $B[\pmb{t}]$ is $C^{p-\mu(t_j)}$-continuous at any knot $t_j$ of multiplicity $\mu(t_j)$.
\end{itemize}
\end{prop}
\begin{thm}[knot insertion]
Given a degree $p$ and a corresponding knot vector $\pmb{t}=(t_1, \ldots, t_{p+2})$, suppose we insert a knot $\hat{t} \in (t_1, t_{p+2})$. We obtain two knot vectors $\pmb{t}_1$ and $\pmb{t}_2$, considering the first and the last $p+2$ knots respectively in $(t_1, \ldots, \hat{t}, \ldots, t_{p+2})$. Then there exist $\alpha_1, \alpha_2 \in [0,1]$ such that \begin{equation}
B[\pmb{t}] = \alpha_1B[\pmb{t}_1] + \alpha_2B[\pmb{t}_2],
\end{equation}
\end{thm}
Now we recall the definition of bivariate B-splines and their basic properties inherited by the univariate B-splines. 
\begin{dfn}
Consider a bidegree $\pmb{p}=(p_1, p_2)$. Let $\pmb{x} = (x_1, \ldots, x_{p_1+2})$ and $\pmb{y} = (y_1, \ldots, y_{p_2+2})$ be nondecreasing knot vectors. We define the \textbf{tensor product B-spline} $B[\pmb{x}, \pmb{y}]: \RR^2 \to \RR$ by \begin{equation}
B[\pmb{x}, \pmb{y}](x, y) := B[\pmb{x}](x)B[\pmb{y}](y),
\end{equation}
where $B[\pmb{x}]$ and $B[\pmb{y}]$ are the univariate B-splines defined on $\pmb{x}$ and $\pmb{y}$ respectively.
\end{dfn}
The pair $\pmb{x}, \pmb{y}$ identifies a tensor product mesh in the element $[x_1, x_{p_1+2}]\times [y_1, y_{p_2+2}]$, $\cM(\pmb{x}, \pmb{y})$. Hence, a knot in the $x$-direction $x_i$ of multiplicity $\mu_1(x_i)$ correspond to a 1-split $\gamma =\{x_i\}\times [y_1, y_{p_2+2}]$ of multiplicity $\mu(\gamma) = \mu_1(x_i)$ in $\cM(\pmb{x}, \pmb{y})$. The same holds for the knots $y_j$ in the $y$-direction. Such splits will be called \textbf{splits of }$B[\pmb{x}, \pmb{y}]$.

The properties of univariate B-splines are conserved by the tensor product B-splines:\begin{itemize}
\item $\mbox{supp}\,B[\pmb{x}, \pmb{y}] = [x_1, x_{p_1+2}]\times [y_1, y_{p_2+2}]$.
\item $B[\pmb{x}, \pmb{y}]$ is a piecewise bivariate polynomial of bidegree $\pmb{p}$.
\item $B[\pmb{x}, \pmb{y}]$ is $C^{p_k-\mu(\gamma)}$-continuous across each $k$-split $\gamma$ of $B[\pmb{x}, \pmb{y}]$.
\end{itemize}
As in the univariate case, after the insertion of a knot $\hat{x}$ in $\pmb{x}$, we can write $B[\pmb{x}, \pmb{y}]$ in terms of two B-splines defined on the two new pairs of knot vectors \begin{equation}
B[\pmb{x}, \pmb{y}]=\alpha_1B[\pmb{x}_1, \pmb{y}] + \alpha_2B[\pmb{x}_2, \pmb{y}].
\end{equation}
The same holds when inserting a knot $\hat{y}$ in $\pmb{y}$.

Given a bidegree $\pmb{p} = (p_1, p_2)$ and a tensor product mesh $\cM(\pmb{x},\pmb{y})$ with $\pmb{x}, \pmb{y}$ of respectively $p_1 + r_1+1$ and $p_2 + r_2+1$ elements with $r_1, r_2 \geq 1$, we can apply the Curry-Schoenberg theorem on each univariate knot vector and state that 
$$
\SSS_{\pmb{p}}(\cM(\pmb{x}_{\gamma_1}, \pmb{y}_{\gamma_2})) = \mbox{span}\,\{B[\pmb{x}_{\gamma_1}^i, \pmb{y}_{\gamma_2}^j]\}\quad\mbox{with }i=1, \ldots, r_1\mbox{ and }j=1, \ldots, r_2,
$$
where $\pmb{x}_{\gamma_1}^i = (x_i, \ldots, x_{i+p_1+1})\subseteq \pmb{x}_{\gamma_1}$ and $\pmb{y}_{\gamma_2}^j = (y_j, \ldots, y_{j+p_2+1})\subseteq \pmb{y}_{\gamma_2}$.

\section{Minimal Support B-splines and LR B-splines}\label{MSLRB}
In this section we define first the Minimal Support B-splines, or MS B-splines, and then the Locally Refined B-splines, or LR B-splines. As we will see the LR B-splines are created algorithmically, splitting, after inserting a split in the mesh, the B-splines traversed by the split through the knot insertion procedure. The main difference with the MS B-splines is that the latter can be created from scratch on the mesh and not only by using knot insertion. 
\begin{dfn}
Given a bivariate B-spline $B[\pmb{x}, \pmb{y}]$ and a split $\gamma$, we say that $\gamma$ \textbf{traverses} $B[\pmb{x}, \pmb{y}]$ if the interior of $supp\,B[\pmb{x},\pmb{y}]$ is divided into two parts by $\gamma$.
\end{dfn} 
\begin{dfn}
Given a pair of knot vectors $\pmb{x}, \pmb{y}$ and an LR-mesh $\cM$ on a rectangle $\Omega$, the tensor product B-spline $B[\pmb{x}, \pmb{y}]$ has \textbf{support} in $\cM$ if\begin{itemize}
\item $supp B[\pmb{x}, \pmb{y}]\subseteq \Omega$, 
\item every split of $B[\pmb{x}, \pmb{y}]$ of multiplicity $m$ is contained in a split $\gamma$ of $\cM$ of multiplicity $\mu(\gamma) \geq m$.
\end{itemize}
$B[\pmb{x}, \pmb{y}]$ has \textbf{minimal support} in $\cM$ if in addition\begin{itemize}
\item a split $\gamma$ of $\cM$ of multiplicity $\mu(\gamma)$ traverses $B[\pmb{x}, \pmb{y}]$ $\sse$ a split of $B[\pmb{x}, \pmb{y}]$ of multiplicity $\mu(\gamma)$ is contained in $\gamma$ 
(See Figure \ref{6} for an example).
\end{itemize}
Given an LR-mesh, we write $\mathcal{B}^{\mathcal{MS}}(\cM)$ for the set of all minimal support B-splines defined on it. 
\end{dfn}
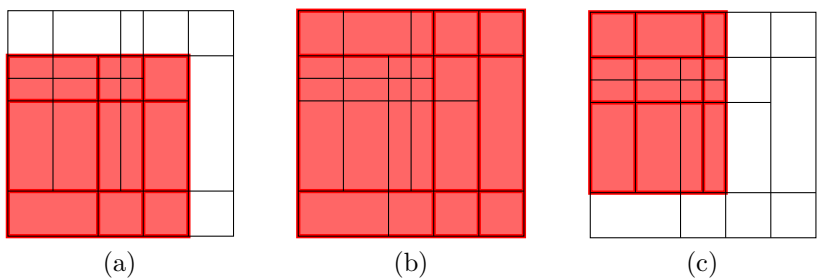
\begin{figure}[ht!]
    \centering
\subfloat[]{\centering
\begin{tikzpicture}[scale = 3]
\filldraw[draw = red,ultra thick, fill =red!60] (0,0) -- (.8,0) -- (.8,.8) -- (0,.8) -- cycle;
\draw[red,ultra thick] (0, .2) -- (.8, .2);
\draw[red,ultra thick] (0,.6) -- (.8, .6);
\draw[red,ultra thick] (.4,0) -- (.4, .8);
\draw[red,ultra thick] (.6, 0) -- (.6, .8);
\draw (0,0) -- (1,0) -- (1,1) -- (0,1) -- cycle;
\draw (0,.2) --(1, .2);
\draw (0,.8) -- (1,.8);
\draw (.6, 0) -- (.6, 1);
\draw (.8, 0) -- (.8, 1);
\draw (0,.6) -- (.8, .6);
\draw (.2, .2) -- (.2, 1);
\draw (0, .7) -- (.6, .7);
\draw (.4, 0) -- (.4, .8);
\draw (.5, .2) -- (.5, 1);
\end{tikzpicture}}\qquad
\subfloat[]{\centering
\begin{tikzpicture}[scale = 3]
\filldraw[draw = red,ultra thick, fill =red!60] (0,0) -- (1,0) -- (1,1) -- (0,1) -- cycle;
\draw[red,ultra thick] (0, .2) -- (1, .2);
\draw[red,ultra thick] (0,.8) -- (1, .8);
\draw[red,ultra thick] (.6,0) -- (.6, 1);
\draw[red,ultra thick] (.8, 0) -- (.8, 1);
\draw (0,0) -- (1,0) -- (1,1) -- (0,1) -- cycle;
\draw (0,.2) --(1, .2);
\draw (0,.8) -- (1,.8);
\draw (.6, 0) -- (.6, 1);
\draw (.8, 0) -- (.8, 1);
\draw (0,.6) -- (.8, .6);
\draw (.2, .2) -- (.2, 1);
\draw (0, .7) -- (.6, .7);
\draw (.4, 0) -- (.4, .8);
\draw (.5, .2) -- (.5, 1);
\end{tikzpicture}}\qquad
\subfloat[]{\centering
\begin{tikzpicture}[scale = 3]
\filldraw[draw = red,ultra thick, fill =red!60] (0,.2) -- (.6,.2) -- (.6,1) -- (0,1) -- cycle;
\draw[red,ultra thick] (0, .6) -- (.6, .6);
\draw[red,ultra thick] (0,.8) -- (.6, .8);
\draw[red,ultra thick] (.2,.2) -- (.2, 1);
\draw[red,ultra thick] (.5, .2) -- (.5, 1);
\draw (0,0) -- (1,0) -- (1,1) -- (0,1) -- cycle;
\draw (0,.2) --(1, .2);
\draw (0,.8) -- (1,.8);
\draw (.6, 0) -- (.6, 1);
\draw (.8, 0) -- (.8, 1);
\draw (0,.6) -- (.8, .6);
\draw (.2, .2) -- (.2, 1);
\draw (0, .7) -- (.6, .7);
\draw (.4, 0) -- (.4, .8);
\draw (.5, .2) -- (.5, 1);
\end{tikzpicture}}
\caption{Support of B-splines of bidegree (2,2) on an LR-mesh of multiplicity 1. (a) and (b) have minimal support whereas (c) has not minimal support: there is an extra horizontal split traversing its support.}\label{6}
\end{figure}
\begin{oss}
Given an LR-mesh $\cM$, in general $span\,\mathcal{B}^{\mathcal{MS}}(\cM) \subseteq \SSS_{\pmb{p}}(\cM)$. If the LR-mesh is a tensor product mesh then the MS B-splines are nothing more than standard bivariate B-splines, and the Curry-Schoenberg Theorem ensures that $\mbox{span}\,\mathcal{B}^{\mathcal{MS}}(\cM) = \SSS_{\pmb{p}}(\cM)$ and the elements of $\mathcal{B}^{\mathcal{MS}}(\cM)$ are linearly independent.
\end{oss}
We now define the set of LR B-splines of degree $\pmb{p}$ on $\cM$, $\mathcal{B}^{\mathcal{LR}}(\cM)$. Recall that any LR-mesh is obtained through a sequence $\cM_{i+1} = \cM_i + \gamma_i$, with $\cM_1$ a tensor product mesh and $\gamma_i$ splits so that  $\cM_{i+1}$ has constant splits. 
%
Whenever a MS B-spline $B$ on $\cM_i$ is traversed by a $k$-split in $\cM_{i+1}$, such split corresponds to a new knot in the $k$th knot vector of $B$. Such MS B-spline could then be refined using the knot insertion procedure. 
\begin{dfn}
Given an LR-mesh $\cM = \cM_N$, the set $\mathcal{B}^{\mathcal{LR}}(\cM)$ is provided algorithmically as follows:\\
\begin{algorithm}[H]
\vspace{.2cm}
Let $\mathcal{B}_1 = \{B[\pmb{x}^j, \pmb{y}^j]\}_j$ be the B-splines of degree $\pmb{p}$ on the tensor product mesh $\cM_1$\;\vspace{.15cm}
\For{every intermediate mesh $\cM_{i+1} = \cM_i + \gamma_i$ with $i=1,\ldots, N-1$}{
$\mathcal{B}_{i+1} = \mathcal{B}_i$\;\vspace{.1cm}
\While{$\exists B[\pmb{x}^j, \pmb{y}^j]\in \mathcal{B}_{i+1}$ with no minimal support on $\cM_{i+1}$}{\vspace{.15cm}
Knot insertion $\Rightarrow \exists B[\pmb{x}_{1}^j, \pmb{y}_{1}^j], B[\pmb{x}_{2}^j, \pmb{y}_{2}^j] : B[\pmb{x}^j, \pmb{y}^j] = \alpha_1B[\pmb{x}_1^{j}, \pmb{y}_1^{j}]+ \alpha_2 B[\pmb{x}_2^{j}, \pmb{y}_2^{j}]$\; \vspace{.15cm}
$\mathcal{B}_{i+1} = (\mathcal{B}_{i+1} \backslash \{B[\pmb{x}^j, \pmb{y}^j]\})\cup\{B[\pmb{x}_1^{j}, \pmb{y}_1^{j}], B[\pmb{x}_2^{j}, \pmb{y}_2^{j}]\}$\; 
}\vspace{.15cm}
}
$\mathcal{B}^{\mathcal{LR}}(\cM) = \mathcal{B}_N$.
\end{algorithm}
\end{dfn}
\begin{oss}
In the while cycle, it could be that we split again $B[\pmb{x}_1^{j}, \pmb{y}_1^{j}]$ or $B[\pmb{x}_2^{j}, \pmb{y}_2^{j}]$ and replace them with minimal support B-splines. 
\end{oss}
\begin{oss}
For any LR-mesh $\cM$, $\mbox{span}\,\mathcal{B}^{\mathcal{LR}}(\cM) \subseteq \mbox{span}\,\mathcal{B}^{\mathcal{MS}}(\cM) \subseteq \SSS_{\pmb{p}}(\cM)$. If $\cM$ is a tensor product mesh then the sets $\mathcal{B}^{\mathcal{LR}}(\cM) = \mathcal{B}^{\mathcal{MS}}(\cM)$ and are nothing more than the standard bivariate B-splines. Again we have, by the Curry-Schoenberg Theorem, $\mbox{span}\,\mathcal{B}^{\mathcal{LR}}(\cM) =\mbox{span}\,\mathcal{B}^{\mathcal{MS}}(\cM) = \SSS_{\pmb{p}}(\cM)$. However, there are other cases where this equality holds; we will see them in the next section.
\end{oss}
\begin{ex}[$\mathcal{B}^{\mathcal{LR}}(\cM)\neq \mathcal{B}^{\mathcal{MS}}(\cM)$]$\,$
\begin{figure}[ht!]
\hspace{.25cm}
    \centering
\subfloat[]{\hspace{-1.75cm}
\begin{tikzpicture}[scale=2.5]
\draw (0,0) -- (1.25,0) -- (1.25,1) -- (0,1) -- cycle;
\draw (0,.55) -- (.75, .55);
\draw (.5, .44) -- (1.25, .44);
\draw (.583, 0) -- (.583, .66);
\draw (.666, .33) -- (.666, 1);
\draw (0,.33) -- (1.25, .33);
\draw (0,.66) -- (1.25, .66);
\draw (.25, 0) -- (.25, 1);
\draw (.5, 0) -- (.5, 1);
\draw (.75, 0) -- (.75, 1);
\draw (1, 0) -- (1, 1);
\draw[white] (-.6,0);
\draw[white] (0,-.08);
\end{tikzpicture}}
\subfloat[]{
\begin{minipage}{8.5cm} \vspace{-2.85cm} 
\begin{tabular}{c}
\begin{tikzpicture}[scale=1.5]
\filldraw[draw=red, ultra thick, fill=red!60] (0,.33) -- (.66, .33) -- (.66,1) -- (0,1) -- cycle;
\draw[red, ultra thick] (0,.55) -- (.66, .55);
\draw[red, ultra thick] (0,.66) -- (.66, .66);
\draw[red, ultra thick] (.25,.33) -- (.25, 1);
\draw[red, ultra thick] (.5,.33) -- (.5, 1);
\draw (0,0) -- (1.25,0) -- (1.25,1) -- (0,1) -- cycle;
\draw (0,.55) -- (.75, .55);
\draw (.5, .44) -- (1.25, .44);
\draw (.583, 0) -- (.583, .66);
\draw (.666, .33) -- (.666, 1);
\draw (0,.33) -- (1.25, .33);
\draw (0,.66) -- (1.25, .66);
\draw (.25, 0) -- (.25, 1);
\draw (.5, 0) -- (.5, 1);
\draw (.75, 0) -- (.75, 1);
\draw (1, 0) -- (1, 1);
\end{tikzpicture},
\begin{tikzpicture}[scale=1.5]
\filldraw[draw=red, ultra thick, fill=red!60] (.25,.33) -- (.75, .33) -- (.75,1) -- (.25,1) -- cycle;
\draw[red, ultra thick] (.25,.55) -- (.75, .55);
\draw[red, ultra thick] (.25,.66) -- (.75, .66);
\draw[red, ultra thick] (.5,.33) -- (.5, 1);
\draw[red, ultra thick] (.66,.33) -- (.66, 1);
\draw (0,0) -- (1.25,0) -- (1.25,1) -- (0,1) -- cycle;
\draw (0,.55) -- (.75, .55);
\draw (.5, .44) -- (1.25, .44);
\draw (.583, 0) -- (.583, .66);
\draw (.666, .33) -- (.666, 1);
\draw (0,.33) -- (1.25, .33);
\draw (0,.66) -- (1.25, .66);
\draw (.25, 0) -- (.25, 1);
\draw (.5, 0) -- (.5, 1);
\draw (.75, 0) -- (.75, 1);
\draw (1, 0) -- (1, 1);
\end{tikzpicture},
\begin{tikzpicture}[scale=1.5]
\filldraw[draw=red, ultra thick, fill=red!60] (.5,.33) -- (1, .33) -- (1,1) -- (.5,1) -- cycle;
\draw[red, ultra thick] (.5,.44) -- (1, .44);
\draw[red, ultra thick] (.5,.66) -- (1, .66);
\draw[red, ultra thick] (.666,.33) -- (.666, 1);
\draw[red, ultra thick] (.75,.33) -- (.75, 1);
\draw (0,0) -- (1.25,0) -- (1.25,1) -- (0,1) -- cycle;
\draw (0,.55) -- (.75, .55);
\draw (.5, .44) -- (1.25, .44);
\draw (.583, 0) -- (.583, .66);
\draw (.666, .33) -- (.666, 1);
\draw (0,.33) -- (1.25, .33);
\draw (0,.66) -- (1.25, .66);
\draw (.25, 0) -- (.25, 1);
\draw (.5, 0) -- (.5, 1);
\draw (.75, 0) -- (.75, 1);
\draw (1, 0) -- (1, 1);
\end{tikzpicture},
\begin{tikzpicture}[scale=1.5]
\filldraw[draw=red, ultra thick, fill=red!60] (.666,.33) -- (1.25, .33) -- (1.25,1) -- (.666,1) -- cycle;
\draw[red, ultra thick] (.666,.44) -- (1.25, .44);
\draw[red, ultra thick] (.666,.66) -- (1.25, .66);
\draw[red, ultra thick] (.75,.33) -- (.75, 1);
\draw[red, ultra thick] (1,.33) -- (1, 1);
\draw (0,0) -- (1.25,0) -- (1.25,1) -- (0,1) -- cycle;
\draw (0,.55) -- (.75, .55);
\draw (.5, .44) -- (1.25, .44);
\draw (.583, 0) -- (.583, .66);
\draw (.666, .33) -- (.666, 1);
\draw (0,.33) -- (1.25, .33);
\draw (0,.66) -- (1.25, .66);
\draw (.25, 0) -- (.25, 1);
\draw (.5, 0) -- (.5, 1);
\draw (.75, 0) -- (.75, 1);
\draw (1, 0) -- (1, 1);
\end{tikzpicture},
\begin{tikzpicture}[scale=1.5]
\filldraw[draw=red, ultra thick, fill=red!60] (0,0) -- (.583, 0) -- (.583,.66) -- (0,.66) -- cycle;
\draw[red, ultra thick] (0,.33) -- (.583, .33);
\draw[red, ultra thick] (0,.55) -- (.583, .55);
\draw[red, ultra thick] (.25,0) -- (.25, .66);
\draw[red, ultra thick] (.5,0) -- (.5, .66);
\draw (0,0) -- (1.25,0) -- (1.25,1) -- (0,1) -- cycle;
\draw (0,.55) -- (.75, .55);
\draw (.5, .44) -- (1.25, .44);
\draw (.583, 0) -- (.583, .66);
\draw (.666, .33) -- (.666, 1);
\draw (0,.33) -- (1.25, .33);
\draw (0,.66) -- (1.25, .66);
\draw (.25, 0) -- (.25, 1);
\draw (.5, 0) -- (.5, 1);
\draw (.75, 0) -- (.75, 1);
\draw (1, 0) -- (1, 1);
\end{tikzpicture}\\
\begin{tikzpicture}[scale=1.5]
\filldraw[draw=red, ultra thick, fill=red!60] (.25,0) -- (.75, 0) -- (.75,.66) -- (.25,.66) -- cycle;
\draw[red, ultra thick] (.25,.33) -- (.75, .33);
\draw[red, ultra thick] (.25,.55) -- (.75, .55);
\draw[red, ultra thick] (.5,0) -- (.5, .66);
\draw[red, ultra thick] (.583,0) -- (.583, .66);
\draw (0,0) -- (1.25,0) -- (1.25,1) -- (0,1) -- cycle;
\draw (0,.55) -- (.75, .55);
\draw (.5, .44) -- (1.25, .44);
\draw (.583, 0) -- (.583, .66);
\draw (.666, .33) -- (.666, 1);
\draw (0,.33) -- (1.25, .33);
\draw (0,.66) -- (1.25, .66);
\draw (.25, 0) -- (.25, 1);
\draw (.5, 0) -- (.5, 1);
\draw (.75, 0) -- (.75, 1);
\draw (1, 0) -- (1, 1);
\end{tikzpicture},
\begin{tikzpicture}[scale=1.5]
\filldraw[draw=red, ultra thick, fill=red!60] (.25,0) -- (1, 0) -- (1,1) -- (.25,1) -- cycle;
\draw[red, ultra thick] (.25,.33) -- (1, .33);
\draw[red, ultra thick] (.25,.66) -- (1, .66);
\draw[red, ultra thick] (.5,0) -- (.5, 1);
\draw[red, ultra thick] (.75,0) -- (.75, 1);
\draw (0,0) -- (1.25,0) -- (1.25,1) -- (0,1) -- cycle;
\draw (0,.55) -- (.75, .55);
\draw (.5, .44) -- (1.25, .44);
\draw (.583, 0) -- (.583, .66);
\draw (.666, .33) -- (.666, 1);
\draw (0,.33) -- (1.25, .33);
\draw (0,.66) -- (1.25, .66);
\draw (.25, 0) -- (.25, 1);
\draw (.5, 0) -- (.5, 1);
\draw (.75, 0) -- (.75, 1);
\draw (1, 0) -- (1, 1);
\end{tikzpicture},
\begin{tikzpicture}[scale=1.5]
\filldraw[draw=red, ultra thick, fill=red!60] (.5,0) -- (1, 0) -- (1,.66) -- (.5,.66) -- cycle;
\draw[red, ultra thick] (.5,.33) -- (1, .33);
\draw[red, ultra thick] (.5,.44) -- (1, .44);
\draw[red, ultra thick] (.583,0) -- (.583, .66);
\draw[red, ultra thick] (.75,0) -- (.75, .66);
\draw (0,0) -- (1.25,0) -- (1.25,1) -- (0,1) -- cycle;
\draw (0,.55) -- (.75, .55);
\draw (.5, .44) -- (1.25, .44);
\draw (.583, 0) -- (.583, .66);
\draw (.666, .33) -- (.666, 1);
\draw (0,.33) -- (1.25, .33);
\draw (0,.66) -- (1.25, .66);
\draw (.25, 0) -- (.25, 1);
\draw (.5, 0) -- (.5, 1);
\draw (.75, 0) -- (.75, 1);
\draw (1, 0) -- (1, 1);
\end{tikzpicture},
\begin{tikzpicture}[scale=1.5]
\filldraw[draw=red, ultra thick, fill=red!60] (.583,0) -- (1.25, 0) -- (1.25,.66) -- (.583,.66) -- cycle;
\draw[red, ultra thick] (.583,.33) -- (1.25, .33);
\draw[red, ultra thick] (.583,.44) -- (1.25, .44);
\draw[red, ultra thick] (.75,0) -- (.75, .66);
\draw[red, ultra thick] (1,0) -- (1, .66);
\draw (0,0) -- (1.25,0) -- (1.25,1) -- (0,1) -- cycle;
\draw (0,.55) -- (.75, .55);
\draw (.5, .44) -- (1.25, .44);
\draw (.583, 0) -- (.583, .66);
\draw (.666, .33) -- (.666, 1);
\draw (0,.33) -- (1.25, .33);
\draw (0,.66) -- (1.25, .66);
\draw (.25, 0) -- (.25, 1);
\draw (.5, 0) -- (.5, 1);
\draw (.75, 0) -- (.75, 1);
\draw (1, 0) -- (1, 1);
\end{tikzpicture}
\end{tabular}
\vspace{-.25cm}

\end{minipage}\hspace{2.5cm}
}
\subfloat[]{\hspace{-.2cm}
\begin{tikzpicture}[scale=2.5]
\filldraw[draw = blue, ultra thick, fill=blue!50] (.5,.33) -- (.75,.33) -- (.75, .66) -- (.5, .66) -- cycle;
\draw[blue, ultra thick] (.5, .44) -- (.75, .44);
\draw[blue, ultra thick] (.5, .55) -- (.75, .55);
\draw[blue, ultra thick] (.583, .33) -- (.583, .66);
\draw[blue, ultra thick] (.666, .33) -- (.666, .66);
\draw (0,0) -- (1.25,0) -- (1.25,1) -- (0,1) -- cycle;
\draw (0,.55) -- (.75, .55);
\draw (.5, .44) -- (1.25, .44);
\draw (.583, 0) -- (.583, .66);
\draw (.666, .33) -- (.666, 1);
\draw (0,.33) -- (1.25, .33);
\draw (0,.66) -- (1.25, .66);
\draw (.25, 0) -- (.25, 1);
\draw (.5, 0) -- (.5, 1);
\draw (.75, 0) -- (.75, 1);
\draw (1, 0) -- (1, 1);
\draw[white] (0,-.08);
\end{tikzpicture}
}
\caption{(a) an LR-mesh $\cM$ of multiplicity 1. (b) Supports of the biquadratic LR B-splines defined on $\cM$. (c) Support of a minimal support B-spline on the mesh but not in $\mathcal{B}^{\mathcal{LR}}(\cM)$.}\label{7}
\end{figure}
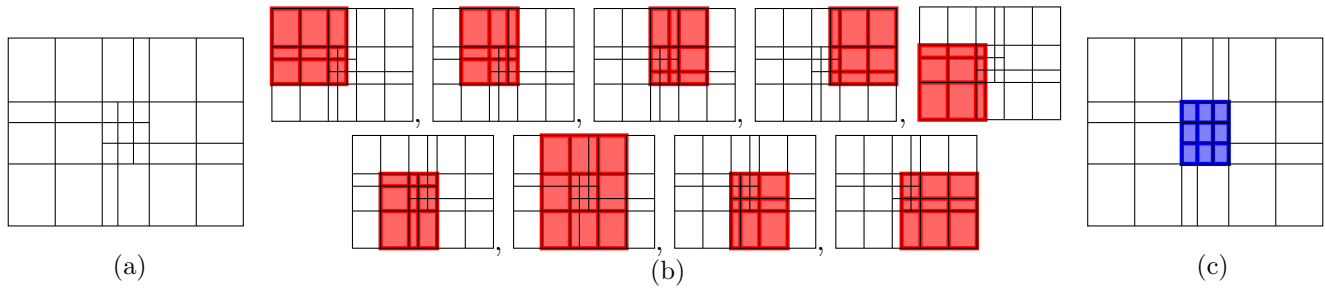\\
  In Figure \ref{7} (a) we have an LR-mesh $\cM$ of multiplicity 1. Suppose $\pmb{p}=(2,2)$. This mesh is obtained by inserting two 2-splits and two 1-splits starting from a (uniform) tensor product mesh $\cM_1$. In (b) we see the supports of LR B-splines on $\cM$, i.e., the elements of $\mathcal{B}^{\mathcal{LR}}(\cM)$ obtained by splitting the B-splines with no minimal support during the insertion of those splits. However if we look at the final mesh $\cM$ in (a), we see that there is one MS B-spline, depicted in (c), not in $\mathcal{B}^{\mathcal{LR}}(\cM)$, appearing in the mesh $\cM$.
\end{ex}
\section{Hand-in-hand principle}\label{HIHP}
In this section we describe the spanning properties of the sets $\mathcal{B}^{\mathcal{MS}}(\cM)$ and $\mathcal{B}^{\mathcal{LR}}(\cM)$. 

Any LR-mesh $\cM = \cM_N$ is defined through a sequence $\cM_{i+1} = \cM_i +\gamma_i$ starting from a tensor product mesh $\cM_1$. We know that on $\cM_1$, $\mbox{span}\,\mathcal{B}^{\mathcal{MS}}(\cM_1) = \SSS_{\pmb{p}}(\cM_1)$ as well as $\mbox{span}\,\mathcal{B}^{\mathcal{LR}}(\cM_1) = \SSS_{\pmb{p}}(\cM_1)$. We want to preserve these equalities throughout the construction of $\cM_N$. 
We require the equalities for two reasons. First, we maximize the approximation power of the considered B-splines, because the full spline space is spanned, and second, since we have a dimension formula for the spline space, we can use it to determine if they are linearly independent or not. Indeed, spanning the whole spline space, if there are more B-splines  than the dimension we can state that they are linearly dependent.   
\begin{dfn}
Let $\cM$ be an LR-mesh and $\pmb{p}$ a given bidegree. Let $\cB^{\mathcal{MS}}(\cM)$ and $\cB^{\mathcal{LR}}(\cM)$ be respectively the corresponding collection of MS and LR B-splines of degree $\pmb{p}$ on $\cM$. Let $\gamma$ be a split and $\cB^{\mathcal{MS}}(\cM+\gamma)$, $\cB^{\mathcal{LR}}(\cM + \gamma)$ be respectively the corresponding collections of MS and LR B-splines of degree $\pmb{p}$ on the LR-mesh $\cM+\gamma$. We say that:\begin{itemize}
\item $\cM+\gamma$ goes \textbf{MS-wise hand-in-hand} with $\cM$ if $$\mbox{span}~ \cB^{\mathcal{MS}}(\cM) = \SSS_{\pmb{p}}(\cM)\R \mbox{span}~ \cB^{\mathcal{MS}}(\cM+\gamma) = \SSS_{\pmb{p}}(\cM+\gamma).$$
\item $\cM+\gamma$ goes \textbf{LR-wise hand-in-hand} with $\cM$ if $$\mbox{span}~\cB^{\mathcal{LR}}(\cM) = \SSS_{\pmb{p}}(\cM)\R \mbox{span}~ \cB^{\mathcal{LR}}(\cM+\gamma) = \SSS_{\pmb{p}}(\cM+\gamma).$$
\end{itemize}
\end{dfn}
In other words, going hand-in-hand means that if on the mesh $\cM$ the span of the considered B-splines fills up the whole spline space $\SSS_{\pmb{p}}(\cM)$, then also the span of the corresponding refined B-splines living on $\cM + \gamma$ will complete the refined spline space $\SSS_{\pmb{p}}(\cM+\gamma)$.
\begin{oss}
If $\cM +\gamma$ goes LR-wise hand-in-hand with $\cM$, then $\cM + \gamma$ also goes MS-wise hand-in-hand with $\cM$. This is trivial because $\cB^{\mathcal{LR}}(\cM + \gamma) \subseteq \cB^{\mathcal{MS}}(\cM+\gamma)$. The converse is not true in general.
\end{oss}
The next result, Theorem \ref{hiht}, gives us a sufficient and necessary condition for going hand-in-hand. First, we recall from \cite{dokken2013polynomial} the restriction of a B-spline to a split.
\begin{dfn}
Let $B[\pmb{x}, \pmb{y}] =B[\pmb{x}]B[\pmb{y}]$ be a tensor product B-spline
and $\gamma$ be a $k$-split for some $k\in \{1,2\}$ traversing $\mbox{supp}\,B[\pmb{x}, \pmb{y}]$. Then we define the \textbf{restriction of $B[\pmb{x},\pmb{y}]$ onto $\gamma$}, $B_\gamma$, as the univariate B-spline $B_\gamma=B[\pmb{x}]$ if $k=2$ or $B_\gamma = B[\pmb{y}]$ if $k = 1$.
\end{dfn}
We are now ready for the sufficient and necessary condition to go hand-in-hand. This is a revisited statement of \cite[Theorem 5.10]{dokken2013polynomial}.
\begin{thm}\label{hiht}
Let $\cM$ be an LR-mesh in $\RR^2$ and $\pmb{p}=(p_1,p_2)$ a given bidegree. Let $\gamma$ be a $k$-split of multiplicity $\mu(\gamma)$ corresponding to an expanded $k$-split $\bar{\gamma}$ of length $p_{3-k}+r+1$ with $r\geq1$ in $\cM + \gamma$. Let \begin{equation*}
\begin{split}
\mathcal{B}^{\mathcal{MS}}(\bar{\gamma})&=\{B[\pmb{x}, \pmb{y}] \in \mathcal{B}^{\mathcal{MS}}(\cM+\gamma): \bar{\gamma}\mbox{ traverses its support } \}\\\vspace{.35cm}
\mathcal{B}^{\mathcal{LR}}(\bar{\gamma}) &=\{B[\pmb{x},\pmb{y}]\in \mathcal{B}^{\mathcal{LR}}(\cM + \gamma):  \bar{\gamma}\mbox{ traverses its support }\}
\end{split}
\end{equation*}
Assume that $span\,\mathcal{B}^{\mathcal{MS}}(\cM) = \SSS_{\pmb{p}}(\cM)=span\,\mathcal{B}^{\mathcal{LR}}(\cM)$. Then $\cM + \gamma$ goes MS-wise hand-in-hand with $\cM$ if and only if the following condition on the B-spline restrictions to $\bar{\gamma}$ holds $$
\dim span\, \{B_{\bar{\gamma}}\}_{B\in \mathcal{B}^{\mathcal{MS}}(\bar{\gamma})} = r.
$$
Similarly, $\cM+\gamma$ goes LR-wise hand-in-hand with $\cM$ if and only if 
$$
\dim span\, \{B_{\bar{\gamma}}\}_{B\in \mathcal{B}^{\mathcal{LR}}(\bar{\gamma})} = r.
$$
\end{thm}
The sets $\mathcal{B}^{\mathcal{MS}}(\bar{\gamma})$ and $\mathcal{B}^{\mathcal{LR}}(\bar{\gamma})$ contain new MS and LR B-splines respectively, created after the insertion of the split $\gamma$, with part of $\bar{\gamma}$ as internal split. Of course, since $\{B_{\bar{\gamma}}\}_{B\in\mathcal{B}^{\mathcal{MS}}(\bar{\gamma})}$ and $\{B_{\bar{\gamma}}\}_{B\in\mathcal{B}^{\mathcal{LR}}(\bar{\gamma})}$ are contained in $\SSS_{p_{3-k}}(\pmb{t}_{\bar{\gamma}})$ with $\pmb{t}_{\bar{\gamma}}$ the knot vector on $\bar{\gamma}$, we have $$
\dim span \{B_{\bar{\gamma}}\}_{B\in\mathcal{B}^{\mathcal{MS}}(\bar{\gamma})} \leq r \quad \mbox{ and }\quad \dim span\{B_{\bar{\gamma}}\}_{B\in\mathcal{B}^{\mathcal{LR}}(\bar{\gamma})} \leq r.
$$
We distinguish two cases when these are strict inequalities: \begin{enumerate}
\item The cardinality of $\mathcal{B}^{\mathcal{MS}}(\bar{\gamma})$, or $\mathcal{B}^{\mathcal{LR}}(\bar{\gamma})$ respectively, is less than $r$.
\item The cardinality of $\mathcal{B}^{\mathcal{MS}}(\bar{\gamma})$, or $\mathcal{B}^{\mathcal{LR}}(\bar{\gamma})$ respectively, is at least $r$ but there are fewer than $r$ linearly independent restrictions $B_{\bar{\gamma}}$.
\end{enumerate}
In the first case, the cardinality of such sets depends on the mutual position of the splits in $\cM + \gamma$. However, we can always guarantee that $\cB^{\mathcal{MS}}(\bar{\gamma})$ and $\cB^{\mathcal{LR}}(\bar{\gamma})$ have at least $r$ elements by slightly modifying the mesh $\cM$ or modifying the length of $\gamma$. 
\begin{ex}\label{case1}
Consider bidegree $\pmb{p} =(2,2)$ and a 1-split $\gamma$ of length $4$ to insert into the LR-mesh $\cM$ of multiplicity 1 shown in Figures \ref{hh1} (a). We recall that $\dim \SSS_{\pmb{p}}(\cM+\gamma) = \dim\SSS_{\pmb{p}}(\cM)+1$. Therefore, $\cM+\gamma$ will go MS-wise or LR-wise hand-in-hand with $\cM$ if we  get a new B-spline of the considered kind.
\begin{figure}[ht!]\centering
\subfloat[]{\centering
\begin{tikzpicture}[scale=2.35]
\draw (0,0) -- (1.25,0) -- (1.25,1) -- (0,1) -- cycle;
\draw (0,.77) -- (.75, .77);
\draw (.5, .44) -- (1.25, .44);
\draw[dashed, red, ultra thick] (.625, .33) -- (.625,.77);
\draw (0,.33) -- (1.25, .33);
\draw (0,.66) -- (1.25, .66);
\draw (.25, 0) -- (.25, 1);
\draw (.5, 0) -- (.5, 1);
\draw (.75, 0) -- (.75, 1);
\draw (1, 0) -- (1, 1);
\end{tikzpicture}}\quad
\subfloat[]{\centering
\begin{tikzpicture}[scale=2.35]
\filldraw[draw =cyan, ultra thick, fill=cyan!50] (.5, .33) -- (1, .33) -- (1, .77) -- (.5, .77) -- cycle;
\fill[pattern=north west lines, pattern color=cyan!50] (.5, .33) -- (1, .33) -- (1, .77) -- (.5, .77) -- cycle;
\draw (0,0) -- (1.25,0) -- (1.25,1) -- (0,1) -- cycle;
\draw (0,.77) -- (.75, .77);
\draw[dashed, ultra thick] (.75, .77) -- (1, .77);
\draw (.5, .44) -- (1.25, .44);
\draw (.625, .33) -- (.625,.77);
\draw (0,.33) -- (1.25, .33);
\draw (0,.66) -- (1.25, .66);
\draw (.25, 0) -- (.25, 1);
\draw (.5, 0) -- (.5, 1);
\draw (.75, 0) -- (.75, 1);
\draw (1, 0) -- (1, 1);
\end{tikzpicture}}\quad
\subfloat[]{\centering
\begin{tikzpicture}[scale=2.35]
\filldraw[draw=green(ryb), ultra thick, fill=green(ryb), opacity=.5] (.5, .33) -- (1,.33) -- (1,.77) -- (.5, .77) -- cycle;
\filldraw[draw=orange, ultra thick, fill = orange, opacity=.5] (.625, .33) -- (1.25,.33) -- (1.25,.77) -- (.625, .77) -- cycle;
\fill[pattern=north west lines, pattern color=green(ryb)]  (.5, .33) -- (1,.33) -- (1,.77) -- (.5, .77) -- cycle;
\fill[pattern=north east lines, pattern color=orange] (.625, .33) -- (1.25,.33) -- (1.25,.77) -- (.625, .77) -- cycle;
\draw (0,0) -- (1.25,0) -- (1.25,1) -- (0,1) -- cycle;
\draw (0,.77) -- (.75, .77);
\draw[dashed, ultra thick] (.75, .77) -- (1.25, .77);
\draw (.5, .44) -- (1.25, .44);
\draw (.625, .33) -- (.625,.77);
\draw (0,.33) -- (1.25, .33);
\draw (0,.66) -- (1.25, .66);
\draw (.25, 0) -- (.25, 1);
\draw (.5, 0) -- (.5, 1);
\draw (.75, 0) -- (.75, 1);
\draw (1, 0) -- (1, 1);
\end{tikzpicture}}\quad
\subfloat[]{\centering
\begin{tikzpicture}[scale=2.35]
\filldraw[draw=green(ryb), ultra thick, fill=green(ryb),opacity=.5] (.25, .33) -- (.75,.33) -- (.75,.77) -- (.25, .77) -- cycle;
\filldraw[draw=orange, ultra thick, fill=orange,opacity=.5] (.5, .33) -- (1,.33) -- (1,.77) -- (.5, .77) -- cycle;
\fill[pattern=north west lines, pattern color=green(ryb)]  (.25, .33) -- (.75,.33) -- (.75,.77) -- (.25, .77) -- cycle;
\fill[pattern=north east lines, pattern color=orange] (.5, .33) -- (1,.33) -- (1,.77) -- (.5, .77) -- cycle;
\draw (0,0) -- (1.25,0) -- (1.25,1) -- (0,1) -- cycle;
\draw (0,.77) -- (.75, .77);
\draw[dashed, ultra thick] (.75, .77) -- (1,.77);
\draw (.5, .44) -- (1.25, .44);
\draw[dashed, ultra thick] (.25, .44) -- (.5, .44);
\draw (.625, .33) -- (.625,.77);
\draw (0,.33) -- (1.25, .33);
\draw (0,.66) -- (1.25, .66);
\draw (.25, 0) -- (.25, 1);
\draw (.5, 0) -- (.5, 1);
\draw (.75, 0) -- (.75, 1);
\draw (1, 0) -- (1, 1);
\end{tikzpicture}}\quad
\subfloat[]{\centering
\begin{tikzpicture}[scale=2.35]
\filldraw[draw=green(ryb), ultra thick, fill=green(ryb), opacity=.5] (.25, .33) -- (.75, .33) -- (.75, 1) -- (.25, 1) -- cycle;
\filldraw[draw=orange, ultra thick, fill=orange,opacity=.5] (.5, .33) -- (1, .33) -- (1, 1) -- (.5, 1) -- cycle;
\fill[pattern=north west lines, pattern color=green(ryb)] (.25, .33) -- (.75, .33) -- (.75, 1) -- (.25, 1) -- cycle;
\fill[pattern=north east lines, pattern color=orange] (.5, .33) -- (1, .33) -- (1, 1) -- (.5, 1) -- cycle;
\draw (0,0) -- (1.25,0) -- (1.25,1) -- (0,1) -- cycle;
\draw (0,.77) -- (.75, .77);
\draw[dashed, ultra thick] (.625, .77) -- (.625, 1);
\draw (.5, .44) -- (1.25, .44);
\draw (.625, .33) -- (.625,.77);
\draw (0,.33) -- (1.25, .33);
\draw (0,.66) -- (1.25, .66);
\draw (.25, 0) -- (.25, 1);
\draw (.5, 0) -- (.5, 1);
\draw (.75, 0) -- (.75, 1);
\draw (1, 0) -- (1, 1);
\end{tikzpicture}}
\caption{(a) LR-mesh $\cM$ of multiplicity 1 and a new split (dashed) to insert. (b) modification of $\cM$ (dashed) to go MS-wise hand-in-hand. (c),(d),(e) modification of $\cM$ (dashed) to go LR-wise hand-in-hand.}\label{hh1}
\end{figure}
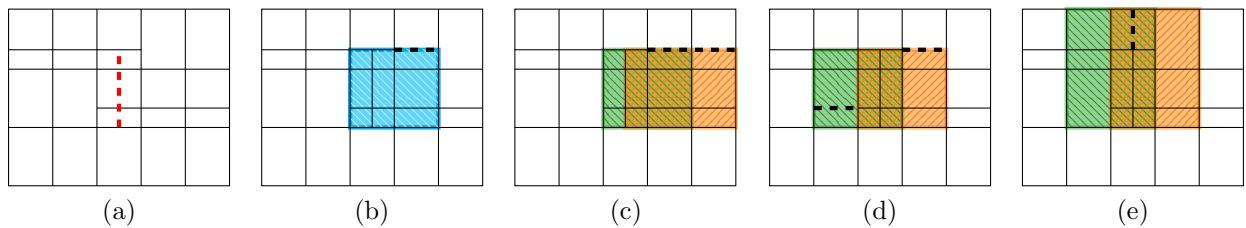\\
However, as a result of the split mutual positions, in $\cM+\gamma$ we do not generate any new B-spline, so that $\cB^{\mathcal{MS}}(\cM+\gamma) = \cB^{\mathcal{MS}}(\cM)$ and $\cB^{\mathcal{LR}}(\cM + \gamma) = \cB^{\mathcal{LR}}(\cM)$. Thus $\cM+\gamma$ does not go neither LR-wise nor MS-wise hand-in-hand with $\cM$.

However, if we extend by one meshline a split on $\cM$ before inserting $\gamma$, see Figure \ref{hh1} (b), we create from scratch a new MS B-spline, whose support is highlighted in blue.

While, if we extend by two the same split, Figure \ref{hh1} (c), or we extend by one both the splits, Figure \ref{hh1} (d), there will be an LR B-spline on these modified meshes to refine after the insertion of $\gamma$.

Another strategy is to insert a longer $\gamma$. Indeed, if we decide to insert $\gamma$ one meshline longer, as in Figure \ref{hh1} (e), we go LR-wise, and so MS-wise, hand-in-hand, splitting the two LR B-splines that were defined in the upper left and upper right corner of $\cM$.
\end{ex}
However, although the cardinality of such sets is sufficiently large, the linearly independent restrictions $B_{\bar{\gamma}}$ can be insufficient for spanning the whole spline space $\SSS_{p_{3-k}}(\pmb{t}_{\bar{\gamma}})$. 
\begin{ex}
Consider bidegree $(2,2)$, the LR-mesh $\cM$ of multiplicity 1 and the new 2-split $\gamma$ depicted in Figure \ref{hh5} (a).
\begin{figure}[ht!]\vspace{1cm}\centering
\subfloat[]{\centering
\begin{tikzpicture}[scale = 4]
\draw (-.25,0) -- (1, 0) -- (1, 1) -- (-.25, 1) --cycle;
\draw (0,0) -- (0,1);
\draw (-.25,.25) -- (1, .25);
\draw (-.25, .5) -- (1, .5);
\draw (-.25, .75) -- (1, .75);
\draw (.25, 0) -- (.25, 1);
\draw (.5, 0) -- (.5, 1); 
\draw (.75, 0) -- (.75, 1);
\draw (.125, .25) -- (.125, 1);
\draw (.25,.44) -- (1, .44);
\draw (.25, .31) -- (1, .31);
\draw (.33, .25) -- (.33, .5);
\draw (.41, .25) -- (.41, .5);
\draw[dashed, thick, red] (-.25,.38) -- (.5,.38);
\draw[white] (-.25, -.3) -- (1,-.3);
\fill (0,.38) circle (.015);
\fill (.125, .38) circle (.015);
\fill (.25,.38) circle (.015);
\fill (.33, .38) circle (.015);
\fill (.41,.38) circle (.015);
\fill (.5, .38) circle (.015);
\fill (-.25,.38) circle (.015);
\end{tikzpicture}}
\subfloat[]{\centering
\begin{minipage}{5.7cm} \vspace{-6cm} 
 \begin{tabular}{c}
\begin{tikzpicture}[scale=2]
\fill[orange!40] (-.25, .25) -- (.25,.25) -- (.25, .75) -- (-.25,.75) -- cycle;
\draw[orange, ultra thick] (-.25, .5) -- (.25,.5);
\draw[orange, ultra thick] (-.25, .38) -- (.25,.38);
\draw[orange, ultra thick] (0, .25) -- (0,.75);
\draw[orange, ultra thick] (.125, .25) -- (.125,.75);
\draw[orange, ultra thick] (-.25, .25) -- (.25,.25) -- (.25, .75) -- (-.25,.75) -- cycle;
\draw (-.25,0) -- (1, 0) -- (1, 1) -- (-.25, 1) --cycle;
\draw (0,0) -- (0,1);
\draw (-.25,.25) -- (1, .25);
\draw (-.25, .5) -- (1, .5);
\draw (-.25, .75) -- (1, .75);
\draw (.25, 0) -- (.25, 1);
\draw (.5, 0) -- (.5, 1); 
\draw (.75, 0) -- (.75, 1);
\draw (.125, .25) -- (.125, 1);
\draw (.25,.44) -- (1, .44);
\draw (.25, .31) -- (1, .31);
\draw (.33, .25) -- (.33, .5);
\draw (.41, .25) -- (.41, .5);
\draw (-.25,.38) -- (.5,.38);
\fill (0,.38) circle (.015);
\fill (.125, .38) circle (.015);
\fill (.25,.38) circle (.015);
\fill (.33, .38) circle (.015);
\fill (.41,.38) circle (.015);
\fill (.5, .38) circle (.015);
\fill (-.25,.38) circle (.015);
\end{tikzpicture}~
\begin{tikzpicture}[scale=2]
\fill[green(ryb)!50] (0, .25) -- (.5,.25) -- (.5, .75) -- (0,.75) -- cycle;
\draw[green(ryb), ultra thick] (0, .5) -- (.5,.5);
\draw[green(ryb), ultra thick] (0, .38) -- (.5,.38);
\draw[green(ryb), ultra thick] (.25, .25) -- (.25,.75);
\draw[green(ryb), ultra thick] (.125, .25) -- (.125,.75);
\draw[green(ryb), ultra thick] (0, .25) -- (.5,.25) -- (.5, .75) -- (0,.75) -- cycle;
\draw (-.25,0) -- (1, 0) -- (1, 1) -- (-.25, 1) --cycle;
\draw (0,0) -- (0,1);
\draw (-.25,.25) -- (1, .25);
\draw (-.25, .5) -- (1, .5);
\draw (-.25, .75) -- (1, .75);
\draw (.25, 0) -- (.25, 1);
\draw (.5, 0) -- (.5, 1); 
\draw (.75, 0) -- (.75, 1);
\draw (.125, .25) -- (.125, 1);
\draw (.25,.44) -- (1, .44);
\draw (.25, .31) -- (1, .31);
\draw (.33, .25) -- (.33, .5);
\draw (.41, .25) -- (.41, .5);
\draw (-.25,.38) -- (.5,.38);
\fill (0,.38) circle (.015);
\fill (.125, .38) circle (.015);
\fill (.25,.38) circle (.015);
\fill (.33, .38) circle (.015);
\fill (.41,.38) circle (.015);
\fill (.5, .38) circle (.015);
\fill (-.25,.38) circle (.015);
\end{tikzpicture}\\
\begin{tikzpicture}[scale=2]
\fill[cyan!50] (-.25, 0) -- (.5,0) -- (.5, .5) -- (-.25,.5) -- cycle;
\draw[cyan, ultra thick] (-.25, .25) -- (.5,.25);
\draw[cyan, ultra thick] (-.25, .38) -- (.5,.38);
\draw[cyan, ultra thick] (0, 0) -- (0,.5);
\draw[cyan, ultra thick] (.25, 0) -- (.25,.5);
\draw[cyan, ultra thick] (-.25, 0) -- (.5,0) -- (.5, .5) -- (-.25,.5) -- cycle;
\draw (-.25,0) -- (1, 0) -- (1, 1) -- (-.25, 1) --cycle;
\draw (0,0) -- (0,1);
\draw (-.25,.25) -- (1, .25);
\draw (-.25, .5) -- (1, .5);
\draw (-.25, .75) -- (1, .75);
\draw (.25, 0) -- (.25, 1);
\draw (.5, 0) -- (.5, 1); 
\draw (.75, 0) -- (.75, 1);
\draw (.125, .25) -- (.125, 1);
\draw (.25,.44) -- (1, .44);
\draw (.25, .31) -- (1, .31);
\draw (.33, .25) -- (.33, .5);
\draw (.41, .25) -- (.41, .5);
\draw (-.25,.38) -- (.5,.38);
\fill (0,.38) circle (.015);
\fill (.125, .38) circle (.015);
\fill (.25,.38) circle (.015);
\fill (.33, .38) circle (.015);
\fill (.41,.38) circle (.015);
\fill (.5, .38) circle (.015);
\fill (-.25,.38) circle (.015);
\end{tikzpicture}\\
\begin{tikzpicture}[scale=2]
\fill[magenta!40] (.25, .25) -- (.5,.25) -- (.5, .44) -- (.25,.44) -- cycle;
\draw[magenta, ultra thick] (.25, .31) -- (.5,.31);
\draw[magenta, ultra thick] (.25, .38) -- (.5,.38);
\draw[magenta, ultra thick] (.33, .25) -- (.33,.44);
\draw[magenta, ultra thick] (.41, .25) -- (.41,.44);
\draw[magenta, ultra thick] (.25, .25) -- (.5,.25) -- (.5, .44) -- (.25,.44) -- cycle;
\draw (-.25,0) -- (1, 0) -- (1, 1) -- (-.25, 1) --cycle;
\draw (0,0) -- (0,1);
\draw (-.25,.25) -- (1, .25);
\draw (-.25, .5) -- (1, .5);
\draw (-.25, .75) -- (1, .75);
\draw (.25, 0) -- (.25, 1);
\draw (.5, 0) -- (.5, 1); 
\draw (.75, 0) -- (.75, 1);
\draw (.125, .25) -- (.125, 1);
\draw (.25,.44) -- (1, .44);
\draw (.25, .31) -- (1, .31);
\draw (.33, .25) -- (.33, .5);
\draw (.41, .25) -- (.41, .5);
\draw (-.25,.38) -- (.5,.38);
\fill (0,.38) circle (.015);
\fill (.125, .38) circle (.015);
\fill (.25,.38) circle (.015);
\fill (.33, .38) circle (.015);
\fill (.41,.38) circle (.015);
\fill (.5, .38) circle (.015);
\fill (-.25,.38) circle (.015);
\end{tikzpicture}~
\begin{tikzpicture}[scale=2]
\fill[magenta!40] (.25, .31) -- (.5,.31) -- (.5, .5) -- (.25,.5) -- cycle;
\draw[magenta, ultra thick] (.25, .38) -- (.5,.38);
\draw[magenta, ultra thick] (.25, .44) -- (.5,.44);
\draw[magenta, ultra thick] (.33, .31) -- (.33,.5);
\draw[magenta, ultra thick] (.41, .31) -- (.41,.5);
\draw[magenta, ultra thick] (.25, .31) -- (.5,.31) -- (.5, .5) -- (.25,.5) -- cycle;
\draw (-.25,0) -- (1, 0) -- (1, 1) -- (-.25, 1) --cycle;
\draw (0,0) -- (0,1);
\draw (-.25,.25) -- (1, .25);
\draw (-.25, .5) -- (1, .5);
\draw (-.25, .75) -- (1, .75);
\draw (.25, 0) -- (.25, 1);
\draw (.5, 0) -- (.5, 1); 
\draw (.75, 0) -- (.75, 1);
\draw (.125, .25) -- (.125, 1);
\draw (.25,.44) -- (1, .44);
\draw (.25, .31) -- (1, .31);
\draw (.33, .25) -- (.33, .5);
\draw (.41, .25) -- (.41, .5);
\draw (-.25,.38) -- (.5,.38);
\fill (0,.38) circle (.015);
\fill (.125, .38) circle (.015);
\fill (.25,.38) circle (.015);
\fill (.33, .38) circle (.015);
\fill (.41,.38) circle (.015);
\fill (.5, .38) circle (.015);
\fill (-.25,.38) circle (.015);
\end{tikzpicture}
\end{tabular}\vspace{-.15cm}
\end{minipage}
}
\subfloat[]{\begin{minipage}{6cm}\vspace{-4.65cm}\includegraphics[width=6cm,trim={2.25cm 0 1.85cm 1.55cm},clip]{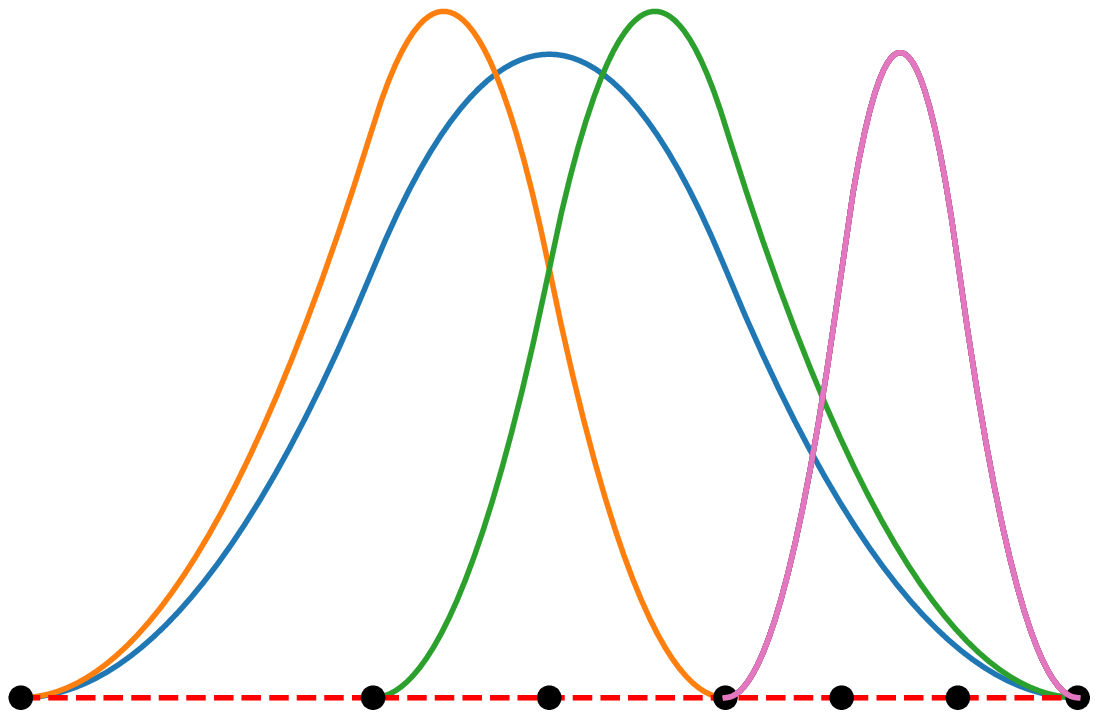}
\begin{minipage}{6cm}
\vspace{-10.5cm} $\left.\right.\qquad\quad~~ \textcolor{orange}{B_{\bar{\gamma}}^1}\quad~ \textcolor{azure}{B_{\bar{\gamma}}^3} \quad~ \textcolor{green(ryb)}{B_{\bar{\gamma}}^2}~~ \textcolor{magenta}{B_{\bar{\gamma}}^4} = \textcolor{magenta}{B_{\bar{\gamma}}^5}$
\end{minipage}
\end{minipage}}
\caption{(a) LR-mesh $\cM$ of multiplicity 1 and a new 2-split $\gamma$ (dashed) with their intersections (black dots). Consider bidegree $(2,2)$. In (b) the supports of the LR B-splines $B^1, B^2$ (top), $B^3$ (center), $B^4, B^5$ (bottom) in $\cM + \gamma$ with part of $\gamma$ as internal split. In (c) the restrictions of these LR B-splines to $\gamma$. 
}\label{hh5}
\end{figure}\\
$\gamma$ has length 7 so that $\dim\SSS_{\pmb{p}}(\cM+\gamma) = \dim \SSS_{\pmb{p}}(\cM) + 4$. Moreover, it easy to check that $\cM$ can be constructed LR-wise hand-in-hand. Therefore $\cB^{\mathcal{LR}}(\cM) = \cB^{\mathcal{MS}}(\cM)$ and their span fulfills the spline space $\SSS_{\pmb{p}}(\cM)$. 

Now, when $\gamma$ is inserted, 
there are 5 LR B-splines, $B^1, B^2,B^3,B^4, B^5$ in $\mathcal{B}^{\mathcal{LR}}(\bar{\gamma})$, see Figure \ref{hh5} (b).  Thus, the cardinalities $|\mathcal{B}^{\mathcal{LR}}(\bar{\gamma})|,|\mathcal{B}^{\mathcal{MS}}(\bar{\gamma})|$ are large enough for $\cM +\gamma$ to go hand-in-hand with $\cM$. However, if we look at the restrictions on $\bar{\gamma}$, depicted in Figure \ref{hh5} (c), we can see that the LR B-splines $B^4, B^5$ have same restriction, while the restriction of $B_{\bar{\gamma}}^3$ can be easily written, via knot insertion, as a linear combination of $B_{\bar{\gamma}}^1, B_{\bar{\gamma}}^2$.

Thus, there are only 3 linearly independent restrictions on $\bar{\gamma}$: two among the $B^1_{\bar{\gamma}}, B_{\bar{\gamma}}^2, B_{\bar{\gamma}}^3$ and one among $B_{\bar{\gamma}}^4, B_{\bar{\gamma}}^5$. The new mesh $\cM + \gamma$, therefore, cannot go neither LR-wise nor MS-wise hand-in-hand with $\cM$.
\end{ex}
However, if the expanded $k$-split $\bar{\gamma}$ has length $p_{3-k}+2$ or $p_{3-k}+3$, this phenomenon cannot happen. Indeed, if $\bar{\gamma}$ is $p_{3-k}+2$ long, the spline space has grown by one and the restrictions on $\bar{\gamma}$ are at least one. While if $\bar{\gamma}$ is $p_{3-k}+3$ long, there are at least two different (and so linearly independent) restrictions to allocate $p_{3-k}+3$ knots into knot vectors of $p_{3-k}+2$ components. 

We will need this observation further on to derive the minimal number of LR B-splines needed for a linear dependence relation.
\section{Characterization of linear dependence in $\mathcal{B}^{\mathcal{MS}}(\cM)$}\label{LDMS}
The purpose of this section is to investigate the minimal number of MS B-splines required for a linear dependence relation on an LR-mesh $\cM$ and features needed in such configurations. 

In particular, the main results of this section are that at least six MS B-splines are necessary for a linear dependence (Proposition \ref{sixMS}) for any bidegree $\pmb{p} = (p_1, p_2)$ with $p_1+p_2 \geq 1$ and that in a configuration of linear dependence with exactly six B-splines, one of them is not an LR B-spline (Proposition \ref{noLR}). We achieve these results by looking at the minimal number of B-splines needed to satisfy necessary conditions for having a linear dependence relation. First we introduce the \emph{nestedness condition} (Proposition \ref{L1}): at any corner of the region of the mesh where we have linear dependence, there is a B-spline in the linear dependence relation whose support is fully contained in the support of another larger B-spline in the linear dependence relation as well. We show that at least five MS B-splines are needed to satisfy this condition (Corollary \ref{fiveMS}). Then we have to prove that it is impossible to have a linear dependence with only these five B-splines. Therefore, first we show the possible arrangements of the supports in the case a linear dependence relation involves only five B-splines (Lemma \ref{cov}). Then we introduce another necessary condition for linear dependencies regarding the T-junctions in the region of the mesh where the linear dependence occurs (Corollary \ref{Tver}). This new condition narrows the possible arrangements of the supports found in Lemma \ref{cov}. Finally, by looking at the position of the splits of the five B-splines in this remaining configurations, one can prove Proposition \ref{sixMS} mentioned above.

We recall that our meaning of linearly dependent functions is slightly different from the standard definition. We consider linearly dependent only those functions that are active, i.e., that have nonzero coefficient, in the dependence relation.  
\begin{dfn}
Given an LR-mesh $\cM$ and two MS B-splines $B[\pmb{x}^1, \pmb{y}^1]$, $B[\pmb{x}^2, \pmb{y}^2]$, defined on $\cM$, we say that $B[\pmb{x}^1, \pmb{y}^1]$ \textbf{is nested into} $B[\pmb{x}^2, \pmb{y}^2]$ if $\mbox{supp}\,B[\pmb{x}^1, \pmb{y}^1] \subset \mbox{supp}\, B[\pmb{x}^2, \pmb{y}^2]$ and $\mbox{supp}\,B[\pmb{x}^1,\pmb{y}^1],\mbox{supp}\, B[\pmb{x}^2, \pmb{y}^2]$ share one, and only one, vertex. 
\end{dfn}
\begin{prop}[Nestedness condition]\label{L1}
Consider the set of MS B-splines of bidegree $\pmb{p}=(p_1,p_2)$ defined on $\cM$, $\mathcal{B}^{\mathcal{MS}}(\cM)$. Let $\mathcal{B} =\{B[\pmb{x}^j, \pmb{y}^j]\}_{j=1}^n\subseteq \mathcal{B}^{\mathcal{MS}}(\cM)$ be a subset of linearly dependent MS B-splines with
$
\pmb{x}^j = (x_1^j, \ldots, x_{p_1+2}^j), \pmb{y}^j = (y_1^j, \ldots, y_{p_2+2}^j),  j=1,\ldots, n.$
Let $\mathcal{R}$ be the region given by the union of their supports. Let $(\bar{x}, \bar{y})$ be any (convex) corner in $\mathcal{R}$ and consider the set $\mathcal{B}' \subseteq \mathcal{B}$ of B-splines that have $(\bar{x}, \bar{y})$ as pair of knots. Then, define the set of MS B-splines with smallest support, in both directions, in the corner $(\bar{x},\bar{y})$ of $\mathcal{R}$: $$
\mathcal{L} = \{B[\pmb{x}, \pmb{y}] \in \mathcal{B}': x_{p_1+2}-\bar{x},~ y_{p_2+2} - \bar{y}\mbox{ are minimal }\}.
$$
Let $\bar{x}^{m}$ and $\bar{y}^{m}$ be the values of the end knots in $\pmb{x}$ and $\pmb{y}$, respectively, of any MS B-spline in $\mathcal{L}$.
Then \begin{enumerate}
\item $\mathcal{L}$ has a unique B-spline $B[\pmb{x}^m, \pmb{y}^m]$ (with $x_1^m = \bar{x},~ y_1^m = \bar{y}$ and $x_{p_1+2}^m = \bar{x}^{m},~ y_{p_2+2}^m = \bar{y}^{m}$).\vspace{.15cm}
\item There exists another B-spline $B[\pmb{x}^\ell, \pmb{y}^\ell] \in \mathcal{B}'$ with $y_{p_2+2}^\ell > y_{p_2+2}^m$, $x_{p_1+2}^\ell > x_{p_1+2}^m$, i.e., 
$B[\pmb{x}^m, \pmb{y}^m]$ is nested into $B[\pmb{x}^\ell, \pmb{y}^\ell]$ and the part of $\mbox{supp}\,B[\pmb{x}^\ell, \pmb{y}^\ell]$ exceeding the support of $B[\pmb{x}^m, \pmb{y}^m]$ is therefore L-shaped.
\end{enumerate}
\begin{proof}
\begin{enumerate}
\item Let us first show that $\mathcal{L} \neq \emptyset$. Since in $\mathcal{R}$ there is a linear dependence relation, every point of it must be inside the support of at least two B-splines.  Consider the element of the box-partition in $\mathcal{R}$ that has $(\bar{x}, \bar{y})$ as vertex. if $\mathcal{L} = \emptyset$, it would mean that such element in the corner of $\mathcal{R}$ is covered by at least the support of two B-splines $B^1 = B[\pmb{x}^1,\pmb{y}^1]$ and $B^2 = B[\pmb{x}^2, \pmb{y}^2]$ such that $B^2$ is taller than $B^1$ but narrower as reported in Figure \ref{9}. 
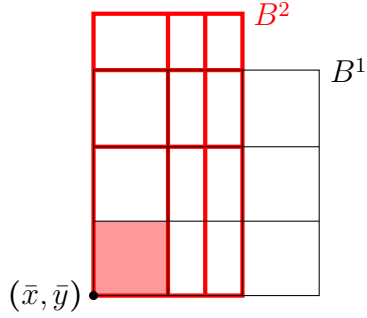
\begin{figure}[h!]\centering
\begin{tikzpicture}[scale=3]
\fill[red!40] (0,0) -- (.33,0) -- (.33,.33) -- (0,.33) -- cycle;
\draw[red, ultra thick] (0,0) -- (.66,0) -- (.66, 1.25) -- (0,1.25) -- cycle;
\draw[red, ultra thick] (.33, 0) -- (.33, 1.25);
\draw[red, ultra thick] (.495, 0) -- (.495, 1.25);
\draw[red, ultra thick] (0,.66) -- (.66,.66);
\draw[red, ultra thick] (0,1) -- (.66, 1);
\draw (0,0) -- (1,0) -- (1,1) -- (0,1) -- cycle;
\draw (0,.33) -- (1,.33);
\draw (0,.66) -- (1,.66);
\draw (.33,0) -- (.33, 1);
\draw (.66, 0) -- (.66,1);
\draw (1, 1) node[right]{$B^1$};
\draw (.66, 1.25) node[right]{\textcolor{red}{$B^2$}};
\fill (0,0) node[left]{$(\bar{x}, \bar{y})$} circle (.02);
\draw[white] (1.66,0) -- (1.66,1);
\end{tikzpicture}
\caption{Consider bidegree (2,2). The element of $\cE$ with $(\bar{x}, \bar{y})$ as one of its vertices is highlighted. If $\mathcal{L}$ is empty, this element is covered by two B-splines $B^1, B^2$ with $B^2$ taller but narrower than $B^1$. Their supports are represented in the Figure. There is at least one horizontal split of $B^1$ traversing $B^2$ that is not a split of it (and at least one vertical split of $B^2$ traversing $B^1$ and not a split of it).}\label{9}
\end{figure}\\
Thus, there are $p_2+1$ horizontal splits of $B^1$ traversing the interior of $\mbox{supp}\,B^2$. Only $p_2$ of them (at most) can be also splits of $B^2$. This is a contradiction because an extra split traverses the support of $B^2$ and so it has not minimal support on the mesh. Hence $|\mathcal{L}| \geq 1$. Let us assume there are two MS  B-splines in $\mathcal{L}$, $B^1 = B[\pmb{x}^1,\pmb{y}^1]$ and $B^2 = B[\pmb{x}^2, \pmb{y}^2]$. So $$
\begin{array}{ll} x_1^1 = x_1^2 = \bar{x} & y_1^1 = y_1^2 = \bar{y}\\\\x_{p_1+2}^1 = x_{p_1+2}^2 = \bar{x}^{m} & y_{p_2+2}^1 = y_{p_2+2}^2 = \bar{y}^{m}.\end{array}
$$
If also the internal knots of $B^1$ and $B^2$ are the same in both directions, it would mean that $B^2 = B^1$ and there is nothing to prove. Thus, let us assume there is at least one different knot in the $x$- or $y$-direction. For instance, suppose there is a different internal knot $x_i^2\in\pmb{x}^2$ for some $i$, with respect to $\pmb{x}^1$.

Then the corresponding vertical split $\{x_i^2\}\times[\bar{y}, \bar{y}^{m}]$ of $B^2$ would traverse the  support of $B^1$ without being a split of $B^1$. This is a contradition because $B^1$ has minimal support.

\item $B^m = B[\pmb{x}^m, \pmb{y}^m]$ is in a linear dependence relation, so every element of the box-partition $\cE$ contained in $\mbox{supp}\,B^m$ must be covered, i.e.,  it is into the support of at least another MS B-spline of $\mathcal{B}$. 

Consider then the element of $\cE$ contained in $\mbox{supp} B^m$ that has $(x_1^m, y_1^m)$ as one of its vertices.  
If there was a MS B-spline $B^i = B[\pmb{x}^i, \pmb{y}^i] \in \mathcal{B}$, for some $i$, fully contained in the support of $B^m$ covering such an element then it would have $x_{p_1+2}^i$ and $y_{p_2+2}^i$ smaller than $x_{p_1+2}^m$ and $y_{p_2+2}^m$ respectively, which is contradicting the minimality of $B^m$. Therefore, such an element can be covered only by a MS B-spline whose support exceeds on the right, or on the top, or both on the right and on the top, the support of $B^m$. Using the argument we used to prove that $|\mathcal{L}| \neq \emptyset$, one shows that only the last case can happen.
Thus, the rectangle in $\mbox{supp}\,B^m$ with $(x_1^m, y_1^m)$ as one of its vertices has to be covered by a MS B-spline $B[\pmb{x}^\ell, \pmb{y}^\ell]\in \mathcal{B}$ in the corner of $\mathcal{R}$ that has larger support in both directions with respect to $B^m$. \qedhere
\end{enumerate}
\end{proof}
\end{prop}
Therefore, in every corner of $\mathcal{R}$ there are at least two MS B-splines, one nested into the other. In particular, when $(p_1,p_2)=(0,0)$, it is not possible to nest a B-spline into another during the LR-mesh building process. This is because it would imply to end a meshline in the middle of a rectangle in order not to split the larger B-spline. Therefore, since Proposition \ref{L1} cannot be satisfied we conclude that the set of MS and LR B-splines of degree $(0,0)$ are linearly independent on any LR-mesh.

Thus, from now on, we assume $(p_1,p_2)\neq (0,0)$.
\begin{cor}\label{fiveMS}
We need at least 5 MS B-splines for a linear dependence relation in $\mathcal{B}^{\mathcal{MS}}(\cM)$.
\begin{proof}
$\mathcal{R}$ has at least four corners and therefore we have 4 different MS B-splines at corners of $\mathcal{R}$ and all of them must be nested into another larger B-spline. The minimal number needed for the nestedness is then 5, considering a MS B-spline whose support coincides with $\mathcal{R}$, containing the 4 MS B-splines at the corners (see Figure \ref{10}).\end{proof} 
\end{cor}
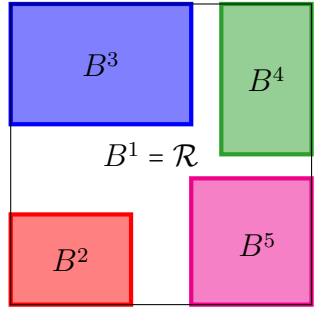
\begin{figure}[h!]\centering
\begin{tikzpicture}[scale=4]
\fill[red, opacity=.5]  (0,0) -- (.4,0) -- (.4,.3) -- (0,.3) -- cycle;
\draw[red, ultra thick] (0,0) -- (.4,0) -- (.4,.3) -- (0,.3) -- cycle;
\fill[magenta, opacity=.5] (.6,0) -- (1,0) -- (1,.42) -- (.6, .42) -- cycle;
\draw[magenta, ultra thick] (.6,0) -- (1,0) -- (1,.42) -- (.6, .42) -- cycle;
\fill[green(ryb), opacity=.5] (.7,.5) -- (1,.5) -- (1,1) -- (.7,1) -- cycle;
\draw[green(ryb), ultra thick] (.7,.5) -- (1,.5) -- (1,1) -- (.7,1) -- cycle;
\fill[blue, opacity=.5] (0,.6) -- (.6,.6) -- (.6,1) -- (0,1) -- cycle;
\draw[blue, ultra thick] (0,.6) -- (.6,.6) -- (.6,1) -- (0,1) -- cycle;
\draw (0,0) -- (1,0) -- (1,1) -- (0,1) -- cycle;
\draw (.66,.5) node[left]{$B^1 = \mathcal{R}$};
\draw (.2,.15) node{$B^2$};
\draw (.3,.8) node{$B^3$};
\draw (.85, .75) node{$B^4$};
\draw (.82,.2) node{$B^5$};
\end{tikzpicture}
\caption{Configuration with 5 MS B-splines satisfying the nested supports condition for linear dependence. With the colors we highlight the support of the B-spline nested into $B^1$. 
}\label{10}
\end{figure}
\noindent
The question now is if 5 minimal support B-splines are enough for a linear dependence relation. From the previous results, we know that if so, the linear dependence region $\mathcal{R}$ has only 4 corners and we have 4 B-splines with supports in the corners of $\mathcal{R}$ and one larger minimal support B-splines with support covering the entire $\mathcal{R}$. The rest of this section is devoted to show that five MS B-splines are not enough. For sake of simplicity, we will keep the notation used in Figure \ref{10}. So $B^1$ will be the larger MS B-spline whose support coincides with $\mathcal{R}$ and $B^2, B^3, B^4, B^5$ are the MS B-splines at the corners ordered clockwise starting from the lower left corner. 

In order to have a linear dependence relation, in every point of $\mathcal{R}$ we must have at least two MS B-splines different from zero. In the following Lemma we present how this fact implies spatial relations of the supports if a linear dependence relation involves only $B^1$, $B^2$, $B^3$, $B^4$ and $B^5$.

\begin{lem}\label{cov}
Suppose only five MS B-splines are in a linear dependence relation on $\mathcal{R}$. Then\begin{enumerate}
\item the supports of $B^2$ and $B^5$ intersect each other as well as the supports of $B^3$ and $B^4$,
\item the supports of $B^2$ and $B^3$ intersect each other as well as the supports of $B^4$ and $B^5$,
\item at least one couple among $\mbox{supp}\,B^2$, $\mbox{supp}\,B^4$ and $\mbox{supp}\,B^3$, $\mbox{supp}\,B^5$ intersect each other.
\end{enumerate}
\begin{proof}
Every point in $\mbox{supp}\,B^1$ must be inside the support of another B-spline in the linear dependence relation, i.e., the supports of $B^2, B^3, B^4, B^5$ must be such that there are no white spots left inside $\mathcal{R}$ in Figure \ref{10}.
\begin{enumerate}
\item We notice that $y^2_1 = y^5_1$ and, by Proposition \ref{L1}, the $y$-widths of the supports of $B^3$ and $B^4$ must be smaller than the $y$-width of $\mathcal{R}$, i.e., $y^{3}_1, y_1^4 > y^{1}_1$. Let $\bar{y} := \min\{y_1^{3}, y_1^{4}\}$. There exists an horizontal band, $[x_1^{2}, x_{p_1+2}^{5}]\times (y_1^{2}, \bar{y})\subset \mathcal{R}$ whose interior cannot be intersected by $\mbox{supp}\,B^{3}$ and $\mbox{supp}\,B^{4}$, see the lower band in Figure \ref{12} (a).  We want to prove that $\mbox{supp}\,B^{2} \cap \mbox{supp}\, B^{5} \neq \emptyset$. If $\mbox{supp}\,B^{2} \cap \mbox{supp}\, B^{5} = \emptyset$ then defined $\bar{\bar{y}}:= \min\{\bar{y}, y_{p_2+2}^2, y_{p_2+2}^5\}$, there would exist a point $(x, y) \in (x_{p_1+2}^{2}, x_1^{5})\times(y_1^{2}, \bar{\bar{y}})$ where none of the four minimal support B-splines in the corners would be different from zero, see Figure \ref{12} (b). $(x,y)$ would be only in the support of $B^1$. This is a contradiction. 
An analogous argument yields that $\mbox{supp}\, B^{3}\cap \mbox{supp}\, B^{4} \neq \emptyset$.
\begin{figure}[h!]\centering
\subfloat[]{
\begin{tikzpicture}[scale=4]
\fill[green(ryb), opacity=.5] (.7,.5) -- (1,.5) -- (1,1) -- (.7,1) -- cycle;
\draw[green(ryb), ultra thick] (.7,.5) -- (1,.5) -- (1,1) -- (.7,1) -- cycle;
\fill[blue, opacity=.5] (0,.6) -- (.6,.6) -- (.6,1) -- (0,1) -- cycle;
\draw[blue, ultra thick] (0,.6) -- (.6,.6) -- (.6,1) -- (0,1) -- cycle;
\fill[orange!50] (0,0) -- (1,0) -- (1,.5) -- (0,.5) -- cycle;
\draw (0,0) -- (1,0) -- (1,1) -- (0,1) -- cycle;
\draw (.3,.8) node{$B^3$};
\draw (.85, .75) node{$B^4$};
\draw[dashed, thick] (0,.5) -- (1,.5);
\draw (1,.5) -- (1.03,.5) node[right]{$\bar{y}$};
\draw (1,0) node[right]{$y_1^2 = y_1^5$};
\draw[white] (-.3,0) -- (-.3,1);
\end{tikzpicture}}\qquad
\subfloat[]{
\begin{tikzpicture}[scale=4]
\fill[orange!50] (.375,0) -- (.625,0) -- (.625, .3) -- (.375,.3) -- cycle;
\fill[red, opacity=.5]  (0,0) -- (.375,0) -- (.375,.3) -- (0,.3) -- cycle;
\draw[red, ultra thick] (0,0) -- (.375,0) -- (.375,.3) -- (0,.3) -- cycle;
\fill[magenta, opacity=.5] (.625,0) -- (1,0) -- (1,.42) -- (.625, .42) -- cycle;
\draw[magenta, ultra thick] (.625,0) -- (1,0) -- (1,.42) -- (.625, .42) -- cycle;
\fill[green(ryb), opacity=.5] (.7,.5) -- (1,.5) -- (1,1) -- (.7,1) -- cycle;
\draw[green(ryb), ultra thick] (.7,.5) -- (1,.5) -- (1,1) -- (.7,1) -- cycle;
\fill[blue, opacity=.5] (0,.6) -- (.6,.6) -- (.6,1) -- (0,1) -- cycle;
\draw[blue, ultra thick] (0,.6) -- (.6,.6) -- (.6,1) -- (0,1) -- cycle;
\draw[dashed] (0,.3) -- (1,.3);
\draw (1,.3) -- (1.03,.3) node[right]{$\bar{\bar{y}}$};
\draw (0,0) -- (1,0) -- (1,1) -- (0,1) -- cycle;
\draw (.2,.15) node{$B^2$};
\draw (.3,.8) node{$B^3$};
\draw (.85, .75) node{$B^4$};
\draw (.82,.2) node{$B^5$};
\draw (1,0) node[right]{$y_1^2 = y_1^5$};
\fill (.502, .2) node[below]{$(x,y)$} circle (.02);
\draw[white] (-.3,0) -- (-.3,1);
\end{tikzpicture}}
\caption{In (a) the horizontal band of $\mathcal{R}$ whose interior is not intersected by $\mbox{supp}\,B^3$ and $\mbox{supp}\,B^4$ is highlighted. In (b), the colored subregion of $\mathcal{R}$ contains a point $(x,y)$ for the proof of the first item of Lemma \ref{cov}.}\label{12}
\end{figure}
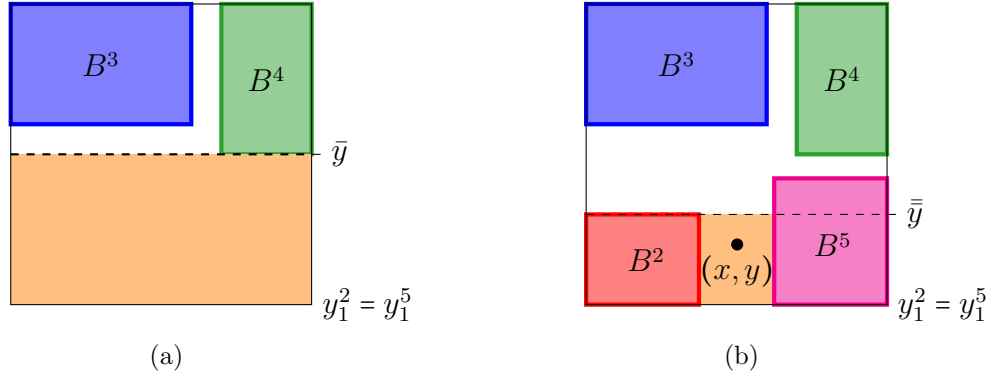
\item Exchanging the axes, we can use the same argument as the previous item.
\item Assume the two couples $B^2, B^4$ and $B^3, B^5$ each do not intersect. Then, since the previous statements are proved, we must have \begin{equation*}\left\{
\begin{array}{l}
x_{p_1+2}^{3} < x_1^{5} \\\\
y_{p_2+2}^{2} < y_1^{4}
\end{array}\right. \quad\mbox{or}\quad\left\{\begin{array}{l}
x_{p_1+2}^{2} < x_1^{4}\\\\
y_{p_2+2}^{5} < y_1^{3}.
\end{array}\right.
\end{equation*}
These two cases, depicted in Figure \ref{dis} (a) and (b) respectively, can be treated in the same way, so we focus only on the first one.
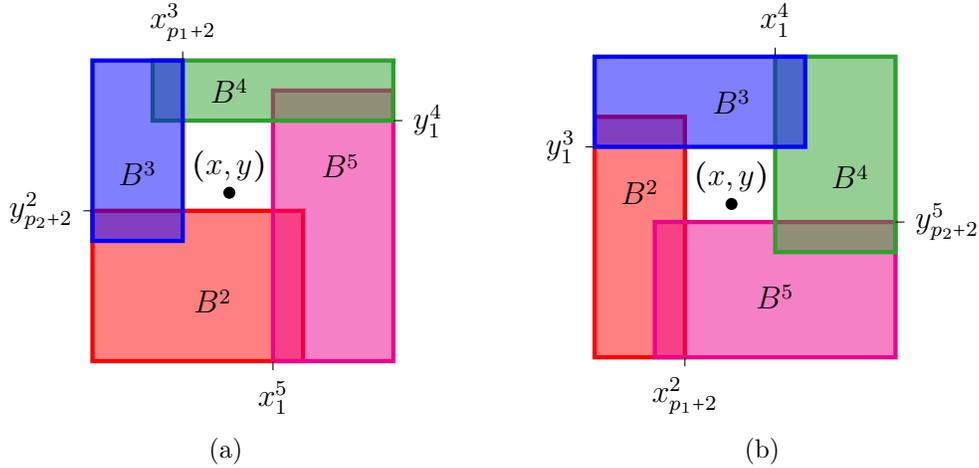
\begin{figure}[h!]\centering
\subfloat[]{
\begin{tikzpicture}[scale=4]
\fill[red, opacity=.5]  (0,0) -- (.7,0) -- (.7,.5) -- (0,.5) -- cycle;
\draw[red, ultra thick] (0,0) -- (.7,0) -- (.7,.5) -- (0,.5) -- cycle;
\fill[magenta, opacity=.5] (.6,0) -- (1,0) -- (1,.9) -- (.6, .9) -- cycle;
\draw[magenta, ultra thick] (.6,0) -- (1,0) -- (1,.9) -- (.6, .9) -- cycle;
\fill[green(ryb), opacity=.5] (.2,.8) -- (1,.8) -- (1,1) -- (.2,1) -- cycle;
\draw[green(ryb), ultra thick] (.2,.8) -- (1,.8) -- (1,1) -- (.2,1) -- cycle;
\fill[blue, opacity=.5] (0,.4) -- (.3,.4) -- (.3,1) -- (0,1) -- cycle;
\draw[blue, ultra thick] (0,.4) -- (.3,.4) -- (.3,1) -- (0,1) -- cycle;
\draw (.4,.2) node{$B^2$};
\draw (.455,.9) node{$B^4$};
\draw (.15, .625) node{$B^3$};
\draw (.825,.65) node{$B^5$};
\draw (-.03,.5) node[left]{$y_{p_2+2}^2$} -- (0,.5);
\draw (1,.8) -- (1.03,.8) node[right]{$y_1^4$};
\draw (.6,-.03) node[below]{$x_1^5$} -- (.6,0);
\draw (.3,1) -- (.3,1.03) node[above]{$x_{p_1+2}^3$};
\fill (.455, .56) node[above]{$(x,y)$} circle (.02);
\end{tikzpicture}
}\qquad
\subfloat[]{
\begin{tikzpicture}[scale=4]
\fill[red, opacity=.5]  (0,0) -- (.3,0) -- (.3,.8) -- (0,.8) -- cycle;
\draw[red, ultra thick] (0,0) -- (.3,0) -- (.3,.8) -- (0,.8) -- cycle;
\fill[magenta, opacity=.5] (.2,0) -- (1,0) -- (1,.45) -- (.2, .45) -- cycle;
\draw[magenta, ultra thick] (.2,0) -- (1,0) -- (1,.45) -- (.2, .45) -- cycle;
\fill[green(ryb), opacity=.5] (.6,.35) -- (1,.35) -- (1,1) -- (.6,1) -- cycle;
\draw[green(ryb), ultra thick] (.6,.35) -- (1,.35) -- (1,1) -- (.6,1) -- cycle;
\fill[blue, opacity=.5] (0,.7) -- (.7,.7) -- (.7,1) -- (0,1) -- cycle;
\draw[blue, ultra thick] (0,.7) -- (.7,.7) -- (.7,1) -- (0,1) -- cycle;
\draw (.15,.55) node{$B^2$};
\draw (.85,.6) node{$B^4$};
\draw (.455, .85) node{$B^3$};
\draw (.6,.2) node{$B^5$};
\draw (-.03,.7) node[left]{$y_1^3$} -- (0,.7);
\draw (1,.45) -- (1.03,.45) node[right]{$y_{p_2+2}^5$};
\draw (.3,-.03) node[below]{$x_{p_1+2}^2$} -- (.3,0);
\draw (.6,1) -- (.6,1.03) node[above]{$x_1^4$};
\fill (.455, .51) node[above]{$(x,y)$} circle (.02);
\end{tikzpicture}
    }
\caption{(a) and (b) are the two possible arrangements of the supports of $B^2, B^3, B^4,B^5$ inside the support of $B^1$, i.e., inside $\mathcal{R}$, when the first two items of Lemma \ref{cov} hold but not the last.}\label{dis}
\end{figure}\\
Consider a point $(x,y) \in (x_{p_1+2}^{3}, x_1^{5})\times (y_{p_2+2}^{2}, y_1^{4})$. 
Since $x \in (x_{p_1+2}^{3}, x_1^{5})$ we have $(x, y) \notin \mbox{supp}\, B_{3}, \mbox{supp}\, B_{5}$. While, since $y \in (y_{p_2+2}^{2}, y_1^{4})$, we have $(x,y) \notin \mbox{supp}\, B^{2}, \mbox{supp}\, B^{4}$. Therefore $(x,y)$ is only in $\mbox{supp}\,B^1 = \mathcal{R}$. This is a contradiction.\qedhere
\end{enumerate}
\end{proof}
\end{lem}
Figure \ref{14} shows possible configurations for the MS B-splines $B^1, B^2, B^3, B^4, B^5$ satisfying Proposition \ref{L1} and Lemma \ref{cov}.
\begin{figure}[h!]\centering
\subfloat[]{
\begin{tikzpicture}[scale=4]
\fill[blue, opacity=.5] (0,.4) -- (.7,.4) -- (.7,1) -- (0,1) -- cycle;
\draw[blue, ultra thick] (0,.4) -- (.7,.4) -- (.7,1) -- (0,1) -- cycle;
\fill[red, opacity=.5]  (0,0) -- (.5,0) -- (.5,.5) -- (0,.5) -- cycle;
\draw[red, ultra thick] (0,0) -- (.5,0) -- (.5,.5) -- (0,.5) -- cycle;
\fill[magenta, opacity=.5] (.3,0) -- (1,0) -- (1,.7) -- (.3, .7) -- cycle;
\draw[magenta, ultra thick] (.3,0) -- (1,0) -- (1,.7) -- (.3, .7) -- cycle;
\fill[green(ryb), opacity=.5] (.4,.6) -- (1,.6) -- (1,1) -- (.4,1) -- cycle;
\draw[green(ryb), ultra thick] (.4,.6) -- (1,.6) -- (1,1) -- (.4,1) -- cycle;
\draw (.15,.2) node{$B^2$};
\draw (.85,.85) node{$B^4$};
\draw (.15, .8) node{$B^3$};
\draw (.8,.2) node{$B^5$};
\end{tikzpicture}}\quad
\subfloat[]{
\begin{tikzpicture}[scale=4]
\fill[blue, opacity=.5] (0,.4) -- (.7,.4) -- (.7,1) -- (0,1) -- cycle;
\draw[blue, ultra thick] (0,.4) -- (.7,.4) -- (.7,1) -- (0,1) -- cycle;
\fill[red, opacity=.5]  (0,0) -- (.5,0) -- (.5,.6) -- (0,.6) -- cycle;
\draw[red, ultra thick] (0,0) -- (.5,0) -- (.5,.6) -- (0,.6) -- cycle;
\fill[magenta, opacity=.5] (.3,0) -- (1,0) -- (1,.7) -- (.3, .7) -- cycle;
\draw[magenta, ultra thick] (.3,0) -- (1,0) -- (1,.7) -- (.3, .7) -- cycle;
\fill[green(ryb), opacity=.5] (.4,.5) -- (1,.5) -- (1,1) -- (.4,1) -- cycle;
\draw[green(ryb), ultra thick] (.4,.5) -- (1,.5) -- (1,1) -- (.4,1) -- cycle;
\draw[green(ryb)!50!blue, ultra thick] (.4,1) -- (.7,1);
\draw[magenta!50!green(ryb), ultra thick] (1,.5) -- (1,.7);
\draw[red!50!blue, ultra thick] (0,.4) -- (0,.6);
\draw[red!50!magenta, ultra thick] (.3,0) -- (.5,0);
\draw (.15,.2) node{$B^2$};
\draw (.85,.85) node{$B^4$};
\draw (.15, .8) node{$B^3$};
\draw (.8,.2) node{$B^5$};
\end{tikzpicture}}\quad
\subfloat[]{
\begin{tikzpicture}[scale=4]
\fill[blue, opacity=.5] (0,.3) -- (.7,.3) -- (.7,1) -- (0,1) -- cycle;
\draw[blue, ultra thick] (0,.3) -- (.7,.3) -- (.7,1) -- (0,1) -- cycle;
\fill[red, opacity=.5]  (0,0) -- (.7,0) -- (.7,.7) -- (0,.7) -- cycle;
\draw[red, ultra thick] (0,0) -- (.7,0) -- (.7,.7) -- (0,.7) -- cycle;
\fill[magenta, opacity=.5] (.3,0) -- (1,0) -- (1,.7) -- (.3, .7) -- cycle;
\draw[magenta, ultra thick] (.3,0) -- (1,0) -- (1,.7) -- (.3, .7) -- cycle;
\fill[green(ryb), opacity=.5] (.3,.3) -- (1,.3) -- (1,1) -- (.3,1) -- cycle;
\draw[green(ryb), ultra thick] (.3,.3) -- (1,.3) -- (1,1) -- (.3,1) -- cycle;
\draw[green(ryb)!50!blue, ultra thick] (.3,1) -- (.7,1);
\draw[blue!50!green(ryb), ultra thick] (.3,.3) -- (.7,.3);
\draw[magenta!50!green(ryb), ultra thick] (1,.3) -- (1,.7);
\draw[magenta!50!green(ryb), ultra thick] (.3,.3) -- (.3,.7);
\draw[red!50!blue, ultra thick] (.7,.3) -- (.7,.7);
\draw[red!50!blue, ultra thick] (0,.3) -- (0,.7);
\draw[red!50!magenta, ultra thick] (.3,0) -- (.7,0);
\draw[red!50!magenta, ultra thick] (.3,.7) -- (.7,.7);
\draw (.15,.15) node{$B^2$};
\draw (.85,.85) node{$B^4$};
\draw (.15, .85) node{$B^3$};
\draw (.85,.15) node{$B^5$};
\end{tikzpicture}}
\caption{Arrangements of the four B-splines at corners of $\mbox{supp} B^1$ satisfying Lemma \ref{cov}. In (a) only the supports $B^3$ and $B^5$ intersect each other, while on the other diagonal of $\mathcal{R}$, the supports of $B^2$ and $B^3$ do not intersect. In (b) and (c) both the pairs at opposing corners of $\mathcal{R}$ intersect each other. In particular, in (c) $B^2$ is as tall as $B^5$, $B^3$ is tall as $B^4$, $B^2$ is as wide as $B^3$ and $B^4$ is as wide as $B^5$.}\label{14}
\end{figure}
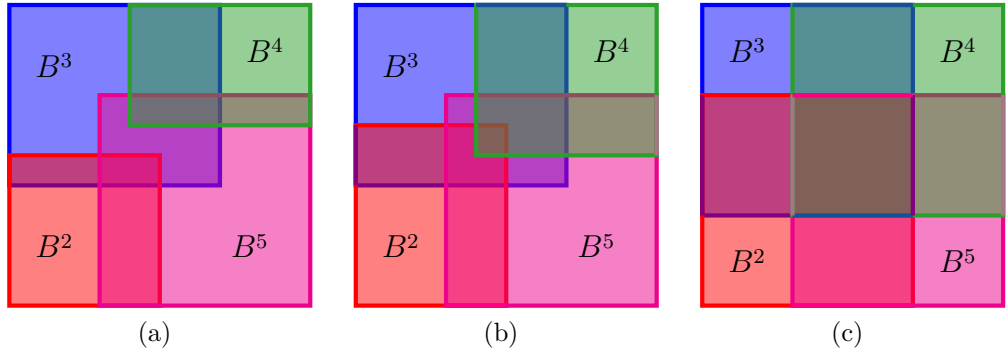
The next results provide further information on the B-splines arrangement inside the region $\mathcal{R}=\mbox{supp}\,B^1$, in the presence of linear dependence. 
Such results hold regardless of the number of MS B-splines involved. Assume again that $\cB = \{B[\pmb{x}^j, \pmb{y}^j]\}_{j=1}^n$, $\mathcal{B} \subseteq \cB^{\mathcal{MS}}(\cM)$, for some $n$, is a set of linearly dependent MS B-splines,
\begin{dfn}
A vertex $(\bar{x}, \bar{y})$ in $\mathcal{R}$ is called \textbf{relevant} if it corresponds to a pair of knots in at least one MS B-spline in the linear dependence relation, see Figure \ref{15}.\\
A meshline $\gamma$ is called \textbf{relevant} if it is part of a split of a MS B-spline in the linear dependence relation.
\end{dfn}
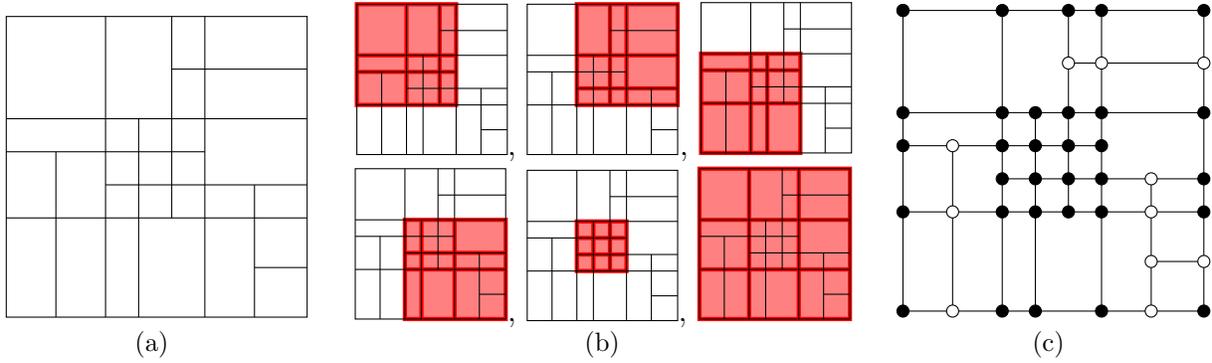
\begin{figure}[h!]\centering
\subfloat[]{
\begin{tikzpicture}[scale=4]
\draw (0,0) -- (1,0) -- (1,1) -- (0,1) -- cycle;
\draw (.33,0) -- (.33,1);
\draw (.66,0) -- (.66,1);
\draw (0,.33) -- (1,.33);
\draw (0,.66) -- (1,.66);
\draw (.44,0) -- (.44,.66);
\draw (.55,.33) -- (.55,1);
\draw (0,.55) -- (.66,.55);
\draw (.33,.44) -- (1,.44);
\draw (.165,0) -- (.165,.55);
\draw (.55,.825) -- (1,.825);
\draw (.825,.165) -- (1,.165);
\draw (.825,0) -- (.825,.44);
\end{tikzpicture}}\quad
\subfloat[]{\centering
\begin{minipage}{6cm} \vspace{-4cm}
\begin{tabular}{c}
\begin{tikzpicture}[scale=2]
\fill[red, opacity=.5] (0,.33) -- (.66,.33) -- (.66,1) -- (0,1) -- cycle;
\draw[red, ultra thick] (0,.33) -- (.66,.33) -- (.66,1) -- (0,1) -- cycle; 
\draw[red, ultra thick] (0,.55) -- (.66,.55);
\draw[red, ultra thick] (0,.66) -- (.66,.66);
\draw[red, ultra thick] (.33,.33) -- (.33,1);
\draw[red, ultra thick] (.55,.33) -- (.55,1);
\draw (0,0) -- (1,0) -- (1,1) -- (0,1) -- cycle;
\draw (.33,0) -- (.33,1);
\draw (.66,0) -- (.66,1);
\draw (0,.33) -- (1,.33);
\draw (0,.66) -- (1,.66);
\draw (.44,0) -- (.44,.66);
\draw (.55,.33) -- (.55,1);
\draw (0,.55) -- (.66,.55);
\draw (.33,.44) -- (1,.44);
\draw (.165,0) -- (.165,.55);
\draw (.55,.825) -- (1,.825);
\draw (.825,.165) -- (1,.165);
\draw (.825,0) -- (.825,.44);
\end{tikzpicture},
\begin{tikzpicture}[scale=2]
\fill[red, opacity=.5] (.33,.33) -- (1,.33) -- (1,1) -- (.33,1) -- cycle;
\draw[red, ultra thick] (.33,.33) -- (1,.33) -- (1,1) -- (.33,1) -- cycle;
\draw[red, ultra thick] (.33,.44) -- (1,.44);
\draw[red, ultra thick] (.33,.66) -- (1,.66);
\draw[red, ultra thick] (.55,.33) -- (.55,1);
\draw[red, ultra thick] (.66,.33) -- (.66,1); 
\draw (0,0) -- (1,0) -- (1,1) -- (0,1) -- cycle;
\draw (.33,0) -- (.33,1);
\draw (.66,0) -- (.66,1);
\draw (0,.33) -- (1,.33);
\draw (0,.66) -- (1,.66);
\draw (.44,0) -- (.44,.66);
\draw (.55,.33) -- (.55,1);
\draw (0,.55) -- (.66,.55);
\draw (.33,.44) -- (1,.44);
\draw (.165,0) -- (.165,.55);
\draw (.55,.825) -- (1,.825);
\draw (.825,.165) -- (1,.165);
\draw (.825,0) -- (.825,.44);
\end{tikzpicture},
\begin{tikzpicture}[scale=2]
\fill[red, opacity=.5] (0,0) -- (.66,0) -- (.66,.66) -- (0,.66) -- cycle;
\draw[red, ultra thick] (0,0) -- (.66,0) -- (.66,.66) -- (0,.66) -- cycle;
\draw[red, ultra thick] (0,.33) -- (.66,.33);
\draw[red, ultra thick] (0,.55) -- (.66,.55);
\draw[red, ultra thick] (.33,0) -- (.33,.66);
\draw[red, ultra thick] (.44,0) -- (.44,.66); 
\draw (0,0) -- (1,0) -- (1,1) -- (0,1) -- cycle;
\draw (.33,0) -- (.33,1);
\draw (.66,0) -- (.66,1);
\draw (0,.33) -- (1,.33);
\draw (0,.66) -- (1,.66);
\draw (.44,0) -- (.44,.66);
\draw (.55,.33) -- (.55,1);
\draw (0,.55) -- (.66,.55);
\draw (.33,.44) -- (1,.44);
\draw (.165,0) -- (.165,.55);
\draw (.55,.825) -- (1,.825);
\draw (.825,.165) -- (1,.165);
\draw (.825,0) -- (.825,.44);
\end{tikzpicture}\\
\begin{tikzpicture}[scale=2]
\fill[red, opacity=.5] (.33,0) -- (1,0) -- (1,.66) -- (.33,.66) -- cycle;
\draw[red, ultra thick] (.33,0) -- (1,0) -- (1,.66) -- (.33,.66) -- cycle;
\draw[red, ultra thick] (.33,.33) -- (1,.33);
\draw[red, ultra thick] (.33,.44) -- (1,.44);
\draw[red, ultra thick] (.44,0) -- (.44,.66);
\draw[red, ultra thick] (.66,0) -- (.66,.66);
\draw (0,0) -- (1,0) -- (1,1) -- (0,1) -- cycle;
\draw (.33,0) -- (.33,1);
\draw (.66,0) -- (.66,1);
\draw (0,.33) -- (1,.33);
\draw (0,.66) -- (1,.66);
\draw (.44,0) -- (.44,.66);
\draw (.55,.33) -- (.55,1);
\draw (0,.55) -- (.66,.55);
\draw (.33,.44) -- (1,.44);
\draw (.165,0) -- (.165,.55);
\draw (.55,.825) -- (1,.825);
\draw (.825,.165) -- (1,.165);
\draw (.825,0) -- (.825,.44);
\end{tikzpicture},
\begin{tikzpicture}[scale=2]
\fill[red, opacity=.5] (.33,.33) -- (.66,.33) -- (.66,.66) -- (.33,.66) -- cycle;
\draw[red, ultra thick] (.33,.33) -- (.66,.33) -- (.66,.66) -- (.33,.66) -- cycle;
\draw[red, ultra thick] (.33,.44) -- (.66,.44);
\draw[red, ultra thick] (.33,.55) -- (.66,.55);
\draw[red, ultra thick] (.44,.33) -- (.44,.66);
\draw[red, ultra thick] (.55,.33) -- (.55,.66);
\draw (0,0) -- (1,0) -- (1,1) -- (0,1) -- cycle;
\draw (.33,0) -- (.33,1);
\draw (.66,0) -- (.66,1);
\draw (0,.33) -- (1,.33);
\draw (0,.66) -- (1,.66);
\draw (.44,0) -- (.44,.66);
\draw (.55,.33) -- (.55,1);
\draw (0,.55) -- (.66,.55);
\draw (.33,.44) -- (1,.44);
\draw (.165,0) -- (.165,.55);
\draw (.55,.825) -- (1,.825);
\draw (.825,.165) -- (1,.165);
\draw (.825,0) -- (.825,.44);
\end{tikzpicture},
\begin{tikzpicture}[scale=2]
\fill[red, opacity=.5] (0,0) -- (1,0) -- (1,1) -- (0,1) -- cycle;
\draw[red, ultra thick] (0,0) -- (1,0) -- (1,1) -- (0,1) -- cycle;
\draw[red, ultra thick] (0,.33) -- (1,.33);
\draw[red, ultra thick] (0,.66) -- (1,.66);
\draw[red, ultra thick] (.33,0) -- (.33,1);
\draw[red, ultra thick] (.66,0) -- (.66,1);
\draw (0,0) -- (1,0) -- (1,1) -- (0,1) -- cycle;
\draw (.33,0) -- (.33,1);
\draw (.66,0) -- (.66,1);
\draw (0,.33) -- (1,.33);
\draw (0,.66) -- (1,.66);
\draw (.44,0) -- (.44,.66);
\draw (.55,.33) -- (.55,1);
\draw (0,.55) -- (.66,.55);
\draw (.33,.44) -- (1,.44);
\draw (.165,0) -- (.165,.55);
\draw (.55,.825) -- (1,.825);
\draw (.825,.165) -- (1,.165);
\draw (.825,0) -- (.825,.44);
\end{tikzpicture}
\end{tabular}
\vspace{-.25cm}
\end{minipage}\hspace{1cm}}~
\subfloat[]{
\begin{tikzpicture}[scale=4]\vspace{.25cm}
\draw (0,0) -- (1,0) -- (1,1) -- (0,1) -- cycle;
\draw (.33,0) -- (.33,1);
\draw (.66,0) -- (.66,1);
\draw (0,.33) -- (1,.33);
\draw (0,.66) -- (1,.66);
\draw (.44,0) -- (.44,.66);
\draw (.55,.33) -- (.55,1);
\draw (0,.55) -- (.66,.55);
\draw (.33,.44) -- (1,.44);
\draw (.165,0) -- (.165,.55);
\draw (.55,.825) -- (1,.825);
\draw (.825,.165) -- (1,.165);
\draw (.825,0) -- (.825,.44);
\filldraw[draw=black,fill=black] (0,0) circle (.02);
\filldraw[draw=black,fill=black] (.33,0) circle (.02);
\filldraw[draw=black,fill=black] (.44,0) circle (.02);
\filldraw[draw=black,fill=black] (.66,0) circle (.02);
\filldraw[draw=black,fill=black] (1,0) circle (.02);
\filldraw[draw=black,fill=black] (0,.33) circle (.02);
\filldraw[draw=black,fill=black] (.33,.33) circle (.02);
\filldraw[draw=black,fill=black] (.44,.33) circle (.02);
\filldraw[draw=black,fill=black] (.55,.33) circle (.02);
\filldraw[draw=black,fill=black] (.66,.33) circle (.02);
\filldraw[draw=black,fill=black] (1,.33) circle (.02);
\filldraw[draw=black,fill=black] (.33,.44) circle (.02);
\filldraw[draw=black,fill=black] (.44,.44) circle (.02);
\filldraw[draw=black,fill=black] (.55,.44) circle (.02);
\filldraw[draw=black,fill=black] (.66,.44) circle (.02);
\filldraw[draw=black,fill=black] (1,.44) circle (.02);
\filldraw[draw=black,fill=black] (0,.55) circle (.02);
\filldraw[draw=black,fill=black] (.33,.55) circle (.02);
\filldraw[draw=black,fill=black] (.44,.55) circle (.02);
\filldraw[draw=black,fill=black] (.55,.55) circle (.02);
\filldraw[draw=black,fill=black] (.66,.55) circle (.02);
\filldraw[draw=black,fill=black] (0,.66) circle (.02);
\filldraw[draw=black,fill=black] (.33,.66) circle (.02);
\filldraw[draw=black,fill=black] (.44,.66) circle (.02);
\filldraw[draw=black,fill=black] (.55,.66) circle (.02);
\filldraw[draw=black,fill=black] (.66,.66) circle (.02);
\filldraw[draw=black,fill=black] (1,.66) circle (.02);
\filldraw[draw=black,fill=black] (0,1) circle (.02);
\filldraw[draw=black,fill=black] (.33,1) circle (.02);
\filldraw[draw=black,fill=black] (.55,1) circle (.02);
\filldraw[draw=black,fill=black] (.66,1) circle (.02);
\filldraw[draw=black,fill=black] (1,1) circle (.02);
\filldraw[draw=black,fill=black] (0,1) circle (.02);
\filldraw[draw=black,fill=white] (.825,0) circle (.02);
\filldraw[draw=black,fill=white] (.165,.33) circle (.02);
\filldraw[draw=black,fill=white] (.825,.165) circle (.02);
\filldraw[draw=black,fill=white] (1,.165) circle (.02);
\filldraw[draw=black,fill=white] (.825,.33) circle (.02);
\filldraw[draw=black,fill=white] (.825,.44) circle (.02);
\filldraw[draw=black,fill=white] (.55,.825) circle (.02);
\filldraw[draw=black,fill=white] (.66,.825) circle (.02);
\filldraw[draw=black,fill=white] (1,.825) circle (.02);
\filldraw[draw=black,fill=white] (.165,.55) circle (.02);
\filldraw[draw=black,fill=white] (.165,0) circle (.02);
\end{tikzpicture}}
\caption{Consider the mesh on the region $\mathcal{R}$ depicted in (a).
Every meshline has multiplicity 1 and consider bidegree (2,2). In (b) we see the MS B-splines on $\mathcal{R}$. We will prove they are linearly dependent in Example \ref{exMS}. In (c) we see the relevant vertices in black and the non-relevant vertices in white.}\label{15}
\end{figure}
\begin{lem}\label{lem1}
Any relevant vertex in $\mathcal{R}$ is the intersection of orthogonal relevant meshlines.
\begin{proof}
Let $(\bar{x}, \bar{y})$ be a relevant vertex in $\mathcal{R}$. Then it corresponds to a pair of knots of $B[\pmb{x}^j, \pmb{y}^j]$ for some $j$. In particular, $(\bar{x}, \bar{y})$ is in the orthogonal splits $[x_1^j, x_{p_1+2}^j] \times \{\bar{y}\}$ and $\{\bar{x}\}\times [y_1^j, y_{p_2+2}^j]$. Therefore there must exist at least 2 relevant orthogonal meshlines 
contained in such splits intersecting in $(\bar{x}, \bar{y})$.
\end{proof}
\end{lem}
\begin{prop}\label{mlsha}
Any relevant meshline is shared by at least two splits of two MS B-splines in the linear dependence relation.
\begin{proof}
Let $\cB = \{B[\pmb{x}^j, \pmb{y}^j]\}_{j=1}^n$ be the set of linearly dependent MS B-splines. Let $\gamma$ be a $k$-meshline of any split of $B[\pmb{x}^1, \pmb{y}^1]$. Assume that $\gamma$ is not part of a split of any other MS B-spline in $\cB$. We know that $B[\pmb{x}^1, \pmb{y}^1]$ is $C^{p_k-\mu(\gamma)}$-continuous on $\gamma$.

The linear dependence relation in $\cB$, $$
\alpha_1B[\pmb{x}^1,\pmb{y}^1](x,y) + \sum_{j=2}^n \alpha_j B[\pmb{x}^j, \pmb{y}^j](x,y)=0 \quad\forall \,(x,y) \in \mathcal{R},
$$
can be rewritten expressing $B[\pmb{x}^1, \pmb{y}^1]$ in terms of the others 
\begin{equation}\label{bexp}
B[\pmb{x}^1, \pmb{y}^1](x,y) = -\frac{1}{\alpha_1}\cdot \sum_{j=2}^n \alpha_j B[\pmb{x}^j, \pmb{y}^j](x,y) \quad\forall \,(x,y) \in \mathcal{R}
\end{equation}
because $\alpha_j \neq 0$ for every $j=1, \ldots, n$.
Consider now $(x,y) \in \gamma \subset \mathcal{R}$. Since $\gamma$ is not a split of any MS B-spline $B[\pmb{x}^j, \pmb{y}^j]$ in $\cB$ with $j\geq 2$, the right-hand side is a function $C^\infty$-continuous on $\gamma$ while the left-hand side is only $C^{p_k-\mu(\gamma)}$-continuous on $\gamma$, which is a contradiction. 
\end{proof}
\end{prop}
\begin{cor}\label{Tver}
Any relevant T-vertex corresponds to a pair of knots shared by at least two MS B-splines active in the linear dependence relation.
\begin{proof}
Let $(\bar{x}, \bar{y})$ be a relevant T-vertex as in Figure \ref{16}. The other three possible cases of T-vertex can be treated similarly.
\begin{figure}[h!]\centering
\begin{tikzpicture}[scale=4]
\draw (0,0) -- (0,1);
\draw (0,.5) --node[below]{$\gamma$} (1,.5);
\filldraw[draw=black,fill=red] (0,.5) node[left]{\textcolor{black}{$(\bar{x}, \bar{y})$}} circle (.02);
\end{tikzpicture}
\caption{The T-vertex used in the proof of Corollary \ref{Tver}.}\label{16}
\end{figure}
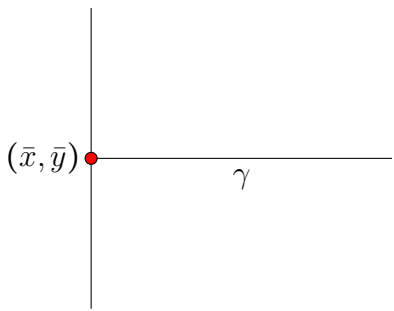\\
Since $(\bar{x}, \bar{y})$ is relevant, $\gamma$ must be relevant from Lemma \ref{lem1}. By  Proposition \ref{mlsha}, $\gamma$ is shared by at least two minimal support B-splines active in the linear dependence relation, say $B[\pmb{x}^1, \pmb{y}^1], B[\pmb{x}^2, \pmb{y}^2]$. This means there are two knots $y_r^1 \in \pmb{y}^1$ and $y_s^2\in \pmb{y}^2$ such that $y_r^1= \bar{y} = y_s^2$ and $[x_1^1, x_{p_1+2}^1]\times \{\bar{y}\}$, $[x_1^2, x_{p_1+2}^2]\times \{\bar{y}\}$ are splits of $B[\pmb{x}^1, \pmb{y}^1]$ and $B[\pmb{x}^2, \pmb{y}^2]$ respectively containing $\gamma$. Since $(\bar{x}, \bar{y})$ is a T-vertex, it ends the splits containing $\gamma$. This means that $(x_1^1, y_r^1) = (\bar{x}, \bar{y}) = (x_1^2, y_s^2)$, i.e., $(\bar{x}, \bar{y})$ is a pair of knots shared by $B[\pmb{x}^1, \pmb{y}^1]$ and $B_2[\pmb{x}^2, \pmb{y}^2]$.
\end{proof}
\end{cor}

\begin{oss} 
In section \ref{peel} we will see that one can use the previous result to improve the Peeling algorithm \cite[Algorithm 6.3]{dokken2013polynomial}, a tool to check if the LR B-splines considered are linearly independent. 
\end{oss}
\begin{dfn}
Any T-vertex $\pmb{v}$ in an LR-mesh is composed of two colinear meshlines and another meshline $\gamma$ orthogonal to them. We assign an orientation to these vertices in the following way. We say that the T-vertex $\pmb{v}$ is \textbf{downward} if $\gamma$ is below $\pmb{v}$,  it is \textbf{upward} if $\gamma$ is above $\pmb{v}$, \textbf{rightward} if $\gamma$ is on the right of $\pmb{v}$ and \textbf{leftward} if $\gamma$ is on the left of $\pmb{v}$.
\end{dfn}
It might happen that a relevant vertex $\pmb{v}$ in $\mathring{\mathcal{R}}$ is a cross-vertex but one meshline ending in $\pmb{v}$ is not relevant. It means that $\pmb{v}$ behaves as a T-vertex for the B-splines in the linear dependence relation. Therefore, we extend the definition of relevant T-vertex and of its orientation also to these vertices in $\mathring{\mathcal{R}}$.

\begin{thm}\label{lemT}
Assume five MS B-splines are linearly dependent inside the region $\mathcal{R}$. Then there are at least 4 relevant T-vertices in $\mathring{\mathcal{R}}$, one per orientation.
\begin{proof}
For the sake of simplicity and without loss of generality, we can assume there are only relevant meshlines in $\mathcal{R}$.
 
Referring to any of the examples in Figure \ref{14}, let us consider the vertical splits of $B^2$ and $B^5$ in the interior of the support of $B^1$, i.e., in $\mathring{\mathcal{R}}$.  

In order to find the minimal number of relevant T-vertices in $\mathring{\mathcal{R}}$, we assume that the parameter values of such vertical splits are the same for $B^2$ and $B^5$.

We assume the same for $B^3$ and $B^4$: the vertical splits of $B^4$ are included into the vertical splits of $B^3$. 

Now, in $\mathring{\mathcal{R}}$ there are $p_1+1$ vertical splits of $B^5$ and $p_1+1$ for $B^3$. 
If an end vertex of one such split of $B^3$ or $B^5$ corresponds to a relevant cross-vertex, it is contained in a split traversing the entire region $\mathcal{R}$, that is, it is contained in a vertical split of $B^1$. 

There are $p_1$ vertical splits of $B^1$ in $\mathring{\mathcal{R}}$. Therefore, at most $p_1$ vertical splits in $\mathring{\mathcal{R}}$ of $B^3$ and $B^5$ can end with a relevant cross-vertex. 

Thus there are $p_1 +1 - p_1 = 1$ relevant vertices of $B^5$ left on the upper edge of $\mbox{supp}\, B^5$ inside $\mathcal{R}$  that cannot be cross-vertices. The same holds for the relevant vertices in $B^3$.

Finally, applying the same argument to the horizontal splits of $B^3$ and $B^5$ we complete the proof.
\end{proof}
\end{thm}
Proposition \ref{lemT} holds also if the number of B-splines involved in the linear dependence relation is larger because of the necessary presence of nested B-splines at the corners.
 
In order to carry out the proof of the next Proposition \ref{sixMS}, we need the following 
\begin{dfn}
Given an LR-mesh $\cM$, let $\gamma$ be a $k$-split in $\cM$ for some $k\in \{1,2\}$. For instance, assume $k=1$ and $\gamma = \{a\}\times[y_{st}, y_{end}]$. Let $F:\RR^2 \to \RR$ be a spline function in $\SSS_{\pmb{p}}(\cM)$. $F$ is a piecewise polynomial and therefore, for sufficiently small $\eps >0$, the functions $F^+ = F_{|(a,a+\eps) \times \RR}$ and $F^- = F_{|(a-\eps,a) \times \RR}$ are polynomials in $x$ (but splines in $y$), i.e., \begin{equation*}
F^+ = \sum_{i=0}^{p_1} f_i^+(y)\cdot (x-a)^i,\qquad
F^- = \sum_{i=0}^{p_1} f_i^-(y)\cdot (x-a)^i
\end{equation*}
for $f_i^+, f_i^-$ univariate spline functions. Then we can extend the expression of $F^+$ and $F^-$ to $\RR^2$. We define the jump function of $F$ with respect to $\gamma$ as $J(F)(x,y) = F^+ - F^-$.
\end{dfn}
\begin{oss}$\left.\right.$
\begin{itemize}
\item If $\gamma$ is not in a split traversing the support of $F$ and is not on its boundary, then $F$ is $C^\infty(\gamma)$ and in particular $F^+ = F^-$ so that $J(F)(x,y) = 0$.
\item When $F$ is a bivariate B-spline, $F = B[\pmb{x}, \pmb{y}]$ and $\gamma$ corresponds to a knot in $\pmb{x}$, $x_i = a$ for some $i$, $\gamma= \{x_i\}\times [y_1, y_{p_2+2}]$, then $$
J(B)(x,y) = J'(B[\pmb{x}])(x) \cdot B[\pmb{y}](y)$$ where $J'(B[\pmb{x}])(x)$ is a polynomial of the form: $$J'(B[\pmb{x}])(x)= \sum_{i=p_1-\mu(\gamma) +1}^{p_1} c_i (x-a)^i.
$$
\item Let $c_1, c_2$ be real numbers and $F_1, F_2$ be spline functions. Then $J(c_1F_1 + c_2F_2)(x,y) = c_1J(F_1)(x,y) + c_2J(F_2)(x,y)$. 
\end{itemize}
\end{oss}
\begin{prop}\label{sixMS}
We need at least 6 minimal support B-splines for a linear dependence relation in $\mathcal{R}$ for any bidegree $(p_1, p_2) \neq (0,0)$.
\begin{proof}
Referring to any configuration in Figure \ref{14}, consider a relevant T-vertex $\pmb{v}$ in $B^5$. By Corollary \ref{Tver}, it has to be shared with at least another MS B-spline. It cannot be shared with $B^2$ if $B^2$ is shorter than $B^5$, and of course it cannot be shared with $B^3$ or $B^4$ because it would not be a T-vertex. Then we have two cases:
\begin{itemize}
\item There exists a new MS B-spline in the linear dependence relation with support in the $y$-direction covering the space between the support of $B^5$ and $B^2$ and having $\pmb{v}$ as pair of knots, or
\item $B^2$ is as tall as $B^5$. 
\end{itemize}
In the first case we have finished the proof. Let us assume then that $B^2$ is as tall as $B^5$. 
Applying the same procedure to the other relevant T-vertices, we either have at least a new MS B-spline in the linear dependence relation, or it must be that $B^4$ is as tall as $B^3$, $B^2$ is as wide as $B^3$ and $B^4$ is as wide as $B^5$. In the first case we have completed the proof. 

In the second, if no other B-splines are involved, we can write $B^1$ in terms of $B^2, B^3, B^4, B^5$:\begin{equation}\label{lind}
B^1(x,y) = \alpha_2B^2(x,y)+ \alpha_3B^3(x,y) + \alpha_4B^4(x,y) + \alpha_5B^5(x,y)\qquad\mbox{with }\alpha_j \neq 0.
\end{equation}
Now, consider any T-vertex downward, corresponding to a $1$-split $\gamma$ in $B^2$ and $B^5$. 
The jump functions of $B^2 = B[\pmb{x}^2, \pmb{y}^2]$ and $B^5[\pmb{x}^5, \pmb{y}^5]$ corresponding to $\gamma$, 
in order to represent $B^1$ as in equation \eqref{lind}, must satisfy \begin{equation}\label{eqjum}
\alpha_2J'(B[\pmb{x}^2])(x)\cdot B[\pmb{y}^2](y) = - \alpha_5J'(B[\pmb{x}^5])(x)\cdot B[\pmb{y}^5](y)
\end{equation}
because $B^1$ is smooth on $\gamma$ and there are no other MS B-splines in the linear dependence relation with less regularity in the $x$-direction on $\gamma$. 
However, the knots of $\pmb{y}^2$ and $\pmb{y}^5$ are different because of the presence of T-vertices leftward and rightward, and equation \eqref{eqjum} is impossible to achieve because $B[\pmb{y}^2]$ and $B[\pmb{y}^5]$ are defined on different knots and cannot be proportional everywhere.
\end{proof}
\end{prop}
\begin{prop}\label{noLR}
In a linear dependence relation with six MS B-splines, the sixth MS B-spline, $B^6$, is not an LR B-spline. 
\begin{proof}
If $B^6$ is an LR B-spline it has been obtained through knot insertion from a LR B-spline in a coarser mesh. When the knot insertion is applied the size of the refined B-splines is smaller only in the direction where the knot has been inserted. Therefore, for $B^6$, in order to be an LR B-spline and be in the linear dependence relation there would exists another B-spline among $B^2, B^3, B^4$ and $B^5$ whose support is either as tall or as wide as the support of $B^6$ and intersects with the support of $B^6$. Assume we are in the first case of the proof of Proposition \ref{sixMS} and there are exactly six MS B-splines in linear dependence. Then there are 4 relevant T-vertices in $\mathring{\mathcal{R}}$ shared with $B^6$ and indentifying the edges of $\mbox{supp}\,B^6$. Therefore, $\mbox{supp}\,B^6 \subseteq \mathring{\mathcal{R}}$ and cannot be the same as the size of any of $B^2,B^3,B^4,B^5$ in any direction.

In the second case of the proof of Proposition \ref{sixMS}, if $B^6$ is an LR B-spline, we can assume that $B^6$ is as tall as $B^2$ and $B^5$ (the other cases can be treated similarly). Then either $\mbox{supp}\,B^6 \subseteq \mbox{supp}\,B^2$ or $\mbox{supp}\,B^6 \subseteq \mbox{supp}\,B^5$. In both cases, there would exists a vertical split of $B^6$ that traverses the support of $B^2$ or $B^5$ without being a split of it. This is impossible for the minimality of $\mbox{supp}\,B^2$ and $\mbox{supp}\,B^5$. 
\end{proof}
\end{prop}
\begin{ex}[Linear dependence relation in $\mathcal{B}^{\mathcal{MS}}(\cM)$ with 6 minimal support B-splines]\label{exMS}
In this example we prove that 6 active MS B-splines are enough for a linear dependence relation in an LR-mesh. As we will see, the arrangement of the 6 MS B-splines inside $\mathcal{R}$ can be reproduced for any bidegree $\pmb{p}=(p_1,p_2)\neq (0,0)$. However, for the sake of simplicity, we start considering bidegree $\pmb{p}=(2,2)$. The LR-mesh $\cM$ considered is depicted in Figure \ref{7} (a) and the support of the MS B-splines defined on it are represented in Figure \ref{7} (b) and (c), and they are 10. However, using the dimension increasing formula in Theorem \ref{diminc}, knowing that\begin{itemize}
\item the dimension of the MS B-spline span on the underlying tensor mesh is 3,
\item the insertion of the two horizontal splits, $p_2+2=4$ long, increases by 1 the dimension,
\item the insertion of the two vertical splits, $p_1+3=5$ long, increases by 2 the dimension,
\end{itemize}
we easily compute the dimension of the spline space on $\cM$, $$\dim\,\SSS_{(2,2)}(\cM) = 3+1+1+2+2 = 9.$$
Moreover, the construction of $\cM$ went LR-wise, and so MS-wise, hand-in-hand. 
Therefore, we can conclude there is a linear dependence relation in the set $\cB^{\mathcal{MS}}(\cM)$. The elements of the set satisfying the necessary conditions, given in this section, to be in linear dependence are the six MS B-splines whose support is depicted in Figure \ref{22}. Finally, notice that the 9 LR B-splines on $\cM$, reported in Figure \ref{7} (b), are still linearly independent and span the spline space on $\cM$.\\
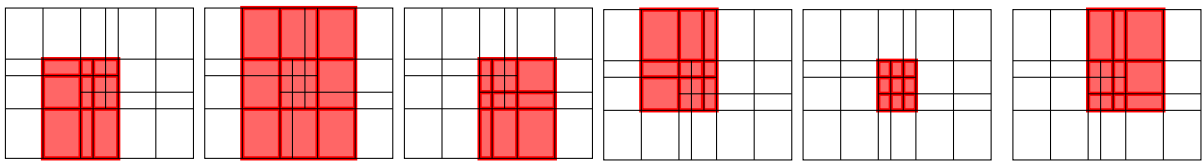
\begin{figure}[h!]\centering
\subfloat{
\begin{tikzpicture}[scale=2]
\filldraw[draw=red, ultra thick, fill=red!60] (.25,0) -- (.75, 0) -- (.75,.66) -- (.25,.66) -- cycle;
\draw[red, ultra thick] (.25,.33) -- (.75, .33);
\draw[red, ultra thick] (.25,.55) -- (.75, .55);
\draw[red, ultra thick] (.5,0) -- (.5, .66);
\draw[red, ultra thick] (.583,0) -- (.583, .66);
\draw (0,0) -- (1.25,0) -- (1.25,1) -- (0,1) -- cycle;
\draw (0,.55) -- (.75, .55);
\draw (.5, .44) -- (1.25, .44);
\draw (.583, 0) -- (.583, .66);
\draw (.666, .33) -- (.666, 1);
\draw (0,.33) -- (1.25, .33);
\draw (0,.66) -- (1.25, .66);
\draw (.25, 0) -- (.25, 1);
\draw (.5, 0) -- (.5, 1);
\draw (.75, 0) -- (.75, 1);
\draw (1, 0) -- (1, 1);
\end{tikzpicture}}
\subfloat{
\begin{tikzpicture}[scale=2]
\filldraw[draw=red, ultra thick, fill=red!60] (.25,0) -- (1, 0) -- (1,1) -- (.25,1) -- cycle;
\draw[red, ultra thick] (.25,.33) -- (1, .33);
\draw[red, ultra thick] (.25,.66) -- (1, .66);
\draw[red, ultra thick] (.5,0) -- (.5, 1);
\draw[red, ultra thick] (.75,0) -- (.75, 1);
\draw (0,0) -- (1.25,0) -- (1.25,1) -- (0,1) -- cycle;
\draw (0,.55) -- (.75, .55);
\draw (.5, .44) -- (1.25, .44);
\draw (.583, 0) -- (.583, .66);
\draw (.666, .33) -- (.666, 1);
\draw (0,.33) -- (1.25, .33);
\draw (0,.66) -- (1.25, .66);
\draw (.25, 0) -- (.25, 1);
\draw (.5, 0) -- (.5, 1);
\draw (.75, 0) -- (.75, 1);
\draw (1, 0) -- (1, 1);
\end{tikzpicture}}
\subfloat{
\begin{tikzpicture}[scale=2]
\filldraw[draw=red, ultra thick, fill=red!60] (.5,0) -- (1, 0) -- (1,.66) -- (.5,.66) -- cycle;
\draw[red, ultra thick] (.5,.33) -- (1, .33);
\draw[red, ultra thick] (.5,.44) -- (1, .44);
\draw[red, ultra thick] (.583,0) -- (.583, .66);
\draw[red, ultra thick] (.75,0) -- (.75, .66);
\draw (0,0) -- (1.25,0) -- (1.25,1) -- (0,1) -- cycle;
\draw (0,.55) -- (.75, .55);
\draw (.5, .44) -- (1.25, .44);
\draw (.583, 0) -- (.583, .66);
\draw (.666, .33) -- (.666, 1);
\draw (0,.33) -- (1.25, .33);
\draw (0,.66) -- (1.25, .66);
\draw (.25, 0) -- (.25, 1);
\draw (.5, 0) -- (.5, 1);
\draw (.75, 0) -- (.75, 1);
\draw (1, 0) -- (1, 1);
\end{tikzpicture}}
\subfloat{
\begin{tikzpicture}[scale=2]
\filldraw[draw=red, ultra thick, fill=red!60] (.25,.33) -- (.75, .33) -- (.75,1) -- (.25,1) -- cycle;
\draw[red, ultra thick] (.25,.55) -- (.75, .55);
\draw[red, ultra thick] (.25,.66) -- (.75, .66);
\draw[red, ultra thick] (.5,.33) -- (.5, 1);
\draw[red, ultra thick] (.66,.33) -- (.66, 1);
\draw (0,0) -- (1.25,0) -- (1.25,1) -- (0,1) -- cycle;
\draw (0,.55) -- (.75, .55);
\draw (.5, .44) -- (1.25, .44);
\draw (.583, 0) -- (.583, .66);
\draw (.666, .33) -- (.666, 1);
\draw (0,.33) -- (1.25, .33);
\draw (0,.66) -- (1.25, .66);
\draw (.25, 0) -- (.25, 1);
\draw (.5, 0) -- (.5, 1);
\draw (.75, 0) -- (.75, 1);
\draw (1, 0) -- (1, 1);
\end{tikzpicture}}
\subfloat{
\begin{tikzpicture}[scale=2]
\filldraw[draw = red, ultra thick, fill=red!60] (.5,.33) -- (.75,.33) -- (.75, .66) -- (.5, .66) -- cycle;
\draw[red, ultra thick] (.5, .44) -- (.75, .44);
\draw[red, ultra thick] (.5, .55) -- (.75, .55);
\draw[red, ultra thick] (.583, .33) -- (.583, .66);
\draw[red, ultra thick] (.666, .33) -- (.666, .66);
\draw (0,0) -- (1.25,0) -- (1.25,1) -- (0,1) -- cycle;
\draw (0,.55) -- (.75, .55);
\draw (.5, .44) -- (1.25, .44);
\draw (.583, 0) -- (.583, .66);
\draw (.666, .33) -- (.666, 1);
\draw (0,.33) -- (1.25, .33);
\draw (0,.66) -- (1.25, .66);
\draw (.25, 0) -- (.25, 1);
\draw (.5, 0) -- (.5, 1);
\draw (.75, 0) -- (.75, 1);
\draw (1, 0) -- (1, 1);
\end{tikzpicture}
}
\subfloat{
\begin{tikzpicture}[scale=2]
\filldraw[draw=red, ultra thick, fill=red!60] (.5,.33) -- (1, .33) -- (1,1) -- (.5,1) -- cycle;
\draw[red, ultra thick] (.5,.44) -- (1, .44);
\draw[red, ultra thick] (.5,.66) -- (1, .66);
\draw[red, ultra thick] (.666,.33) -- (.666, 1);
\draw[red, ultra thick] (.75,.33) -- (.75, 1);
\draw (0,0) -- (1.25,0) -- (1.25,1) -- (0,1) -- cycle;
\draw (0,.55) -- (.75, .55);
\draw (.5, .44) -- (1.25, .44);
\draw (.583, 0) -- (.583, .66);
\draw (.666, .33) -- (.666, 1);
\draw (0,.33) -- (1.25, .33);
\draw (0,.66) -- (1.25, .66);
\draw (.25, 0) -- (.25, 1);
\draw (.5, 0) -- (.5, 1);
\draw (.75, 0) -- (.75, 1);
\draw (1, 0) -- (1, 1);
\end{tikzpicture}}
\caption{The supports of the six MS B-splines of degree (2,2) in linear a dependence relation on the LR-mesh depicted in Figure \ref{7} (a).}\label{22}
\end{figure}\\
For any other bidegree $(p_1,p_2)\neq (0,0)$, one can build an LR-mesh preserving the same structure of Figure \ref{7} (a). Figure \ref{deg34} shows the cases for $(p_1,p_2)=(3,3),(4,4),(1,1),(1,0), (3,1)$. The insertions are the same as for bidegree $(2,2)$ if $p_k \geq 2$ for some $k \in \{1,2\}$ while if $(p_1,p_2) = (1,0), (0,1), (1,1)$, then it is necessery to use some extensions to get an equivalent arrangement (see the dashed meshlines in the mesh (e) and (g) of Figure \ref{deg34}). Again the dimension of the spline space is 9 while there are 10 MS B-splines in all the cases.\\
\begin{figure}\centering
\subfloat[]{
\begin{tikzpicture}[scale = 3]
\draw (0,0) -- (1.5,0) -- (1.5, 1) -- (0,1) -- cycle;
\draw (0,.25) -- (1.5, .25);
\draw (0,.5) -- (1.5, .5);
\draw (0,.75) -- (1.5,.75);
\draw (.25,0) -- (.25, 1);
\draw (.5,0) -- (.5, 1);
\draw (.75,0) -- (.75, 1);
\draw (1,0) -- (1, 1);
\draw (1.25,0) -- (1.25, 1);
\draw (0,.625) -- (1,.625);
\draw (.5,.375) -- (1.5,.375);
\draw (.625,1) -- (.625, .25);
\draw (.875,0) -- (.875,.75);
\end{tikzpicture}
}~
\subfloat[]{
\begin{minipage}{6cm}\vspace{-3cm}
\begin{tabular}{c}
\begin{tikzpicture}[scale=1.5]
\filldraw[draw = red, ultra thick, fill=red!60] (.25,0) -- (1,0) -- (1,.75) -- (.25,.75) -- cycle;
\draw[red, ultra thick] (.25,.25) -- (1,.25);
\draw[red, ultra thick] (.25,.5) -- (1,.5);
\draw[red, ultra thick] (.25,.625) -- (1,.625);
\draw[red, ultra thick] (.5,0) -- (.5,.75);
\draw[red, ultra thick] (.75,0) -- (.75,.75);
\draw[red, ultra thick] (.875,0) -- (.875,.75);
\draw (0,0) -- (1.5,0) -- (1.5, 1) -- (0,1) -- cycle;
\draw (0,.25) -- (1.5, .25);
\draw (0,.5) -- (1.5, .5);
\draw (0,.75) -- (1.5,.75);
\draw (.25,0) -- (.25, 1);
\draw (.5,0) -- (.5, 1);
\draw (.75,0) -- (.75, 1);
\draw (1,0) -- (1, 1);
\draw (1.25,0) -- (1.25, 1);
\draw (0,.625) -- (1,.625);
\draw (.5,.375) -- (1.5,.375);
\draw (.625,1) -- (.625, .25);
\draw (.875,0) -- (.875,.75);
\end{tikzpicture},
\begin{tikzpicture}[scale=1.5]
\filldraw[draw = red, ultra thick, fill=red!60] (.25,0) -- (1.25,0) -- (1.25,1) -- (.25,1) -- cycle;
\draw[red, ultra thick] (.25,.25) -- (1.25,.25);
\draw[red, ultra thick] (.25,.5) -- (1.25,.5);
\draw[red, ultra thick] (.25,.75) -- (1.25,.75);
\draw[red, ultra thick] (.5,0) -- (.5,1);
\draw[red, ultra thick] (.75,0) -- (.75,1);
\draw[red, ultra thick] (1,0) -- (1,1);
\draw (0,0) -- (1.5,0) -- (1.5, 1) -- (0,1) -- cycle;
\draw (0,.25) -- (1.5, .25);
\draw (0,.5) -- (1.5, .5);
\draw (0,.75) -- (1.5,.75);
\draw (.25,0) -- (.25, 1);
\draw (.5,0) -- (.5, 1);
\draw (.75,0) -- (.75, 1);
\draw (1,0) -- (1, 1);
\draw (1.25,0) -- (1.25, 1);
\draw (0,.625) -- (1,.625);
\draw (.5,.375) -- (1.5,.375);
\draw (.625,1) -- (.625, .25);
\draw (.875,0) -- (.875,.75);
\end{tikzpicture},
\begin{tikzpicture}[scale=1.5]
\filldraw[draw = red, ultra thick, fill=red!60] (.5,0) -- (1.25,0) -- (1.25,.75) -- (.5,.75) -- cycle;
\draw[red, ultra thick] (.5,.25) -- (1.25,.25);
\draw[red, ultra thick] (.5,.375) -- (1.25,.375);
\draw[red, ultra thick] (.5,.5) -- (1.25,.5);
\draw[red, ultra thick] (.75,0) -- (.75,.75);
\draw[red, ultra thick] (.875,0) -- (.875,.75);
\draw[red, ultra thick] (1,0) -- (1,.75);
\draw (0,0) -- (1.5,0) -- (1.5, 1) -- (0,1) -- cycle;
\draw (0,.25) -- (1.5, .25);
\draw (0,.5) -- (1.5, .5);
\draw (0,.75) -- (1.5,.75);
\draw (.25,0) -- (.25, 1);
\draw (.5,0) -- (.5, 1);
\draw (.75,0) -- (.75, 1);
\draw (1,0) -- (1, 1);
\draw (1.25,0) -- (1.25, 1);
\draw (0,.625) -- (1,.625);
\draw (.5,.375) -- (1.5,.375);
\draw (.625,1) -- (.625, .25);
\draw (.875,0) -- (.875,.75);
\end{tikzpicture}\\
\begin{tikzpicture}[scale=1.5]
\filldraw[draw = red, ultra thick, fill=red!60] (.25,.25) -- (1,.25) -- (1,1) -- (.25,1) -- cycle;
\draw[red, ultra thick] (.25,.5) -- (1,.5);
\draw[red, ultra thick] (.25,.625) -- (1,.625);
\draw[red, ultra thick] (.25,.75) -- (1,.75);
\draw[red, ultra thick] (.5,.25) -- (.5,1);
\draw[red, ultra thick] (.625,.25) -- (.625,1);
\draw[red, ultra thick] (.75,.25) -- (.75,1);
\draw (0,0) -- (1.5,0) -- (1.5, 1) -- (0,1) -- cycle;
\draw (0,.25) -- (1.5, .25);
\draw (0,.5) -- (1.5, .5);
\draw (0,.75) -- (1.5,.75);
\draw (.25,0) -- (.25, 1);
\draw (.5,0) -- (.5, 1);
\draw (.75,0) -- (.75, 1);
\draw (1,0) -- (1, 1);
\draw (1.25,0) -- (1.25, 1);
\draw (0,.625) -- (1,.625);
\draw (.5,.375) -- (1.5,.375);
\draw (.625,1) -- (.625, .25);
\draw (.875,0) -- (.875,.75);
\end{tikzpicture},
\begin{tikzpicture}[scale=1.5]
\filldraw[draw = red, ultra thick, fill=red!60] (.5,.25) -- (1,.25) -- (1,.75) -- (.5,.75) -- cycle;
\draw[red, ultra thick] (.5,.375) -- (1,.375);
\draw[red, ultra thick] (.5,.5) -- (1,.5);
\draw[red, ultra thick] (.5,.625) -- (1,.625);
\draw[red, ultra thick] (.625,.25) -- (.625,.75);
\draw[red, ultra thick] (.75,.25) -- (.75,.75);
\draw[red, ultra thick] (.875,.25) -- (.875,.75);
\draw (0,0) -- (1.5,0) -- (1.5, 1) -- (0,1) -- cycle;
\draw (0,.25) -- (1.5, .25);
\draw (0,.5) -- (1.5, .5);
\draw (0,.75) -- (1.5,.75);
\draw (.25,0) -- (.25, 1);
\draw (.5,0) -- (.5, 1);
\draw (.75,0) -- (.75, 1);
\draw (1,0) -- (1, 1);
\draw (1.25,0) -- (1.25, 1);
\draw (0,.625) -- (1,.625);
\draw (.5,.375) -- (1.5,.375);
\draw (.625,1) -- (.625, .25);
\draw (.875,0) -- (.875,.75);
\end{tikzpicture},
\begin{tikzpicture}[scale=1.5]
\filldraw[draw = red, ultra thick, fill=red!60] (.5,.25) -- (1.25,.25) -- (1.25,1) -- (.5,1) -- cycle;
\draw[red, ultra thick] (.5,.375) -- (1.25,.375);
\draw[red, ultra thick] (.5,.5) -- (1.25,.5);
\draw[red, ultra thick] (.5,.75) -- (1.25,.75);
\draw[red, ultra thick] (.625,.25) -- (.625,1);
\draw[red, ultra thick] (.75,.25) -- (.75,1);
\draw[red, ultra thick] (1,.25) -- (1,1);
\draw (0,0) -- (1.5,0) -- (1.5, 1) -- (0,1) -- cycle;
\draw (0,.25) -- (1.5, .25);
\draw (0,.5) -- (1.5, .5);
\draw (0,.75) -- (1.5,.75);
\draw (.25,0) -- (.25, 1);
\draw (.5,0) -- (.5, 1);
\draw (.75,0) -- (.75, 1);
\draw (1,0) -- (1, 1);
\draw (1.25,0) -- (1.25, 1);
\draw (0,.625) -- (1,.625);
\draw (.5,.375) -- (1.5,.375);
\draw (.625,1) -- (.625, .25);
\draw (.875,0) -- (.875,.75);
\end{tikzpicture}
\end{tabular}\vspace{-.15cm}
\end{minipage}\hspace{1.75cm}
}\\
\subfloat[]{
\begin{tikzpicture}[scale=3.1]
\draw (0,0) -- (1.4,0) -- (1.4,1) -- (0,1) -- cycle;
\draw (0,.2) -- (1.4,.2);
\draw (0,.4) -- (1.4,.4);
\draw (0,.6) -- (1.4,.6);
\draw (0,.8) -- (1.4,.8);
\draw (.2,0) -- (.2,1);
\draw (.4,0) -- (.4,1);
\draw (.6,0) -- (.6,1);
\draw (.8,0) -- (.8,1);
\draw (1,0) -- (1,1);
\draw (1.2,0) -- (1.2,1);
\draw (0,.7) -- (1,.7);
\draw (.4,.3) -- (1.4,.3);
\draw (.5,.2) -- (.5,1);
\draw (.9,0) -- (.9,.8);
\end{tikzpicture}
}~
\subfloat[]{
\begin{minipage}{6cm}\vspace{-3.15cm}
\begin{tabular}{c}
\begin{tikzpicture}[scale=1.6]
\filldraw[draw = red, ultra thick, fill=red!60] (.2,0) -- (1,0) -- (1,.8) -- (.2,.8) -- cycle;
\draw[red, ultra thick] (.2,.2) -- (1,.2);
\draw[red, ultra thick] (.2,.4) -- (1,.4);
\draw[red, ultra thick] (.2,.6) -- (1,.6);
\draw[red, ultra thick] (.2,.7) -- (1,.7);
\draw[red, ultra thick] (.4,0) -- (.4,.8);
\draw[red, ultra thick] (.6,0) -- (.6,.8);
\draw[red, ultra thick] (.8,0) -- (.8,.8);
\draw[red, ultra thick] (.9,0) -- (.9,.8);
\draw (0,0) -- (1.4,0) -- (1.4,1) -- (0,1) -- cycle;
\draw (0,.2) -- (1.4,.2);
\draw (0,.4) -- (1.4,.4);
\draw (0,.6) -- (1.4,.6);
\draw (0,.8) -- (1.4,.8);
\draw (.2,0) -- (.2,1);
\draw (.4,0) -- (.4,1);
\draw (.6,0) -- (.6,1);
\draw (.8,0) -- (.8,1);
\draw (1,0) -- (1,1);
\draw (1.2,0) -- (1.2,1);
\draw (0,.7) -- (1,.7);
\draw (.4,.3) -- (1.4,.3);
\draw (.5,.2) -- (.5,1);
\draw (.9,0) -- (.9,.8);
\end{tikzpicture},
\begin{tikzpicture}[scale=1.6]
\filldraw[draw = red, ultra thick, fill=red!60] (.2,0) -- (1.2,0) -- (1.2,1) -- (.2,1) -- cycle;
\draw[red, ultra thick] (.2,.2) -- (1.2,.2);
\draw[red, ultra thick] (.2,.4) -- (1.2,.4);
\draw[red, ultra thick] (.2,.6) -- (1.2,.6);
\draw[red, ultra thick] (.2,.8) -- (1.2,.8);
\draw[red, ultra thick] (.4,0) -- (.4,1);
\draw[red, ultra thick] (.6,0) -- (.6,1);
\draw[red, ultra thick] (.8,0) -- (.8,1);
\draw[red, ultra thick] (1,0) -- (1,1);
\draw (0,0) -- (1.4,0) -- (1.4,1) -- (0,1) -- cycle;
\draw (0,.2) -- (1.4,.2);
\draw (0,.4) -- (1.4,.4);
\draw (0,.6) -- (1.4,.6);
\draw (0,.8) -- (1.4,.8);
\draw (.2,0) -- (.2,1);
\draw (.4,0) -- (.4,1);
\draw (.6,0) -- (.6,1);
\draw (.8,0) -- (.8,1);
\draw (1,0) -- (1,1);
\draw (1.2,0) -- (1.2,1);
\draw (0,.7) -- (1,.7);
\draw (.4,.3) -- (1.4,.3);
\draw (.5,.2) -- (.5,1);
\draw (.9,0) -- (.9,.8);
\end{tikzpicture},
\begin{tikzpicture}[scale=1.6]
\filldraw[draw = red, ultra thick, fill=red!60] (.4,0) -- (1.2,0) -- (1.2,.8) -- (.4,.8) -- cycle;
\draw[red, ultra thick] (.4,.2) -- (1.2,.2);
\draw[red, ultra thick] (.4,.3) -- (1.2,.3);
\draw[red, ultra thick] (.4,.4) -- (1.2,.4);
\draw[red, ultra thick] (.4,.6) -- (1.2,.6);
\draw[red, ultra thick] (.6,0) -- (.6,.8);
\draw[red, ultra thick] (.8,0) -- (.8,.8);
\draw[red, ultra thick] (.9,0) -- (.9,.8);
\draw[red, ultra thick] (1,0) -- (1,.8);
\draw (0,0) -- (1.4,0) -- (1.4,1) -- (0,1) -- cycle;
\draw (0,.2) -- (1.4,.2);
\draw (0,.4) -- (1.4,.4);
\draw (0,.6) -- (1.4,.6);
\draw (0,.8) -- (1.4,.8);
\draw (.2,0) -- (.2,1);
\draw (.4,0) -- (.4,1);
\draw (.6,0) -- (.6,1);
\draw (.8,0) -- (.8,1);
\draw (1,0) -- (1,1);
\draw (1.2,0) -- (1.2,1);
\draw (0,.7) -- (1,.7);
\draw (.4,.3) -- (1.4,.3);
\draw (.5,.2) -- (.5,1);
\draw (.9,0) -- (.9,.8);
\end{tikzpicture}\\
\begin{tikzpicture}[scale=1.6]
\filldraw[draw = red, ultra thick, fill=red!60] (.2,.2) -- (1,.2) -- (1,1) -- (.2,1) -- cycle;
\draw[red, ultra thick] (.2,.8) -- (1,.8);
\draw[red, ultra thick] (.2,.4) -- (1,.4);
\draw[red, ultra thick] (.2,.6) -- (1,.6);
\draw[red, ultra thick] (.2,.7) -- (1,.7);
\draw[red, ultra thick] (.4,.2) -- (.4,1);
\draw[red, ultra thick] (.6,.2) -- (.6,1);
\draw[red, ultra thick] (.8,.2) -- (.8,1);
\draw[red, ultra thick] (.5,.2) -- (.5,1);
\draw (0,0) -- (1.4,0) -- (1.4,1) -- (0,1) -- cycle;
\draw (0,.2) -- (1.4,.2);
\draw (0,.4) -- (1.4,.4);
\draw (0,.6) -- (1.4,.6);
\draw (0,.8) -- (1.4,.8);
\draw (.2,0) -- (.2,1);
\draw (.4,0) -- (.4,1);
\draw (.6,0) -- (.6,1);
\draw (.8,0) -- (.8,1);
\draw (1,0) -- (1,1);
\draw (1.2,0) -- (1.2,1);
\draw (0,.7) -- (1,.7);
\draw (.4,.3) -- (1.4,.3);
\draw (.5,.2) -- (.5,1);
\draw (.9,0) -- (.9,.8);
\end{tikzpicture},
\begin{tikzpicture}[scale=1.6]
\filldraw[draw = red, ultra thick, fill=red!60] (.4,.2) -- (1,.2) -- (1,.8) -- (.4,.8) -- cycle;
\draw[red, ultra thick] (.4,.4) -- (1,.4);
\draw[red, ultra thick] (.4,.3) -- (1,.3);
\draw[red, ultra thick] (.4,.6) -- (1,.6);
\draw[red, ultra thick] (.4,.7) -- (1,.7);
\draw[red, ultra thick] (.5,.2) -- (.5,.8);
\draw[red, ultra thick] (.6,.2) -- (.6,.8);
\draw[red, ultra thick] (.8,.2) -- (.8,.8);
\draw[red, ultra thick] (.9,.2) -- (.9,.8);
\draw (0,0) -- (1.4,0) -- (1.4,1) -- (0,1) -- cycle;
\draw (0,.2) -- (1.4,.2);
\draw (0,.4) -- (1.4,.4);
\draw (0,.6) -- (1.4,.6);
\draw (0,.8) -- (1.4,.8);
\draw (.2,0) -- (.2,1);
\draw (.4,0) -- (.4,1);
\draw (.6,0) -- (.6,1);
\draw (.8,0) -- (.8,1);
\draw (1,0) -- (1,1);
\draw (1.2,0) -- (1.2,1);
\draw (0,.7) -- (1,.7);
\draw (.4,.3) -- (1.4,.3);
\draw (.5,.2) -- (.5,1);
\draw (.9,0) -- (.9,.8);
\end{tikzpicture},
\begin{tikzpicture}[scale=1.6]
\filldraw[draw = red, ultra thick, fill=red!60] (.4,.2) -- (1.2,.2) -- (1.2,1) -- (.4,1) -- cycle;
\draw[red, ultra thick] (.4,.4) -- (1.2,.4);
\draw[red, ultra thick] (.4,.3) -- (1.2,.3);
\draw[red, ultra thick] (.4,.6) -- (1.2,.6);
\draw[red, ultra thick] (.4,.8) -- (1.2,.8);
\draw[red, ultra thick] (.5,.2) -- (.5,1);
\draw[red, ultra thick] (.6,.2) -- (.6,1);
\draw[red, ultra thick] (.8,.2) -- (.8,1);
\draw[red, ultra thick] (1,.2) -- (1,1);
\draw (0,0) -- (1.4,0) -- (1.4,1) -- (0,1) -- cycle;
\draw (0,.2) -- (1.4,.2);
\draw (0,.4) -- (1.4,.4);
\draw (0,.6) -- (1.4,.6);
\draw (0,.8) -- (1.4,.8);
\draw (.2,0) -- (.2,1);
\draw (.4,0) -- (.4,1);
\draw (.6,0) -- (.6,1);
\draw (.8,0) -- (.8,1);
\draw (1,0) -- (1,1);
\draw (1.2,0) -- (1.2,1);
\draw (0,.7) -- (1,.7);
\draw (.4,.3) -- (1.4,.3);
\draw (.5,.2) -- (.5,1);
\draw (.9,0) -- (.9,.8);
\end{tikzpicture}
\end{tabular}\vspace{-.15cm}
\end{minipage}\hspace{1.75cm}
}\\
\subfloat[]{
\begin{tikzpicture}[scale=5.5]
\draw (0,0) -- (.8,0) -- (.8,.6) -- (0,.6) -- cycle;
\draw (0,.3) -- (.8,.3);
\draw (.2,0) -- (.2,.6);
\draw (.4,0) -- (.4,.6);
\draw (.6,0) -- (.6,.6);
\draw (0,.45) -- (.4,.45);
\draw[dashed] (.4,.45) -- (.5, .45);
\draw (.4,.15) -- (.8,.15);
\draw[dashed] (.3,.15) -- (.4,.15);
\draw (.3,0) -- (.3,.45);
\draw (.5,.6) -- (.5,.15);
\end{tikzpicture}}\quad~
\subfloat[]{
\begin{minipage}{6cm}\vspace{-3.15cm}
\begin{tabular}{c}
\begin{tikzpicture}[scale=2.6]
\filldraw[draw = red, ultra thick, fill=red!60] (.2,0) -- (.4,0) -- (.4,.45) -- (.2,.45) -- cycle;
\draw[red, ultra thick] (.2,.3) -- (.4,.3);
\draw[red, ultra thick] (.3,0) -- (.3,.45);
\draw (0,0) -- (.8,0) -- (.8,.6) -- (0,.6) -- cycle;
\draw (0,.3) -- (.8,.3);
\draw (.2,0) -- (.2,.6);
\draw (.4,0) -- (.4,.6);
\draw (.6,0) -- (.6,.6);
\draw (0,.45) -- (.5,.45);
\draw (.3,.15) -- (.8,.15);
\draw (.3,0) -- (.3,.45);
\draw (.5,.6) -- (.5,.15);
\end{tikzpicture},
\begin{tikzpicture}[scale=2.6]
\filldraw[draw = red, ultra thick, fill=red!60] (.2,0) -- (.6,0) -- (.6,.6) -- (.2,.6) -- cycle;
\draw[red, ultra thick] (.2,.3) -- (.6,.3);
\draw[red, ultra thick] (.4,0) -- (.4,.6);
\draw (0,0) -- (.8,0) -- (.8,.6) -- (0,.6) -- cycle;
\draw (0,.3) -- (.8,.3);
\draw (.2,0) -- (.2,.6);
\draw (.4,0) -- (.4,.6);
\draw (.6,0) -- (.6,.6);
\draw (0,.45) -- (.5,.45);
\draw (.3,.15) -- (.8,.15);
\draw (.3,0) -- (.3,.45);
\draw (.5,.6) -- (.5,.15);
\end{tikzpicture},
\begin{tikzpicture}[scale=2.6]
\filldraw[draw = red, ultra thick, fill=red!60] (.3,0) -- (.6,0) -- (.6,.3) -- (.3,.3) -- cycle;
\draw[red, ultra thick] (.3,.15) -- (.6,.15);
\draw[red, ultra thick] (.4,0) -- (.4,.3);
\draw (0,0) -- (.8,0) -- (.8,.6) -- (0,.6) -- cycle;
\draw (0,.3) -- (.8,.3);
\draw (.2,0) -- (.2,.6);
\draw (.4,0) -- (.4,.6);
\draw (.6,0) -- (.6,.6);
\draw (0,.45) -- (.5,.45);
\draw (.3,.15) -- (.8,.15);
\draw (.3,0) -- (.3,.45);
\draw (.5,.6) -- (.5,.15);
\end{tikzpicture}\\
\begin{tikzpicture}[scale=2.6]
\filldraw[draw = red, ultra thick, fill=red!60] (.2,.3) -- (.5,.3) -- (.5,.6) -- (.2,.6) -- cycle;
\draw[red, ultra thick] (.2,.45) -- (.5,.45);
\draw[red, ultra thick] (.4,.3) -- (.4,.6);
\draw (0,0) -- (.8,0) -- (.8,.6) -- (0,.6) -- cycle;
\draw (0,.3) -- (.8,.3);
\draw (.2,0) -- (.2,.6);
\draw (.4,0) -- (.4,.6);
\draw (.6,0) -- (.6,.6);
\draw (0,.45) -- (.5,.45);
\draw (.3,.15) -- (.8,.15);
\draw (.3,0) -- (.3,.45);
\draw (.5,.6) -- (.5,.15);
\end{tikzpicture},
\begin{tikzpicture}[scale=2.6]
\filldraw[draw = red, ultra thick, fill=red!60] (.3,.15) -- (.5,.15) -- (.5,.45) -- (.3,.45) -- cycle;
\draw[red, ultra thick] (.3,.3) -- (.5,.3);
\draw[red, ultra thick] (.4,.15) -- (.4,.45);
\draw (0,0) -- (.8,0) -- (.8,.6) -- (0,.6) -- cycle;
\draw (0,.3) -- (.8,.3);
\draw (.2,0) -- (.2,.6);
\draw (.4,0) -- (.4,.6);
\draw (.6,0) -- (.6,.6);
\draw (0,.45) -- (.5,.45);
\draw (.3,.15) -- (.8,.15);
\draw (.3,0) -- (.3,.45);
\draw (.5,.6) -- (.5,.15);
\end{tikzpicture},
\begin{tikzpicture}[scale=2.6]
\filldraw[draw = red, ultra thick, fill=red!60] (.4,.15) -- (.6,.15) -- (.6,.6) -- (.4,.6) -- cycle;
\draw[red, ultra thick] (.4,.3) -- (.6,.3);
\draw[red, ultra thick] (.5,.6) -- (.5,.15);
\draw (0,0) -- (.8,0) -- (.8,.6) -- (0,.6) -- cycle;
\draw (0,.3) -- (.8,.3);
\draw (.2,0) -- (.2,.6);
\draw (.4,0) -- (.4,.6);
\draw (.6,0) -- (.6,.6);
\draw (0,.45) -- (.5,.45);
\draw (.3,.15) -- (.8,.15);
\draw (.3,0) -- (.3,.45);
\draw (.5,.6) -- (.5,.15);
\end{tikzpicture}
\end{tabular}\vspace{-.15cm}
\end{minipage}\hspace{1.25cm}
}\\
\subfloat[]{
\begin{tikzpicture}[scale=5.5]
\draw (0,0) -- (.8,0) -- (.8,.6) -- (0,.6) -- cycle;
\draw (.2,0) -- (.2,.6);
\draw (.4,0) -- (.4,.6);
\draw (.6,0) -- (.6,.6);
\draw (0,.45) -- (.4,.45);
\draw[dashed] (.4,.45) -- (.5, .45);
\draw (.4,.15) -- (.8,.15);
\draw[dashed] (.3,.15) -- (.4,.15);
\draw (.3,0) -- (.3,.45);
\draw (.5,.6) -- (.5,.15);
\end{tikzpicture}}\quad~
\subfloat[]{
\begin{minipage}{6cm}\vspace{-3.15cm}
\begin{tabular}{c}
\begin{tikzpicture}[scale=2.6]
\filldraw[draw = red, ultra thick, fill=red!60] (.2,0) -- (.4,0) -- (.4,.45) -- (.2,.45) -- cycle;
\draw[red, ultra thick] (.3,0) -- (.3,.45);
\draw (0,0) -- (.8,0) -- (.8,.6) -- (0,.6) -- cycle;
\draw (.2,0) -- (.2,.6);
\draw (.4,0) -- (.4,.6);
\draw (.6,0) -- (.6,.6);
\draw (0,.45) -- (.5,.45);
\draw (.3,.15) -- (.8,.15);
\draw (.3,0) -- (.3,.45);
\draw (.5,.6) -- (.5,.15);
\end{tikzpicture},
\begin{tikzpicture}[scale=2.6]
\filldraw[draw = red, ultra thick, fill=red!60] (.2,0) -- (.6,0) -- (.6,.6) -- (.2,.6) -- cycle;
\draw[red, ultra thick] (.4,0) -- (.4,.6);
\draw (0,0) -- (.8,0) -- (.8,.6) -- (0,.6) -- cycle;
\draw (.2,0) -- (.2,.6);
\draw (.4,0) -- (.4,.6);
\draw (.6,0) -- (.6,.6);
\draw (0,.45) -- (.5,.45);
\draw (.3,.15) -- (.8,.15);
\draw (.3,0) -- (.3,.45);
\draw (.5,.6) -- (.5,.15);
\end{tikzpicture},
\begin{tikzpicture}[scale=2.6]
\filldraw[draw = red, ultra thick, fill=red!60] (.3,0) -- (.6,0) -- (.6,.15) -- (.3,.15) -- cycle;
\draw[red, ultra thick] (.4,0) -- (.4,.15);
\draw (0,0) -- (.8,0) -- (.8,.6) -- (0,.6) -- cycle;
\draw (.2,0) -- (.2,.6);
\draw (.4,0) -- (.4,.6);
\draw (.6,0) -- (.6,.6);
\draw (0,.45) -- (.5,.45);
\draw (.3,.15) -- (.8,.15);
\draw (.3,0) -- (.3,.45);
\draw (.5,.6) -- (.5,.15);
\end{tikzpicture}\\
\begin{tikzpicture}[scale=2.6]
\filldraw[draw = red, ultra thick, fill=red!60] (.2,.45) -- (.5,.45) -- (.5,.6) -- (.2,.6) -- cycle;
\draw[red, ultra thick] (.4,.45) -- (.4,.6);
\draw (0,0) -- (.8,0) -- (.8,.6) -- (0,.6) -- cycle;
\draw (0,.3) -- (.8,.3);
\draw (.2,0) -- (.2,.6);
\draw (.4,0) -- (.4,.6);
\draw (.6,0) -- (.6,.6);
\draw (0,.45) -- (.5,.45);
\draw (.3,.15) -- (.8,.15);
\draw (.3,0) -- (.3,.45);
\draw (.5,.6) -- (.5,.15);
\end{tikzpicture},
\begin{tikzpicture}[scale=2.6]
\filldraw[draw = red, ultra thick, fill=red!60] (.3,.15) -- (.5,.15) -- (.5,.45) -- (.3,.45) -- cycle;
\draw[red, ultra thick] (.4,.15) -- (.4,.45);
\draw (0,0) -- (.8,0) -- (.8,.6) -- (0,.6) -- cycle;
\draw (.2,0) -- (.2,.6);
\draw (.4,0) -- (.4,.6);
\draw (.6,0) -- (.6,.6);
\draw (0,.45) -- (.5,.45);
\draw (.3,.15) -- (.8,.15);
\draw (.3,0) -- (.3,.45);
\draw (.5,.6) -- (.5,.15);
\end{tikzpicture},
\begin{tikzpicture}[scale=2.6]
\filldraw[draw = red, ultra thick, fill=red!60] (.4,.15) -- (.6,.15) -- (.6,.6) -- (.4,.6) -- cycle;
\draw[red, ultra thick] (.5,.6) -- (.5,.15);
\draw (0,0) -- (.8,0) -- (.8,.6) -- (0,.6) -- cycle;
\draw (.2,0) -- (.2,.6);
\draw (.4,0) -- (.4,.6);
\draw (.6,0) -- (.6,.6);
\draw (0,.45) -- (.5,.45);
\draw (.3,.15) -- (.8,.15);
\draw (.3,0) -- (.3,.45);
\draw (.5,.6) -- (.5,.15);
\end{tikzpicture}
\end{tabular}\vspace{-.15cm}
\end{minipage}\hspace{1.25cm}
}\\
\subfloat[]{
\begin{tikzpicture}[scale = 3]
\draw (0,0) -- (1.5,0) -- (1.5, 1) -- (0,1) -- cycle;
\draw (0,.5) -- (1.5, .5);
\draw (.25,0) -- (.25, 1);
\draw (.5,0) -- (.5, 1);
\draw (.75,0) -- (.75, 1);
\draw (1,0) -- (1, 1);
\draw (1.25,0) -- (1.25, 1);
\draw (0,.75) -- (1,.75);
\draw (.5,.25) -- (1.5,.25);
\draw (.625,1) -- (.625, .25);
\draw (.875,0) -- (.875,.75);
\end{tikzpicture}
}~
\subfloat[]{
\begin{minipage}{6cm}\vspace{-3cm}
\begin{tabular}{c}
\begin{tikzpicture}[scale=1.5]
\filldraw[draw = red, ultra thick, fill=red!60] (.25,0) -- (1,0) -- (1,.75) -- (.25,.75) -- cycle;
\draw[red, ultra thick] (.25,.5) -- (1,.5);
\draw[red, ultra thick] (.5,0) -- (.5,.75);
\draw[red, ultra thick] (.75,0) -- (.75,.75);
\draw[red, ultra thick] (.875,0) -- (.875,.75);
\draw (0,0) -- (1.5,0) -- (1.5, 1) -- (0,1) -- cycle;
\draw (0,.5) -- (1.5, .5);
\draw (.25,0) -- (.25, 1);
\draw (.5,0) -- (.5, 1);
\draw (.75,0) -- (.75, 1);
\draw (1,0) -- (1, 1);
\draw (1.25,0) -- (1.25, 1);
\draw (0,.75) -- (1,.75);
\draw (.5,.25) -- (1.5,.25);
\draw (.625,1) -- (.625, .25);
\draw (.875,0) -- (.875,.75);
\end{tikzpicture},
\begin{tikzpicture}[scale=1.5]
\filldraw[draw = red, ultra thick, fill=red!60] (.25,0) -- (1.25,0) -- (1.25,1) -- (.25,1) -- cycle;
\draw[red, ultra thick] (.25,.5) -- (1.25,.5);
\draw[red, ultra thick] (.5,0) -- (.5,1);
\draw[red, ultra thick] (.75,0) -- (.75,1);
\draw[red, ultra thick] (1,0) -- (1,1);
\draw (0,0) -- (1.5,0) -- (1.5, 1) -- (0,1) -- cycle;
\draw (0,.5) -- (1.5, .5);
\draw (.25,0) -- (.25, 1);
\draw (.5,0) -- (.5, 1);
\draw (.75,0) -- (.75, 1);
\draw (1,0) -- (1, 1);
\draw (1.25,0) -- (1.25, 1);
\draw (0,.75) -- (1,.75);
\draw (.5,.25) -- (1.5,.25);
\draw (.625,1) -- (.625, .25);
\draw (.875,0) -- (.875,.75);
\end{tikzpicture},
\begin{tikzpicture}[scale=1.5]
\filldraw[draw = red, ultra thick, fill=red!60] (.5,0) -- (1.25,0) -- (1.25,.5) -- (.5,.5) -- cycle;
\draw[red, ultra thick] (.5,.25) -- (1.25,.25);
\draw[red, ultra thick] (.75,0) -- (.75,.5);
\draw[red, ultra thick] (.875,0) -- (.875,.5);
\draw[red, ultra thick] (1,0) -- (1,.5);
\draw (0,0) -- (1.5,0) -- (1.5, 1) -- (0,1) -- cycle;
\draw (0,.5) -- (1.5, .5);
\draw (.25,0) -- (.25, 1);
\draw (.5,0) -- (.5, 1);
\draw (.75,0) -- (.75, 1);
\draw (1,0) -- (1, 1);
\draw (1.25,0) -- (1.25, 1);
\draw (0,.75) -- (1,.75);
\draw (.5,.25) -- (1.5,.25);
\draw (.625,1) -- (.625, .25);
\draw (.875,0) -- (.875,.75);
\end{tikzpicture}\\
\begin{tikzpicture}[scale=1.5]
\filldraw[draw = red, ultra thick, fill=red!60] (.25,.5) -- (1,.5) -- (1,1) -- (.25,1) -- cycle;
\draw[red, ultra thick] (.25,.75) -- (1,.75);
\draw[red, ultra thick] (.5,.5) -- (.5,1);
\draw[red, ultra thick] (.625,.5) -- (.625,1);
\draw[red, ultra thick] (.75,.5) -- (.75,1);
\draw (0,0) -- (1.5,0) -- (1.5, 1) -- (0,1) -- cycle;
\draw (0,.5) -- (1.5, .5);
\draw (.25,0) -- (.25, 1);
\draw (.5,0) -- (.5, 1);
\draw (.75,0) -- (.75, 1);
\draw (1,0) -- (1, 1);
\draw (1.25,0) -- (1.25, 1);
\draw (0,.75) -- (1,.75);
\draw (.5,.25) -- (1.5,.25);
\draw (.625,1) -- (.625, .25);
\draw (.875,0) -- (.875,.75);
\end{tikzpicture},
\begin{tikzpicture}[scale=1.5]
\filldraw[draw = red, ultra thick, fill=red!60] (.5,.25) -- (1,.25) -- (1,.75) -- (.5,.75) -- cycle;
\draw[red, ultra thick] (.5,.5) -- (1,.5);
\draw[red, ultra thick] (.625,.25) -- (.625,.75);
\draw[red, ultra thick] (.75,.25) -- (.75,.75);
\draw[red, ultra thick] (.875,.25) -- (.875,.75);
\draw (0,0) -- (1.5,0) -- (1.5, 1) -- (0,1) -- cycle;
\draw (0,.5) -- (1.5, .5);
\draw (.25,0) -- (.25, 1);
\draw (.5,0) -- (.5, 1);
\draw (.75,0) -- (.75, 1);
\draw (1,0) -- (1, 1);
\draw (1.25,0) -- (1.25, 1);
\draw (0,.75) -- (1,.75);
\draw (.5,.25) -- (1.5,.25);
\draw (.625,1) -- (.625, .25);
\draw (.875,0) -- (.875,.75);
\end{tikzpicture},
\begin{tikzpicture}[scale=1.5]
\filldraw[draw = red, ultra thick, fill=red!60] (.5,.25) -- (1.25,.25) -- (1.25,1) -- (.5,1) -- cycle;
\draw[red, ultra thick] (.5,.5) -- (1.25,.5);
\draw[red, ultra thick] (.625,.25) -- (.625,1);
\draw[red, ultra thick] (.75,.25) -- (.75,1);
\draw[red, ultra thick] (1,.25) -- (1,1);
\draw (0,0) -- (1.5,0) -- (1.5, 1) -- (0,1) -- cycle;
\draw (0,.5) -- (1.5, .5);
\draw (.25,0) -- (.25, 1);
\draw (.5,0) -- (.5, 1);
\draw (.75,0) -- (.75, 1);
\draw (1,0) -- (1, 1);
\draw (1.25,0) -- (1.25, 1);
\draw (0,.75) -- (1,.75);
\draw (.5,.25) -- (1.5,.25);
\draw (.625,1) -- (.625, .25);
\draw (.875,0) -- (.875,.75);
\end{tikzpicture}
\end{tabular}\vspace{-.15cm}
\end{minipage}\hspace{1.75cm}
}
\caption{In (a) an LR-mesh providing MS B-splines of bidegree (3,3) in an equivalent arrengement of the MS B-splines of bidegree (2,2) on the mesh in Figure \ref{7} (a). In (b) are shown the supports of the six B-splines in the linear dependence relation. In (c)-(d) the same for bidegree (4,4). In (e)-(f) we have the same  for bidegree (1,1). Note that we have used two extensions (dashed meshlines) to obtain an equivalent arrangement as for the other bidegrees. Finally in (g)-(h) and (i)-(j) we show the equivalent configuration for bidegrees (1,0) and (3,1). Again the dashed meshlines are extensions.}\label{deg34}
\end{figure}
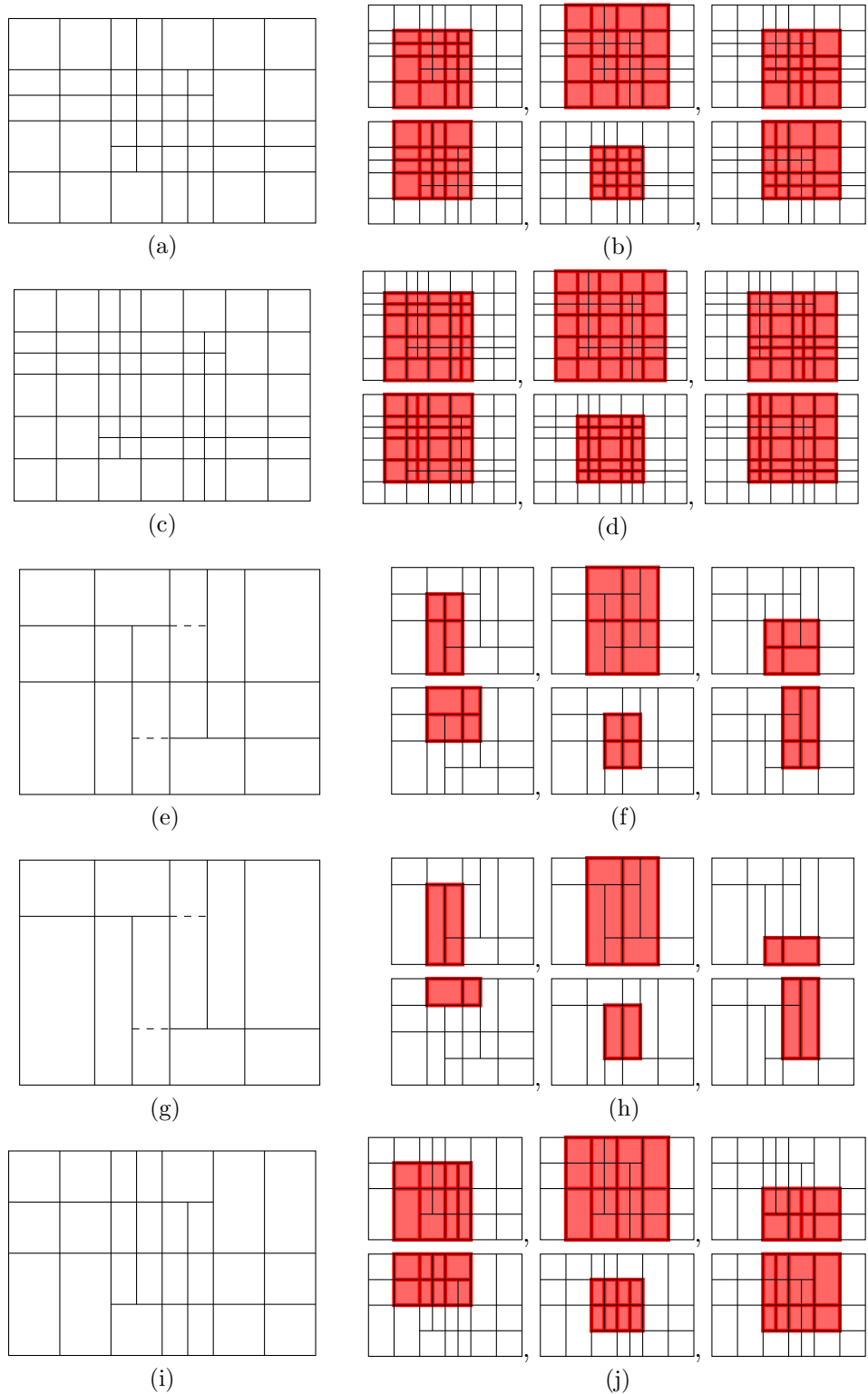\\
\end{ex}

\section{Minimal number of LR B-splines for a linear dependence relation}\label{LDLR}
In this section we show that at least eight B-splines must be involved in a linear dependence relation in $\mathcal{B}^{\mathcal{LR}}(\cM)$. Then we provide examples for any bidegree $\pmb{p} = (p_1,p_2)$ with at least one $p_k \geq 2$ for some $k\in \{1,2\}$ where the LR B-splines in linear dependence are exactly eight. In such examples the meshes will be refinements of the meshes presented in Example \ref{exMS}. As we pointed out in Proposition \ref{noLR}, the sixth MS B-spline $B^6$ in Example \ref{exMS}  is not an LR B-spline on $\cM$. In the new examples here we show how to refine $\cM$ in order to split $B^6$ into two B-splines that can be now obtained through the knot insertion algorithm from LR B-splines on coarser meshes. This will move the number of MS B-splines involved in the linear dependence from 6 to 8 but all of them will now be LR B-splines.
\begin{lem}\label{p3}
Given an LR-mesh $\cM$ and a bidegree $(p_1,p_2)$, assume the elements in $\cB^{\mathcal{LR}}(\cM)$ are linearly independent. If the insertion of a $k$-split $\gamma$ causes a linear dependence relation in  $\cB^{\mathcal{LR}}(\cM+\gamma)$, then the expanded $k$-split $\bar{\gamma}$ is at least $p_{3-k}+3$ long and the growth of cardinality is $|\cB^{\mathcal{LR}}(\cM+\gamma)| - |\cB^{\mathcal{LR}}(\cM)| > 2$.
\begin{proof}
Theorem 5.2 of \cite{dokken2013polynomial} ensures that if $\bar{\gamma}$ is $p_{3-k}+2$ long then the elements in $\cB^{\mathcal{LR}}(\cM+\gamma)$ are linearly independent. Assume that $\bar{\gamma}$ is $p_{3-k}+3$ long. From the end of Section \ref{HIHP}, the refinement goes hand-in-hand only if  $|\cB^{\mathcal{LR}}(\cM+\gamma)| - |\cB^{\mathcal{LR}}(\cM)| \geq  2$ and there is a linear dependence relation if it is a strict inequality. If the refinement does not go hand-in-hand then it must be $|\cB^{\mathcal{LR}}(\cM+\gamma)| - |\cB^{\mathcal{LR}}(\cM)| \leq 1$ and the new B-spline (if exists) is linearly independent of the B-splines in $\cB^{\mathcal{LR}}(\cM)$ as it has a split, given by the insertion of $\gamma$, not in $\cM$.
\end{proof}
\end{lem}
\begin{prop}\label{eightLR}
We need at least 8 LR B-splines for a linear dependence relation in $\cB^{\mathcal{LR}}(\cM)$.
\begin{proof}
By Proposition \ref{L1}, we must have four nested B-splines at the four corners of $\mathcal{R}$. In order to keep the number of B-splines needed as low as possible we assume as in the proof of Corollary \ref{fiveMS} that we only need, for the nestedness condition, five B-splines: $B^2, B^3, B^4, B^5$ contained in $B^1$.
$B^2,B^3,B^4,B^5$ are created through the knot insertion algorithm splitting LR B-splines on coarser meshes during the construction of $\cM$. It means that the size in one direction of each of them must be the same as of another LR B-spline in the linear dependence relation. Moreover, every point $(x,y) \in \mathcal{R}$ must be contained in at least another support than $\mbox{supp}\,B^1$. Therefore either \begin{enumerate}
\item The supports of $B^2, B^3, B^4, B^5$ intersect each other and have the same size in one direction, as the MS B-splines depicted in Figure \ref{14} (c), or
\item the supports of at least two nested B-splines are separated.
\end{enumerate}
In case 1, as we have seen in the proof of Proposition \ref{noLR}, if there are no more relevant meshlines apart from those in the splits of $B^2, B^3, B^4, B^5$, the other MS B-splines that can be generated in $\mathcal{R}$ using relevant meshlines are not LR B-splines.

Therefore, in order to make a linear dependence relation in $\mathcal{R}$, there must exist at least another relevant split that has provided, from Lemma \ref{p3}, a growth in the LR B-spline set of at least three, bringing the number of LR B-splines involved to at least eight.  

In case 2, suppose that the supports of two B-splines among $B^2, B^3, B^4, B^5$ do not intersect each other. 
Then there are T-vertices not shared by two B-splines nested at the corners. There must exist other LR B-splines sharing these T-vertices and bringing linear dependence. Hence, there must exist at least another split, aside from those needed for the contruction of the nested B-splines, that has provided, from Lemma \ref{p3}, a growth of three in the LR B-spline set, moving the total number to at least eight.
\end{proof}
\end{prop}
In the following example we show meshes where there are exactly eight LR B-splines in a linear dependence relation for any bidegree $\pmb{p}= (p_1,p_2)$ with $p_k \geq2$ for some $k \in \{1,2\}$. Such meshes are refinements of the meshes presented in Example \ref{exMS} at the end of the previous section.
\begin{ex}\label{exLR}
Let us consider first bidegree $\pmb{p}=(2,2)$. Consider the mesh $\cM$ of Figure \ref{7} (a). As we have shown in Example \ref{exMS}, $\dim\,\SSS_{(2,2)}(\cM) = 9$ and the construction of $\cM$ went LR-wise hand-in-hand. Let us now insert a new split $\gamma$ of length $p_2+3 = 5$ in $\cM$ to get the mesh $\cM + \gamma$ as shown in Figure \ref{8LR} (a). Then, $\dim\, \SSS_{(2,2)} (\cM + \gamma) = \dim\, \SSS_{(2,2)}(\cM) + 2 = 11$ and $\cM + \gamma$ went LR-wise hand-in-hand with $\cM$. Furthermore, the LR B-spline set grows by three, $|B^{\mathcal{LR}}(\cM+\gamma)| = |\cB^{\mathcal{LR}}(\cM)| + 3 = 12$ as shown in Figure \ref{8LR} (b). Therefore, there is a linear dependence relation. The only eight LR B-splines satisfying Proposition \ref{L1} and Corollary \ref{Tver} are depicted in Figure \ref{8LR} (c).
\begin{figure}[h!]\centering
\subfloat[]{
\begin{tikzpicture}[scale=3]
\draw (0,0) -- (1.25,0) -- (1.25,1) -- (0,1) -- cycle;
\draw (0,.55) -- (.75, .55);
\draw (.5, .44) -- (1.25, .44);
\draw (.562, .33) -- (.562, 1);
\draw (.686, .33) -- (.686, 1);
\draw (0,.33) -- (1.25, .33);
\draw (0,.66) -- (1.25, .66);
\draw (.25, 0) -- (.25, 1);
\draw (.5, 0) -- (.5, 1);
\draw (.75, 0) -- (.75, 1);
\draw (1, 0) -- (1, 1);
\draw (.624,0) -- (.624, .66);
\end{tikzpicture}}
\subfloat[]{
\begin{minipage}{12cm}\vspace{-3cm}
\begin{tabular}{c}
\begin{tikzpicture}[scale=1.5]
\filldraw[draw=red, ultra thick, fill=red!60] (0,0) -- (.624, 0) -- (.624,.66) -- (0,.66) -- cycle;
\draw[red, ultra thick] (0,.33) -- (.624, .33);
\draw[red, ultra thick] (0,.55) -- (.624, .55);
\draw[red, ultra thick] (.25,0) -- (.25, .66);
\draw[red, ultra thick] (.5,0) -- (.5, .66);
\draw (0,0) -- (1.25,0) -- (1.25,1) -- (0,1) -- cycle;
\draw (0,.55) -- (.75, .55);
\draw (.5, .44) -- (1.25, .44);
\draw (.562, .33) -- (.562, 1);
\draw (.686, .33) -- (.686, 1);
\draw (0,.33) -- (1.25, .33);
\draw (0,.66) -- (1.25, .66);
\draw (.25, 0) -- (.25, 1);
\draw (.5, 0) -- (.5, 1);
\draw (.75, 0) -- (.75, 1);
\draw (1, 0) -- (1, 1);
\draw (.624,0) -- (.624, .66);
\end{tikzpicture},
\begin{tikzpicture}[scale=1.5]
\filldraw[draw=red, ultra thick, fill=red!60] (.25,0) -- (1, 0) -- (1,1) -- (.25,1) -- cycle;
\draw[red, ultra thick] (.25,.33) -- (1, .33);
\draw[red, ultra thick] (.25,.66) -- (1, .66);
\draw[red, ultra thick] (.5,0) -- (.5, 1);
\draw[red, ultra thick] (.75,0) -- (.75, 1);
\draw (0,0) -- (1.25,0) -- (1.25,1) -- (0,1) -- cycle;
\draw (0,.55) -- (.75, .55);
\draw (.5, .44) -- (1.25, .44);
\draw (.562, .33) -- (.562, 1);
\draw (.686, .33) -- (.686, 1);
\draw (0,.33) -- (1.25, .33);
\draw (0,.66) -- (1.25, .66);
\draw (.25, 0) -- (.25, 1);
\draw (.5, 0) -- (.5, 1);
\draw (.75, 0) -- (.75, 1);
\draw (1, 0) -- (1, 1);
\draw (.624,0) -- (.624, .66);
\end{tikzpicture},
\begin{tikzpicture}[scale=1.5]
\filldraw[draw=red, ultra thick, fill=red!60] (.25,0) -- (.75, 0) -- (.75,.66) -- (.25,.66) -- cycle;
\draw[red, ultra thick] (.25,.33) -- (.75, .33);
\draw[red, ultra thick] (.25,.55) -- (.75, .55);
\draw[red, ultra thick] (.5,0) -- (.5, .66);
\draw[red, ultra thick] (.624,0) -- (.624, .66);
\draw (0,0) -- (1.25,0) -- (1.25,1) -- (0,1) -- cycle;
\draw (0,.55) -- (.75, .55);
\draw (.5, .44) -- (1.25, .44);
\draw (.562, .33) -- (.562, 1);
\draw (.686, .33) -- (.686, 1);
\draw (0,.33) -- (1.25, .33);
\draw (0,.66) -- (1.25, .66);
\draw (.25, 0) -- (.25, 1);
\draw (.5, 0) -- (.5, 1);
\draw (.75, 0) -- (.75, 1);
\draw (1, 0) -- (1, 1);
\draw (.624,0) -- (.624, .66);
\end{tikzpicture},
\begin{tikzpicture}[scale=1.5]
\filldraw[draw=red, ultra thick, fill=red!60] (.5,0) -- (1, 0) -- (1,.66) -- (.5,.66) -- cycle;
\draw[red, ultra thick] (.5,.33) -- (1, .33);
\draw[red, ultra thick] (.5,.44) -- (1, .44);
\draw[red, ultra thick] (.625,0) -- (.625, .66);
\draw[red, ultra thick] (.75,0) -- (.75, .66);
\draw (0,0) -- (1.25,0) -- (1.25,1) -- (0,1) -- cycle;
\draw (0,.55) -- (.75, .55);
\draw (.5, .44) -- (1.25, .44);
\draw (.562, .33) -- (.562, 1);
\draw (.686, .33) -- (.686, 1);
\draw (0,.33) -- (1.25, .33);
\draw (0,.66) -- (1.25, .66);
\draw (.25, 0) -- (.25, 1);
\draw (.5, 0) -- (.5, 1);
\draw (.75, 0) -- (.75, 1);
\draw (1, 0) -- (1, 1);
\draw (.624,0) -- (.624, .66);
\end{tikzpicture},
\begin{tikzpicture}[scale=1.5]
\filldraw[draw=red, ultra thick, fill=red!60] (.625,0) -- (1.25, 0) -- (1.25,.66) -- (.625,.66) -- cycle;
\draw[red, ultra thick] (.625,.33) -- (1.25, .33);
\draw[red, ultra thick] (.625,.44) -- (1.25, .44);
\draw[red, ultra thick] (.75,0) -- (.75, .66);
\draw[red, ultra thick] (1,0) -- (1, .66);
\draw (0,0) -- (1.25,0) -- (1.25,1) -- (0,1) -- cycle;
\draw (0,.55) -- (.75, .55);
\draw (.5, .44) -- (1.25, .44);
\draw (.562, .33) -- (.562, 1);
\draw (.686, .33) -- (.686, 1);
\draw (0,.33) -- (1.25, .33);
\draw (0,.66) -- (1.25, .66);
\draw (.25, 0) -- (.25, 1);
\draw (.5, 0) -- (.5, 1);
\draw (.75, 0) -- (.75, 1);
\draw (1, 0) -- (1, 1);
\draw (.624,0) -- (.624, .66);
\end{tikzpicture},
\begin{tikzpicture}[scale=1.5]
\filldraw[draw=red, ultra thick, fill=red!60] (0,.33) -- (.562, .33) -- (.562,1) -- (0,1) -- cycle;
\draw[red, ultra thick] (0,.55) -- (.562, .55);
\draw[red, ultra thick] (0,.66) -- (.562, .66);
\draw[red, ultra thick] (.25,.33) -- (.25, 1);
\draw[red, ultra thick] (.5,.33) -- (.5, 1);
\draw (0,0) -- (1.25,0) -- (1.25,1) -- (0,1) -- cycle;
\draw (0,.55) -- (.75, .55);
\draw (.5, .44) -- (1.25, .44);
\draw (.562, .33) -- (.562, 1);
\draw (.686, .33) -- (.686, 1);
\draw (0,.33) -- (1.25, .33);
\draw (0,.66) -- (1.25, .66);
\draw (.25, 0) -- (.25, 1);
\draw (.5, 0) -- (.5, 1);
\draw (.75, 0) -- (.75, 1);
\draw (1, 0) -- (1, 1);
\draw (.624,0) -- (.624, .66);
\end{tikzpicture}\\
\begin{tikzpicture}[scale=1.5]
\filldraw[draw=red, ultra thick, fill=red!60] (.25,.33) -- (.686, .33) -- (.686,1) -- (.25,1) -- cycle;
\draw[red, ultra thick] (.25,.55) -- (.686, .55);
\draw[red, ultra thick] (.25,.66) -- (.686, .66);
\draw[red, ultra thick] (.5,.33) -- (.5, 1);
\draw[red, ultra thick] (.562,.33) -- (.562, 1);
\draw (0,0) -- (1.25,0) -- (1.25,1) -- (0,1) -- cycle;
\draw (0,.55) -- (.75, .55);
\draw (.5, .44) -- (1.25, .44);
\draw (.562, .33) -- (.562, 1);
\draw (.686, .33) -- (.686, 1);
\draw (0,.33) -- (1.25, .33);
\draw (0,.66) -- (1.25, .66);
\draw (.25, 0) -- (.25, 1);
\draw (.5, 0) -- (.5, 1);
\draw (.75, 0) -- (.75, 1);
\draw (1, 0) -- (1, 1);
\draw (.624,0) -- (.624, .66);
\end{tikzpicture},
\begin{tikzpicture}[scale=1.5]
\filldraw[draw=red, ultra thick, fill=red!60] (.5,.33) -- (.686, .33) -- (.686,.66) -- (.5,.66) -- cycle;
\draw[red, ultra thick] (.5,.55) -- (.686, .55);
\draw[red, ultra thick] (.5,.44) -- (.686, .44);
\draw[red, ultra thick] (.562,.33) -- (.562, .66);
\draw[red, ultra thick] (.624,.33) -- (.624, .66);
\draw (0,0) -- (1.25,0) -- (1.25,1) -- (0,1) -- cycle;
\draw (0,.55) -- (.75, .55);
\draw (.5, .44) -- (1.25, .44);
\draw (.562, .33) -- (.562, 1);
\draw (.686, .33) -- (.686, 1);
\draw (0,.33) -- (1.25, .33);
\draw (0,.66) -- (1.25, .66);
\draw (.25, 0) -- (.25, 1);
\draw (.5, 0) -- (.5, 1);
\draw (.75, 0) -- (.75, 1);
\draw (1, 0) -- (1, 1);
\draw (.624,0) -- (.624, .66);
\end{tikzpicture},
\begin{tikzpicture}[scale=1.5]
\filldraw[draw=red, ultra thick, fill=red!60] (.562,.33) -- (.75, .33) -- (.75,.66) -- (.562,.66) -- cycle;
\draw[red, ultra thick] (.562,.55) -- (.75, .55);
\draw[red, ultra thick] (.562,.44) -- (.75, .44);
\draw[red, ultra thick] (.686,.33) -- (.686, .66);
\draw[red, ultra thick] (.624,.33) -- (.624, .66);
\draw (0,0) -- (1.25,0) -- (1.25,1) -- (0,1) -- cycle;
\draw (0,.55) -- (.75, .55);
\draw (.5, .44) -- (1.25, .44);
\draw (.562, .33) -- (.562, 1);
\draw (.686, .33) -- (.686, 1);
\draw (0,.33) -- (1.25, .33);
\draw (0,.66) -- (1.25, .66);
\draw (.25, 0) -- (.25, 1);
\draw (.5, 0) -- (.5, 1);
\draw (.75, 0) -- (.75, 1);
\draw (1, 0) -- (1, 1);
\draw (.624,0) -- (.624, .66);
\end{tikzpicture},
\begin{tikzpicture}[scale=1.5]
\filldraw[draw=red, ultra thick, fill=red!60] (.562,.33) -- (1, .33) -- (1,1) -- (.562,1) -- cycle;
\draw[red, ultra thick] (.562,.66) -- (1, .66);
\draw[red, ultra thick] (.562,.44) -- (1, .44);
\draw[red, ultra thick] (.686,.33) -- (.686, 1);
\draw[red, ultra thick] (.75,.33) -- (.75, 1);
\draw (0,0) -- (1.25,0) -- (1.25,1) -- (0,1) -- cycle;
\draw (0,.55) -- (.75, .55);
\draw (.5, .44) -- (1.25, .44);
\draw (.562, .33) -- (.562, 1);
\draw (.686, .33) -- (.686, 1);
\draw (0,.33) -- (1.25, .33);
\draw (0,.66) -- (1.25, .66);
\draw (.25, 0) -- (.25, 1);
\draw (.5, 0) -- (.5, 1);
\draw (.75, 0) -- (.75, 1);
\draw (1, 0) -- (1, 1);
\draw (.624,0) -- (.624, .66);
\end{tikzpicture},
\begin{tikzpicture}[scale=1.5]
\filldraw[draw=red, ultra thick, fill=red!60] (.686,.33) -- (1.25, .33) -- (1.25,1) -- (.686,1) -- cycle;
\draw[red, ultra thick] (.686,.66) -- (1.25, .66);
\draw[red, ultra thick] (.686,.44) -- (1.25, .44);
\draw[red, ultra thick] (.75,.33) -- (.75, 1);
\draw[red, ultra thick] (1,.33) -- (1, 1);
\draw (0,0) -- (1.25,0) -- (1.25,1) -- (0,1) -- cycle;
\draw (0,.55) -- (.75, .55);
\draw (.5, .44) -- (1.25, .44);
\draw (.562, .33) -- (.562, 1);
\draw (.686, .33) -- (.686, 1);
\draw (0,.33) -- (1.25, .33);
\draw (0,.66) -- (1.25, .66);
\draw (.25, 0) -- (.25, 1);
\draw (.5, 0) -- (.5, 1);
\draw (.75, 0) -- (.75, 1);
\draw (1, 0) -- (1, 1);
\draw (.624,0) -- (.624, .66);
\end{tikzpicture},
\begin{tikzpicture}[scale=1.5]
\filldraw[draw=red, ultra thick, fill=red!60] (.5,.44) -- (.75, .44) -- (.75,1) -- (.5,1) -- cycle;
\draw[red, ultra thick] (.5,.55) -- (.75, .55);
\draw[red, ultra thick] (.5,.66) -- (.75, .66);
\draw[red, ultra thick] (.562,.44) -- (.562, 1);
\draw[red, ultra thick] (.686,.44) -- (.686, 1);
\draw (0,0) -- (1.25,0) -- (1.25,1) -- (0,1) -- cycle;
\draw (0,.55) -- (.75, .55);
\draw (.5, .44) -- (1.25, .44);
\draw (.562, .33) -- (.562, 1);
\draw (.686, .33) -- (.686, 1);
\draw (0,.33) -- (1.25, .33);
\draw (0,.66) -- (1.25, .66);
\draw (.25, 0) -- (.25, 1);
\draw (.5, 0) -- (.5, 1);
\draw (.75, 0) -- (.75, 1);
\draw (1, 0) -- (1, 1);
\draw (.624,0) -- (.624, .66);
\end{tikzpicture}
\end{tabular}\vspace{-.15cm}
\end{minipage}\hspace{1.25cm}
}\\
\subfloat[]{
\hspace{-.5cm}
\begin{tabular}{c}
\begin{tikzpicture}[scale=1.5]
\filldraw[draw=red, ultra thick, fill=red!60] (.25,0) -- (1, 0) -- (1,1) -- (.25,1) -- cycle;
\draw[red, ultra thick] (.25,.33) -- (1, .33);
\draw[red, ultra thick] (.25,.66) -- (1, .66);
\draw[red, ultra thick] (.5,0) -- (.5, 1);
\draw[red, ultra thick] (.75,0) -- (.75, 1);
\draw (0,0) -- (1.25,0) -- (1.25,1) -- (0,1) -- cycle;
\draw (0,.55) -- (.75, .55);
\draw (.5, .44) -- (1.25, .44);
\draw (.562, .33) -- (.562, 1);
\draw (.686, .33) -- (.686, 1);
\draw (0,.33) -- (1.25, .33);
\draw (0,.66) -- (1.25, .66);
\draw (.25, 0) -- (.25, 1);
\draw (.5, 0) -- (.5, 1);
\draw (.75, 0) -- (.75, 1);
\draw (1, 0) -- (1, 1);
\draw (.624,0) -- (.624, .66);
\end{tikzpicture},
\begin{tikzpicture}[scale=1.5]
\filldraw[draw=red, ultra thick, fill=red!60] (.25,0) -- (.75, 0) -- (.75,.66) -- (.25,.66) -- cycle;
\draw[red, ultra thick] (.25,.33) -- (.75, .33);
\draw[red, ultra thick] (.25,.55) -- (.75, .55);
\draw[red, ultra thick] (.5,0) -- (.5, .66);
\draw[red, ultra thick] (.624,0) -- (.624, .66);
\draw (0,0) -- (1.25,0) -- (1.25,1) -- (0,1) -- cycle;
\draw (0,.55) -- (.75, .55);
\draw (.5, .44) -- (1.25, .44);
\draw (.562, .33) -- (.562, 1);
\draw (.686, .33) -- (.686, 1);
\draw (0,.33) -- (1.25, .33);
\draw (0,.66) -- (1.25, .66);
\draw (.25, 0) -- (.25, 1);
\draw (.5, 0) -- (.5, 1);
\draw (.75, 0) -- (.75, 1);
\draw (1, 0) -- (1, 1);
\draw (.624,0) -- (.624, .66);
\end{tikzpicture},
\begin{tikzpicture}[scale=1.5]
\filldraw[draw=red, ultra thick, fill=red!60] (.5,0) -- (1, 0) -- (1,.66) -- (.5,.66) -- cycle;
\draw[red, ultra thick] (.5,.33) -- (1, .33);
\draw[red, ultra thick] (.5,.44) -- (1, .44);
\draw[red, ultra thick] (.625,0) -- (.625, .66);
\draw[red, ultra thick] (.75,0) -- (.75, .66);
\draw (0,0) -- (1.25,0) -- (1.25,1) -- (0,1) -- cycle;
\draw (0,.55) -- (.75, .55);
\draw (.5, .44) -- (1.25, .44);
\draw (.562, .33) -- (.562, 1);
\draw (.686, .33) -- (.686, 1);
\draw (0,.33) -- (1.25, .33);
\draw (0,.66) -- (1.25, .66);
\draw (.25, 0) -- (.25, 1);
\draw (.5, 0) -- (.5, 1);
\draw (.75, 0) -- (.75, 1);
\draw (1, 0) -- (1, 1);
\draw (.624,0) -- (.624, .66);
\end{tikzpicture},
\begin{tikzpicture}[scale=1.5]
\filldraw[draw=red, ultra thick, fill=red!60] (.5,.44) -- (.75, .44) -- (.75,1) -- (.5,1) -- cycle;
\draw[red, ultra thick] (.5,.55) -- (.75, .55);
\draw[red, ultra thick] (.5,.66) -- (.75, .66);
\draw[red, ultra thick] (.562,.44) -- (.562, 1);
\draw[red, ultra thick] (.686,.44) -- (.686, 1);
\draw (0,0) -- (1.25,0) -- (1.25,1) -- (0,1) -- cycle;
\draw (0,.55) -- (.75, .55);
\draw (.5, .44) -- (1.25, .44);
\draw (.562, .33) -- (.562, 1);
\draw (.686, .33) -- (.686, 1);
\draw (0,.33) -- (1.25, .33);
\draw (0,.66) -- (1.25, .66);
\draw (.25, 0) -- (.25, 1);
\draw (.5, 0) -- (.5, 1);
\draw (.75, 0) -- (.75, 1);
\draw (1, 0) -- (1, 1);
\draw (.624,0) -- (.624, .66);
\end{tikzpicture},
\begin{tikzpicture}[scale=1.5]
\filldraw[draw=red, ultra thick, fill=red!60] (.25,.33) -- (.686, .33) -- (.686,1) -- (.25,1) -- cycle;
\draw[red, ultra thick] (.25,.55) -- (.686, .55);
\draw[red, ultra thick] (.25,.66) -- (.686, .66);
\draw[red, ultra thick] (.5,.33) -- (.5, 1);
\draw[red, ultra thick] (.562,.33) -- (.562, 1);
\draw (0,0) -- (1.25,0) -- (1.25,1) -- (0,1) -- cycle;
\draw (0,.55) -- (.75, .55);
\draw (.5, .44) -- (1.25, .44);
\draw (.562, .33) -- (.562, 1);
\draw (.686, .33) -- (.686, 1);
\draw (0,.33) -- (1.25, .33);
\draw (0,.66) -- (1.25, .66);
\draw (.25, 0) -- (.25, 1);
\draw (.5, 0) -- (.5, 1);
\draw (.75, 0) -- (.75, 1);
\draw (1, 0) -- (1, 1);
\draw (.624,0) -- (.624, .66);
\end{tikzpicture},
\begin{tikzpicture}[scale=1.5]
\filldraw[draw=red, ultra thick, fill=red!60] (.5,.33) -- (.686, .33) -- (.686,.66) -- (.5,.66) -- cycle;
\draw[red, ultra thick] (.5,.55) -- (.686, .55);
\draw[red, ultra thick] (.5,.44) -- (.686, .44);
\draw[red, ultra thick] (.562,.33) -- (.562, .66);
\draw[red, ultra thick] (.624,.33) -- (.624, .66);
\draw (0,0) -- (1.25,0) -- (1.25,1) -- (0,1) -- cycle;
\draw (0,.55) -- (.75, .55);
\draw (.5, .44) -- (1.25, .44);
\draw (.562, .33) -- (.562, 1);
\draw (.686, .33) -- (.686, 1);
\draw (0,.33) -- (1.25, .33);
\draw (0,.66) -- (1.25, .66);
\draw (.25, 0) -- (.25, 1);
\draw (.5, 0) -- (.5, 1);
\draw (.75, 0) -- (.75, 1);
\draw (1, 0) -- (1, 1);
\draw (.624,0) -- (.624, .66);
\end{tikzpicture},
\begin{tikzpicture}[scale=1.5]
\filldraw[draw=red, ultra thick, fill=red!60] (.562,.33) -- (.75, .33) -- (.75,.66) -- (.562,.66) -- cycle;
\draw[red, ultra thick] (.562,.55) -- (.75, .55);
\draw[red, ultra thick] (.562,.44) -- (.75, .44);
\draw[red, ultra thick] (.686,.33) -- (.686, .66);
\draw[red, ultra thick] (.624,.33) -- (.624, .66);
\draw (0,0) -- (1.25,0) -- (1.25,1) -- (0,1) -- cycle;
\draw (0,.55) -- (.75, .55);
\draw (.5, .44) -- (1.25, .44);
\draw (.562, .33) -- (.562, 1);
\draw (.686, .33) -- (.686, 1);
\draw (0,.33) -- (1.25, .33);
\draw (0,.66) -- (1.25, .66);
\draw (.25, 0) -- (.25, 1);
\draw (.5, 0) -- (.5, 1);
\draw (.75, 0) -- (.75, 1);
\draw (1, 0) -- (1, 1);
\draw (.624,0) -- (.624, .66);
\end{tikzpicture},
\begin{tikzpicture}[scale=1.5]
\filldraw[draw=red, ultra thick, fill=red!60] (.562,.33) -- (1, .33) -- (1,1) -- (.562,1) -- cycle;
\draw[red, ultra thick] (.562,.66) -- (1, .66);
\draw[red, ultra thick] (.562,.44) -- (1, .44);
\draw[red, ultra thick] (.686,.33) -- (.686, 1);
\draw[red, ultra thick] (.75,.33) -- (.75, 1);
\draw (0,0) -- (1.25,0) -- (1.25,1) -- (0,1) -- cycle;
\draw (0,.55) -- (.75, .55);
\draw (.5, .44) -- (1.25, .44);
\draw (.562, .33) -- (.562, 1);
\draw (.686, .33) -- (.686, 1);
\draw (0,.33) -- (1.25, .33);
\draw (0,.66) -- (1.25, .66);
\draw (.25, 0) -- (.25, 1);
\draw (.5, 0) -- (.5, 1);
\draw (.75, 0) -- (.75, 1);
\draw (1, 0) -- (1, 1);
\draw (.624,0) -- (.624, .66);
\end{tikzpicture}
\end{tabular}\vspace{-.15cm}
}
\caption{(a) the LR-mesh of multiplicity 1, (b) the LR B-splines of degree (2,2) on it, (c) the LR B-splines in a linear dependence relation.}\label{8LR}
\end{figure}
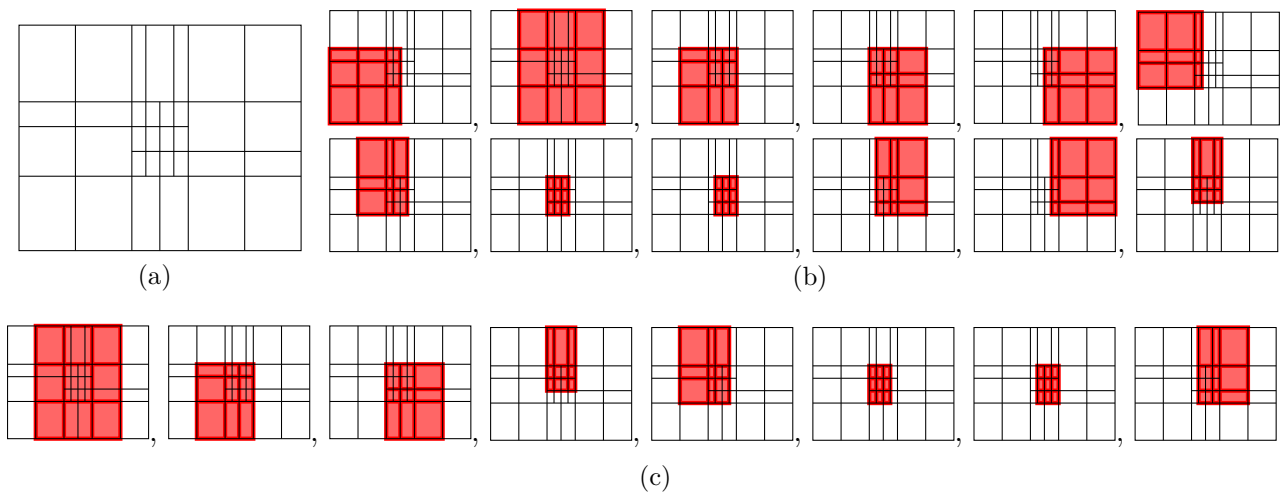
For what concerns general bidegrees $(p_1,p_2)$ with $p_k \geq 2$ for some $k \in \{1,2\}$, it is always possible to arrange the LR B-splines in the same way as for bidegree $(2,2)$. For instance, in Figure \ref{deg34LR} are reported the cases for $(p_1,p_2) =(3,3), (4,4), (3,1), (2,0)$. Also here $\dim \SSS_{(p_1,p_2)} (\cM+\gamma) = 11$ while $|\mathcal{B}^{\mathcal{LR}}(\cM+\gamma)| = 12$.
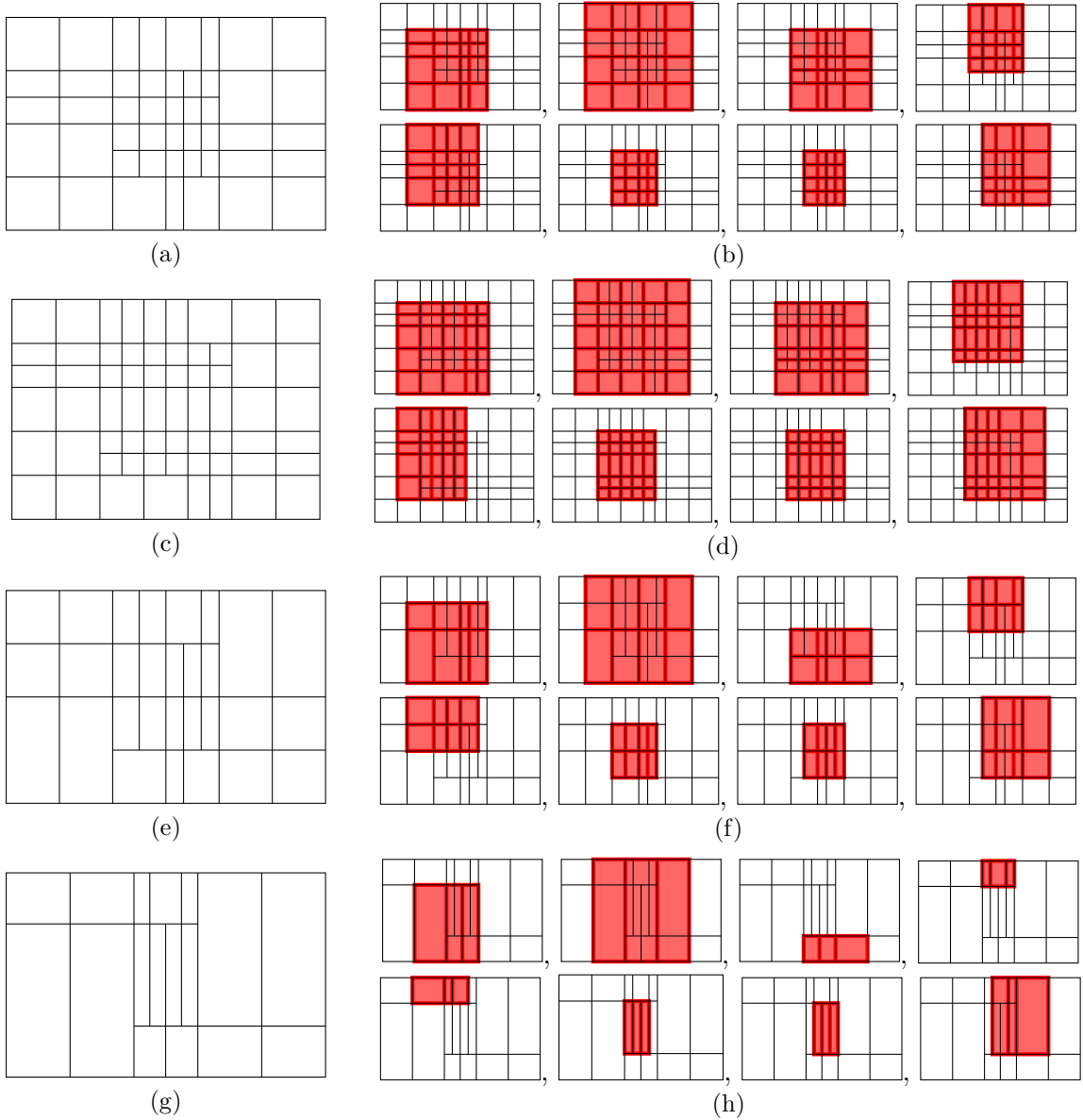
\begin{figure}[h!]\centering
\subfloat[]{
\begin{tikzpicture}[scale = 3]
\draw (0,0) -- (1.5,0) -- (1.5, 1) -- (0,1) -- cycle;
\draw (0,.25) -- (1.5, .25);
\draw (0,.5) -- (1.5, .5);
\draw (0,.75) -- (1.5,.75);
\draw (.25,0) -- (.25, 1);
\draw (.5,0) -- (.5, 1);
\draw (.75,0) -- (.75, 1);
\draw (1,0) -- (1, 1);
\draw (1.25,0) -- (1.25, 1);
\draw (0,.625) -- (1,.625);
\draw (.5,.375) -- (1.5,.375);
\draw (.625,1) -- (.625, .25);
\draw (.833,0) -- (.833,.75);
\draw (.916,.25) -- (.916,1);
\end{tikzpicture}
}~
\subfloat[]{
\begin{minipage}{8.5cm}\vspace{-3cm}
\begin{tabular}{c}
\begin{tikzpicture}[scale=1.5]
\filldraw[draw = red, ultra thick, fill=red!60] (.25,0) -- (1,0) -- (1,.75) -- (.25,.75) -- cycle;
\draw[red, ultra thick] (.25,.25) -- (1,.25);
\draw[red, ultra thick] (.25,.5) -- (1,.5);
\draw[red, ultra thick] (.25,.625) -- (1,.625);
\draw[red, ultra thick] (.5,0) -- (.5,.75);
\draw[red, ultra thick] (.75,0) -- (.75,.75);
\draw[red, ultra thick] (.833,0) -- (.833,.75);
\draw (0,0) -- (1.5,0) -- (1.5, 1) -- (0,1) -- cycle;
\draw (0,.25) -- (1.5, .25);
\draw (0,.5) -- (1.5, .5);
\draw (0,.75) -- (1.5,.75);
\draw (.25,0) -- (.25, 1);
\draw (.5,0) -- (.5, 1);
\draw (.75,0) -- (.75, 1);
\draw (1,0) -- (1, 1);
\draw (1.25,0) -- (1.25, 1);
\draw (0,.625) -- (1,.625);
\draw (.5,.375) -- (1.5,.375);
\draw (.625,1) -- (.625, .25);
\draw (.833,0) -- (.833,.75);
\draw (.916,.25) -- (.916,1);
\end{tikzpicture},
\begin{tikzpicture}[scale=1.5]
\filldraw[draw = red, ultra thick, fill=red!60] (.25,0) -- (1.25,0) -- (1.25,1) -- (.25,1) -- cycle;
\draw[red, ultra thick] (.25,.25) -- (1.25,.25);
\draw[red, ultra thick] (.25,.5) -- (1.25,.5);
\draw[red, ultra thick] (.25,.75) -- (1.25,.75);
\draw[red, ultra thick] (.5,0) -- (.5,1);
\draw[red, ultra thick] (.75,0) -- (.75,1);
\draw[red, ultra thick] (1,0) -- (1,1);
\draw (0,0) -- (1.5,0) -- (1.5, 1) -- (0,1) -- cycle;
\draw (0,.25) -- (1.5, .25);
\draw (0,.5) -- (1.5, .5);
\draw (0,.75) -- (1.5,.75);
\draw (.25,0) -- (.25, 1);
\draw (.5,0) -- (.5, 1);
\draw (.75,0) -- (.75, 1);
\draw (1,0) -- (1, 1);
\draw (1.25,0) -- (1.25, 1);
\draw (0,.625) -- (1,.625);
\draw (.5,.375) -- (1.5,.375);
\draw (.625,1) -- (.625, .25);
\draw (.833,0) -- (.833,.75);
\draw (.916,.25) -- (.916,1);
\end{tikzpicture},
\begin{tikzpicture}[scale=1.5]
\filldraw[draw = red, ultra thick, fill=red!60] (.5,0) -- (1.25,0) -- (1.25,.75) -- (.5,.75) -- cycle;
\draw[red, ultra thick] (.5,.25) -- (1.25,.25);
\draw[red, ultra thick] (.5,.375) -- (1.25,.375);
\draw[red, ultra thick] (.5,.5) -- (1.25,.5);
\draw[red, ultra thick] (.75,0) -- (.75,.75);
\draw[red, ultra thick] (.833,0) -- (.833,.75);
\draw[red, ultra thick] (1,0) -- (1,.75);
\draw (0,0) -- (1.5,0) -- (1.5, 1) -- (0,1) -- cycle;
\draw (0,.25) -- (1.5, .25);
\draw (0,.5) -- (1.5, .5);
\draw (0,.75) -- (1.5,.75);
\draw (.25,0) -- (.25, 1);
\draw (.5,0) -- (.5, 1);
\draw (.75,0) -- (.75, 1);
\draw (1,0) -- (1, 1);
\draw (1.25,0) -- (1.25, 1);
\draw (0,.625) -- (1,.625);
\draw (.5,.375) -- (1.5,.375);
\draw (.625,1) -- (.625, .25);
\draw (.833,0) -- (.833,.75);
\draw (.916,.25) -- (.916,1);
\end{tikzpicture},
\begin{tikzpicture}[scale=1.5]
\filldraw[draw = red, ultra thick, fill=red!60] (.5,.375) -- (1,.375) -- (1,1) -- (.5,1) -- cycle;
\draw[red, ultra thick] (.5,.5) -- (1,.5);
\draw[red, ultra thick] (.5,.625) -- (1,.625);
\draw[red, ultra thick] (.5,.75) -- (1,.75);
\draw[red, ultra thick] (.625,.375) -- (.625,1);
\draw[red, ultra thick] (.75,.375) -- (.75,1);
\draw[red, ultra thick] (.916,.375) -- (.916,1);
\draw (0,0) -- (1.5,0) -- (1.5, 1) -- (0,1) -- cycle;
\draw (0,.25) -- (1.5, .25);
\draw (0,.5) -- (1.5, .5);
\draw (0,.75) -- (1.5,.75);
\draw (.25,0) -- (.25, 1);
\draw (.5,0) -- (.5, 1);
\draw (.75,0) -- (.75, 1);
\draw (1,0) -- (1, 1);
\draw (1.25,0) -- (1.25, 1);
\draw (0,.625) -- (1,.625);
\draw (.5,.375) -- (1.5,.375);
\draw (.625,1) -- (.625, .25);
\draw (.833,0) -- (.833,.75);
\draw (.916,.25) -- (.916,1);
\end{tikzpicture}
\\
\begin{tikzpicture}[scale=1.5]
\filldraw[draw = red, ultra thick, fill=red!60] (.25,.25) -- (.916,.25) -- (.916,1) -- (.25,1) -- cycle;
\draw[red, ultra thick] (.25,.5) -- (.916,.5);
\draw[red, ultra thick] (.25,.625) -- (.916,.625);
\draw[red, ultra thick] (.25,.75) -- (.916,.75);
\draw[red, ultra thick] (.5,.25) -- (.5,1);
\draw[red, ultra thick] (.625,.25) -- (.625,1);
\draw[red, ultra thick] (.75,.25) -- (.75,1);
\draw (0,0) -- (1.5,0) -- (1.5, 1) -- (0,1) -- cycle;
\draw (0,.25) -- (1.5, .25);
\draw (0,.5) -- (1.5, .5);
\draw (0,.75) -- (1.5,.75);
\draw (.25,0) -- (.25, 1);
\draw (.5,0) -- (.5, 1);
\draw (.75,0) -- (.75, 1);
\draw (1,0) -- (1, 1);
\draw (1.25,0) -- (1.25, 1);
\draw (0,.625) -- (1,.625);
\draw (.5,.375) -- (1.5,.375);
\draw (.625,1) -- (.625, .25);
\draw (.833,0) -- (.833,.75);
\draw (.916,.25) -- (.916,1);
\end{tikzpicture},
\begin{tikzpicture}[scale=1.5]
\filldraw[draw = red, ultra thick, fill=red!60] (.5,.25) -- (.916,.25) -- (.916,.75) -- (.5,.75) -- cycle;
\draw[red, ultra thick] (.5,.375) -- (.916,.375);
\draw[red, ultra thick] (.5,.5) -- (.916,.5);
\draw[red, ultra thick] (.5,.625) -- (.916,.625);
\draw[red, ultra thick] (.625,.25) -- (.625,.75);
\draw[red, ultra thick] (.75,.25) -- (.75,.75);
\draw[red, ultra thick] (.833,.25) -- (.833,.75);
\draw (0,0) -- (1.5,0) -- (1.5, 1) -- (0,1) -- cycle;
\draw (0,.25) -- (1.5, .25);
\draw (0,.5) -- (1.5, .5);
\draw (0,.75) -- (1.5,.75);
\draw (.25,0) -- (.25, 1);
\draw (.5,0) -- (.5, 1);
\draw (.75,0) -- (.75, 1);
\draw (1,0) -- (1, 1);
\draw (1.25,0) -- (1.25, 1);
\draw (0,.625) -- (1,.625);
\draw (.5,.375) -- (1.5,.375);
\draw (.625,1) -- (.625, .25);
\draw (.833,0) -- (.833,.75);
\draw (.916,.25) -- (.916,1);
\end{tikzpicture},
\begin{tikzpicture}[scale=1.5]
\filldraw[draw = red, ultra thick, fill=red!60] (.625,.25) -- (1,.25) -- (1,.75) -- (.625,.75) -- cycle;
\draw[red, ultra thick] (.625,.375) -- (1,.375);
\draw[red, ultra thick] (.625,.5) -- (1,.5);
\draw[red, ultra thick] (.625,.625) -- (1,.625);
\draw[red, ultra thick] (.916,.25) -- (.916,.75);
\draw[red, ultra thick] (.75,.25) -- (.75,.75);
\draw[red, ultra thick] (.833,.25) -- (.833,.75);
\draw (0,0) -- (1.5,0) -- (1.5, 1) -- (0,1) -- cycle;
\draw (0,.25) -- (1.5, .25);
\draw (0,.5) -- (1.5, .5);
\draw (0,.75) -- (1.5,.75);
\draw (.25,0) -- (.25, 1);
\draw (.5,0) -- (.5, 1);
\draw (.75,0) -- (.75, 1);
\draw (1,0) -- (1, 1);
\draw (1.25,0) -- (1.25, 1);
\draw (0,.625) -- (1,.625);
\draw (.5,.375) -- (1.5,.375);
\draw (.625,1) -- (.625, .25);
\draw (.833,0) -- (.833,.75);
\draw (.916,.25) -- (.916,1);
\end{tikzpicture},
\begin{tikzpicture}[scale=1.5]
\filldraw[draw = red, ultra thick, fill=red!60] (.625,.25) -- (1.25,.25) -- (1.25,1) -- (.625,1) -- cycle;
\draw[red, ultra thick] (.625,.375) -- (1.25,.375);
\draw[red, ultra thick] (.625,.5) -- (1.25,.5);
\draw[red, ultra thick] (.625,.75) -- (1.25,.75);
\draw[red, ultra thick] (.916,.25) -- (.916,1);
\draw[red, ultra thick] (.75,.25) -- (.75,1);
\draw[red, ultra thick] (1,.25) -- (1,1);
\draw (0,0) -- (1.5,0) -- (1.5, 1) -- (0,1) -- cycle;
\draw (0,.25) -- (1.5, .25);
\draw (0,.5) -- (1.5, .5);
\draw (0,.75) -- (1.5,.75);
\draw (.25,0) -- (.25, 1);
\draw (.5,0) -- (.5, 1);
\draw (.75,0) -- (.75, 1);
\draw (1,0) -- (1, 1);
\draw (1.25,0) -- (1.25, 1);
\draw (0,.625) -- (1,.625);
\draw (.5,.375) -- (1.5,.375);
\draw (.625,1) -- (.625, .25);
\draw (.833,0) -- (.833,.75);
\draw (.916,.25) -- (.916,1);
\end{tikzpicture}
\end{tabular}\vspace{-.15cm}
\end{minipage}\hspace{1.75cm}
}\\
\subfloat[]{
\begin{tikzpicture}[scale=3.1]
\draw (0,0) -- (1.4,0) -- (1.4,1) -- (0,1) -- cycle;
\draw (0,.2) -- (1.4,.2);
\draw (0,.4) -- (1.4,.4);
\draw (0,.6) -- (1.4,.6);
\draw (0,.8) -- (1.4,.8);
\draw (.2,0) -- (.2,1);
\draw (.4,0) -- (.4,1);
\draw (.6,0) -- (.6,1);
\draw (.8,0) -- (.8,1);
\draw (1,0) -- (1,1);
\draw (1.2,0) -- (1.2,1);
\draw (0,.7) -- (1,.7);
\draw (.4,.3) -- (1.4,.3);
\draw (.5,.2) -- (.5,1);
\draw (.9,0) -- (.9,.8);
\draw (.7, .2) -- (.7,1);
\end{tikzpicture}
}~
\subfloat[]{
\begin{minipage}{8.5cm}\vspace{-3.15cm}
\begin{tabular}{c}
\begin{tikzpicture}[scale=1.6]
\filldraw[draw = red, ultra thick, fill=red!60] (.2,0) -- (1,0) -- (1,.8) -- (.2,.8) -- cycle;
\draw[red, ultra thick] (.2,.2) -- (1,.2);
\draw[red, ultra thick] (.2,.4) -- (1,.4);
\draw[red, ultra thick] (.2,.6) -- (1,.6);
\draw[red, ultra thick] (.2,.7) -- (1,.7);
\draw[red, ultra thick] (.4,0) -- (.4,.8);
\draw[red, ultra thick] (.6,0) -- (.6,.8);
\draw[red, ultra thick] (.8,0) -- (.8,.8);
\draw[red, ultra thick] (.9,0) -- (.9,.8);
\draw (0,0) -- (1.4,0) -- (1.4,1) -- (0,1) -- cycle;
\draw (0,.2) -- (1.4,.2);
\draw (0,.4) -- (1.4,.4);
\draw (0,.6) -- (1.4,.6);
\draw (0,.8) -- (1.4,.8);
\draw (.2,0) -- (.2,1);
\draw (.4,0) -- (.4,1);
\draw (.6,0) -- (.6,1);
\draw (.8,0) -- (.8,1);
\draw (1,0) -- (1,1);
\draw (1.2,0) -- (1.2,1);
\draw (0,.7) -- (1,.7);
\draw (.4,.3) -- (1.4,.3);
\draw (.5,.2) -- (.5,1);
\draw (.9,0) -- (.9,.8);
\draw (.7, .2) -- (.7,1);
\end{tikzpicture},
\begin{tikzpicture}[scale=1.6]
\filldraw[draw = red, ultra thick, fill=red!60] (.2,0) -- (1.2,0) -- (1.2,1) -- (.2,1) -- cycle;
\draw[red, ultra thick] (.2,.2) -- (1.2,.2);
\draw[red, ultra thick] (.2,.4) -- (1.2,.4);
\draw[red, ultra thick] (.2,.6) -- (1.2,.6);
\draw[red, ultra thick] (.2,.8) -- (1.2,.8);
\draw[red, ultra thick] (.4,0) -- (.4,1);
\draw[red, ultra thick] (.6,0) -- (.6,1);
\draw[red, ultra thick] (.8,0) -- (.8,1);
\draw[red, ultra thick] (1,0) -- (1,1);
\draw (0,0) -- (1.4,0) -- (1.4,1) -- (0,1) -- cycle;
\draw (0,.2) -- (1.4,.2);
\draw (0,.4) -- (1.4,.4);
\draw (0,.6) -- (1.4,.6);
\draw (0,.8) -- (1.4,.8);
\draw (.2,0) -- (.2,1);
\draw (.4,0) -- (.4,1);
\draw (.6,0) -- (.6,1);
\draw (.8,0) -- (.8,1);
\draw (1,0) -- (1,1);
\draw (1.2,0) -- (1.2,1);
\draw (0,.7) -- (1,.7);
\draw (.4,.3) -- (1.4,.3);
\draw (.5,.2) -- (.5,1);
\draw (.9,0) -- (.9,.8);
\draw (.7, .2) -- (.7,1);
\end{tikzpicture},
\begin{tikzpicture}[scale=1.6]
\filldraw[draw = red, ultra thick, fill=red!60] (.4,0) -- (1.2,0) -- (1.2,.8) -- (.4,.8) -- cycle;
\draw[red, ultra thick] (.4,.2) -- (1.2,.2);
\draw[red, ultra thick] (.4,.3) -- (1.2,.3);
\draw[red, ultra thick] (.4,.4) -- (1.2,.4);
\draw[red, ultra thick] (.4,.6) -- (1.2,.6);
\draw[red, ultra thick] (.6,0) -- (.6,.8);
\draw[red, ultra thick] (.8,0) -- (.8,.8);
\draw[red, ultra thick] (.9,0) -- (.9,.8);
\draw[red, ultra thick] (1,0) -- (1,.8);
\draw (0,0) -- (1.4,0) -- (1.4,1) -- (0,1) -- cycle;
\draw (0,.2) -- (1.4,.2);
\draw (0,.4) -- (1.4,.4);
\draw (0,.6) -- (1.4,.6);
\draw (0,.8) -- (1.4,.8);
\draw (.2,0) -- (.2,1);
\draw (.4,0) -- (.4,1);
\draw (.6,0) -- (.6,1);
\draw (.8,0) -- (.8,1);
\draw (1,0) -- (1,1);
\draw (1.2,0) -- (1.2,1);
\draw (0,.7) -- (1,.7);
\draw (.4,.3) -- (1.4,.3);
\draw (.5,.2) -- (.5,1);
\draw (.9,0) -- (.9,.8);
\draw (.7, .2) -- (.7,1);
\end{tikzpicture},
\begin{tikzpicture}[scale=1.6]
\filldraw[draw = red, ultra thick, fill=red!60] (.4,.3) -- (1,.3) -- (1,1) -- (.4,1) -- cycle;
\draw[red, ultra thick] (.4,.8) -- (1,.8);
\draw[red, ultra thick] (.4,.4) -- (1,.4);
\draw[red, ultra thick] (.4,.6) -- (1,.6);
\draw[red, ultra thick] (.4,.7) -- (1,.7);
\draw[red, ultra thick] (.5,.3) -- (.5,1);
\draw[red, ultra thick] (.6,.3) -- (.6,1);
\draw[red, ultra thick] (.8,.3) -- (.8,1);
\draw[red, ultra thick] (.7,.3) -- (.7,1);
\draw (0,0) -- (1.4,0) -- (1.4,1) -- (0,1) -- cycle;
\draw (0,.2) -- (1.4,.2);
\draw (0,.4) -- (1.4,.4);
\draw (0,.6) -- (1.4,.6);
\draw (0,.8) -- (1.4,.8);
\draw (.2,0) -- (.2,1);
\draw (.4,0) -- (.4,1);
\draw (.6,0) -- (.6,1);
\draw (.8,0) -- (.8,1);
\draw (1,0) -- (1,1);
\draw (1.2,0) -- (1.2,1);
\draw (0,.7) -- (1,.7);
\draw (.4,.3) -- (1.4,.3);
\draw (.5,.2) -- (.5,1);
\draw (.9,0) -- (.9,.8);
\draw (.7, .2) -- (.7,1);
\end{tikzpicture}\\
\begin{tikzpicture}[scale=1.6]
\filldraw[draw = red, ultra thick, fill=red!60] (.2,.2) -- (.8,.2) -- (.8,1) -- (.2,1) -- cycle;
\draw[red, ultra thick] (.2,.8) -- (.8,.8);
\draw[red, ultra thick] (.2,.4) -- (.8,.4);
\draw[red, ultra thick] (.2,.6) -- (.8,.6);
\draw[red, ultra thick] (.2,.7) -- (.8,.7);
\draw[red, ultra thick] (.4,.2) -- (.4,1);
\draw[red, ultra thick] (.6,.2) -- (.6,1);
\draw[red, ultra thick] (.7,.2) -- (.7,1);
\draw[red, ultra thick] (.5,.2) -- (.5,1);
\draw (0,0) -- (1.4,0) -- (1.4,1) -- (0,1) -- cycle;
\draw (0,.2) -- (1.4,.2);
\draw (0,.4) -- (1.4,.4);
\draw (0,.6) -- (1.4,.6);
\draw (0,.8) -- (1.4,.8);
\draw (.2,0) -- (.2,1);
\draw (.4,0) -- (.4,1);
\draw (.6,0) -- (.6,1);
\draw (.8,0) -- (.8,1);
\draw (1,0) -- (1,1);
\draw (1.2,0) -- (1.2,1);
\draw (0,.7) -- (1,.7);
\draw (.4,.3) -- (1.4,.3);
\draw (.5,.2) -- (.5,1);
\draw (.9,0) -- (.9,.8);
\draw (.7, .2) -- (.7,1);
\end{tikzpicture},
\begin{tikzpicture}[scale=1.6]
\filldraw[draw = red, ultra thick, fill=red!60] (.4,.2) -- (.9,.2) -- (.9,.8) -- (.4,.8) -- cycle;
\draw[red, ultra thick] (.4,.4) -- (.9,.4);
\draw[red, ultra thick] (.4,.3) -- (.9,.3);
\draw[red, ultra thick] (.4,.6) -- (.9,.6);
\draw[red, ultra thick] (.4,.7) -- (.9,.7);
\draw[red, ultra thick] (.5,.2) -- (.5,.8);
\draw[red, ultra thick] (.6,.2) -- (.6,.8);
\draw[red, ultra thick] (.8,.2) -- (.8,.8);
\draw[red, ultra thick] (.7,.2) -- (.7,.8);
\draw (0,0) -- (1.4,0) -- (1.4,1) -- (0,1) -- cycle;
\draw (0,.2) -- (1.4,.2);
\draw (0,.4) -- (1.4,.4);
\draw (0,.6) -- (1.4,.6);
\draw (0,.8) -- (1.4,.8);
\draw (.2,0) -- (.2,1);
\draw (.4,0) -- (.4,1);
\draw (.6,0) -- (.6,1);
\draw (.8,0) -- (.8,1);
\draw (1,0) -- (1,1);
\draw (1.2,0) -- (1.2,1);
\draw (0,.7) -- (1,.7);
\draw (.4,.3) -- (1.4,.3);
\draw (.5,.2) -- (.5,1);
\draw (.9,0) -- (.9,.8);
\draw (.7, .2) -- (.7,1);
\end{tikzpicture},
\begin{tikzpicture}[scale=1.6]
\filldraw[draw = red, ultra thick, fill=red!60] (.5,.2) -- (1,.2) -- (1,.8) -- (.5,.8) -- cycle;
\draw[red, ultra thick] (.5,.4) -- (1,.4);
\draw[red, ultra thick] (.5,.3) -- (1,.3);
\draw[red, ultra thick] (.5,.6) -- (1,.6);
\draw[red, ultra thick] (.5,.7) -- (1,.7);
\draw[red, ultra thick] (.9,.2) -- (.9,.8);
\draw[red, ultra thick] (.6,.2) -- (.6,.8);
\draw[red, ultra thick] (.8,.2) -- (.8,.8);
\draw[red, ultra thick] (.7,.2) -- (.7,.8);
\draw (0,0) -- (1.4,0) -- (1.4,1) -- (0,1) -- cycle;
\draw (0,.2) -- (1.4,.2);
\draw (0,.4) -- (1.4,.4);
\draw (0,.6) -- (1.4,.6);
\draw (0,.8) -- (1.4,.8);
\draw (.2,0) -- (.2,1);
\draw (.4,0) -- (.4,1);
\draw (.6,0) -- (.6,1);
\draw (.8,0) -- (.8,1);
\draw (1,0) -- (1,1);
\draw (1.2,0) -- (1.2,1);
\draw (0,.7) -- (1,.7);
\draw (.4,.3) -- (1.4,.3);
\draw (.5,.2) -- (.5,1);
\draw (.9,0) -- (.9,.8);
\draw (.7, .2) -- (.7,1);
\end{tikzpicture},
\begin{tikzpicture}[scale=1.6]
\filldraw[draw = red, ultra thick, fill=red!60] (.5,.2) -- (1.2,.2) -- (1.2,1) -- (.5,1) -- cycle;
\draw[red, ultra thick] (.5,.4) -- (1.2,.4);
\draw[red, ultra thick] (.5,.3) -- (1.2,.3);
\draw[red, ultra thick] (.5,.6) -- (1.2,.6);
\draw[red, ultra thick] (.5,.8) -- (1.2,.8);
\draw[red, ultra thick] (.7,.2) -- (.7,1);
\draw[red, ultra thick] (.6,.2) -- (.6,1);
\draw[red, ultra thick] (.8,.2) -- (.8,1);
\draw[red, ultra thick] (1,.2) -- (1,1);
\draw (0,0) -- (1.4,0) -- (1.4,1) -- (0,1) -- cycle;
\draw (0,.2) -- (1.4,.2);
\draw (0,.4) -- (1.4,.4);
\draw (0,.6) -- (1.4,.6);
\draw (0,.8) -- (1.4,.8);
\draw (.2,0) -- (.2,1);
\draw (.4,0) -- (.4,1);
\draw (.6,0) -- (.6,1);
\draw (.8,0) -- (.8,1);
\draw (1,0) -- (1,1);
\draw (1.2,0) -- (1.2,1);
\draw (0,.7) -- (1,.7);
\draw (.4,.3) -- (1.4,.3);
\draw (.5,.2) -- (.5,1);
\draw (.9,0) -- (.9,.8);
\draw (.7, .2) -- (.7,1);
\end{tikzpicture}
\end{tabular}\vspace{-.15cm}
\end{minipage}\hspace{1.75cm}
}\\
\subfloat[]{
\begin{tikzpicture}[scale = 3]
\draw (0,0) -- (1.5,0) -- (1.5, 1) -- (0,1) -- cycle;
\draw (0,.5) -- (1.5, .5);
\draw (.25,0) -- (.25, 1);
\draw (.5,0) -- (.5, 1);
\draw (.75,0) -- (.75, 1);
\draw (1,0) -- (1, 1);
\draw (1.25,0) -- (1.25, 1);
\draw (0,.75) -- (1,.75);
\draw (.5,.25) -- (1.5,.25);
\draw (.625,1) -- (.625, .25);
\draw (.833,0) -- (.833,.75);
\draw (.916,.25) -- (.916,1);
\end{tikzpicture}
}~
\subfloat[]{
\begin{minipage}{6cm}\vspace{-3cm}
\begin{tabular}{c}
\begin{tikzpicture}[scale=1.5]
\filldraw[draw = red, ultra thick, fill=red!60] (.25,0) -- (1,0) -- (1,.75) -- (.25,.75) -- cycle;
\draw[red, ultra thick] (.25,.5) -- (1,.5);
\draw[red, ultra thick] (.5,0) -- (.5,.75);
\draw[red, ultra thick] (.75,0) -- (.75,.75);
\draw[red, ultra thick] (.833,0) -- (.833,.75);
\draw (0,0) -- (1.5,0) -- (1.5, 1) -- (0,1) -- cycle;
\draw (0,.5) -- (1.5, .5);
\draw (.25,0) -- (.25, 1);
\draw (.5,0) -- (.5, 1);
\draw (.75,0) -- (.75, 1);
\draw (1,0) -- (1, 1);
\draw (1.25,0) -- (1.25, 1);
\draw (0,.75) -- (1,.75);
\draw (.5,.25) -- (1.5,.25);
\draw (.625,1) -- (.625, .25);
\draw (.833,0) -- (.833,.75);
\draw (.916,.25) -- (.916,1);
\end{tikzpicture},
\begin{tikzpicture}[scale=1.5]
\filldraw[draw = red, ultra thick, fill=red!60] (.25,0) -- (1.25,0) -- (1.25,1) -- (.25,1) -- cycle;
\draw[red, ultra thick] (.25,.5) -- (1.25,.5);
\draw[red, ultra thick] (.5,0) -- (.5,1);
\draw[red, ultra thick] (.75,0) -- (.75,1);
\draw[red, ultra thick] (1,0) -- (1,1);
\draw (0,0) -- (1.5,0) -- (1.5, 1) -- (0,1) -- cycle;
\draw (0,.5) -- (1.5, .5);
\draw (.25,0) -- (.25, 1);
\draw (.5,0) -- (.5, 1);
\draw (.75,0) -- (.75, 1);
\draw (1,0) -- (1, 1);
\draw (1.25,0) -- (1.25, 1);
\draw (0,.75) -- (1,.75);
\draw (.5,.25) -- (1.5,.25);
\draw (.625,1) -- (.625, .25);
\draw (.833,0) -- (.833,.75);
\draw (.916,.25) -- (.916,1);
\end{tikzpicture},
\begin{tikzpicture}[scale=1.5]
\filldraw[draw = red, ultra thick, fill=red!60] (.5,0) -- (1.25,0) -- (1.25,.5) -- (.5,.5) -- cycle;
\draw[red, ultra thick] (.5,.25) -- (1.25,.25);
\draw[red, ultra thick] (.75,0) -- (.75,.5);
\draw[red, ultra thick] (.833,0) -- (.833,.5);
\draw[red, ultra thick] (1,0) -- (1,.5);
\draw (0,0) -- (1.5,0) -- (1.5, 1) -- (0,1) -- cycle;
\draw (0,.5) -- (1.5, .5);
\draw (.25,0) -- (.25, 1);
\draw (.5,0) -- (.5, 1);
\draw (.75,0) -- (.75, 1);
\draw (1,0) -- (1, 1);
\draw (1.25,0) -- (1.25, 1);
\draw (0,.75) -- (1,.75);
\draw (.5,.25) -- (1.5,.25);
\draw (.625,1) -- (.625, .25);
\draw (.833,0) -- (.833,.75);
\draw (.916,.25) -- (.916,1);
\end{tikzpicture},
\begin{tikzpicture}[scale=1.5]
\filldraw[draw = red, ultra thick, fill=red!60] (.5,.5) -- (1,.5) -- (1,1) -- (.5,1) -- cycle;
\draw[red, ultra thick] (.5,.75) -- (1,.75);
\draw[red, ultra thick] (.625,.5) -- (.625,1);
\draw[red, ultra thick] (.75,.5) -- (.75,1);
\draw[red, ultra thick] (.916,.5) -- (.916,1);
\draw (0,0) -- (1.5,0) -- (1.5, 1) -- (0,1) -- cycle;
\draw (0,.5) -- (1.5, .5);
\draw (.25,0) -- (.25, 1);
\draw (.5,0) -- (.5, 1);
\draw (.75,0) -- (.75, 1);
\draw (1,0) -- (1, 1);
\draw (1.25,0) -- (1.25, 1);
\draw (0,.75) -- (1,.75);
\draw (.5,.25) -- (1.5,.25);
\draw (.625,1) -- (.625, .25);
\draw (.833,0) -- (.833,.75);
\draw (.916,.25) -- (.916,1);
\end{tikzpicture}\\
\begin{tikzpicture}[scale=1.5]
\filldraw[draw = red, ultra thick, fill=red!60] (.25,.5) -- (.916,.5) -- (.916,1) -- (.25,1) -- cycle;
\draw[red, ultra thick] (.25,.75) -- (.916,.75);
\draw[red, ultra thick] (.5,.5) -- (.5,1);
\draw[red, ultra thick] (.625,.5) -- (.625,1);
\draw[red, ultra thick] (.75,.5) -- (.75,1);
\draw (0,0) -- (1.5,0) -- (1.5, 1) -- (0,1) -- cycle;
\draw (0,.5) -- (1.5, .5);
\draw (.25,0) -- (.25, 1);
\draw (.5,0) -- (.5, 1);
\draw (.75,0) -- (.75, 1);
\draw (1,0) -- (1, 1);
\draw (1.25,0) -- (1.25, 1);
\draw (0,.75) -- (1,.75);
\draw (.5,.25) -- (1.5,.25);
\draw (.625,1) -- (.625, .25);
\draw (.833,0) -- (.833,.75);
\draw (.916,.25) -- (.916,1);
\end{tikzpicture},
\begin{tikzpicture}[scale=1.5]
\filldraw[draw = red, ultra thick, fill=red!60] (.5,.25) -- (.916,.25) -- (.916,.75) -- (.5,.75) -- cycle;
\draw[red, ultra thick] (.5,.5) -- (.916,.5);
\draw[red, ultra thick] (.625,.25) -- (.625,.75);
\draw[red, ultra thick] (.75,.25) -- (.75,.75);
\draw[red, ultra thick] (.833,.25) -- (.833,.75);
\draw (0,0) -- (1.5,0) -- (1.5, 1) -- (0,1) -- cycle;
\draw (0,.5) -- (1.5, .5);
\draw (.25,0) -- (.25, 1);
\draw (.5,0) -- (.5, 1);
\draw (.75,0) -- (.75, 1);
\draw (1,0) -- (1, 1);
\draw (1.25,0) -- (1.25, 1);
\draw (0,.75) -- (1,.75);
\draw (.5,.25) -- (1.5,.25);
\draw (.625,1) -- (.625, .25);
\draw (.833,0) -- (.833,.75);
\draw (.916,.25) -- (.916,1);
\end{tikzpicture},
\begin{tikzpicture}[scale=1.5]
\filldraw[draw = red, ultra thick, fill=red!60] (.625,.25) -- (1,.25) -- (1,.75) -- (.625,.75) -- cycle;
\draw[red, ultra thick] (.625,.5) -- (.916,.5);
\draw[red, ultra thick] (.916,.25) -- (.916,.75);
\draw[red, ultra thick] (.75,.25) -- (.75,.75);
\draw[red, ultra thick] (.833,.25) -- (.833,.75);
\draw (0,0) -- (1.5,0) -- (1.5, 1) -- (0,1) -- cycle;
\draw (0,.5) -- (1.5, .5);
\draw (.25,0) -- (.25, 1);
\draw (.5,0) -- (.5, 1);
\draw (.75,0) -- (.75, 1);
\draw (1,0) -- (1, 1);
\draw (1.25,0) -- (1.25, 1);
\draw (0,.75) -- (1,.75);
\draw (.5,.25) -- (1.5,.25);
\draw (.625,1) -- (.625, .25);
\draw (.833,0) -- (.833,.75);
\draw (.916,.25) -- (.916,1);
\end{tikzpicture},
\begin{tikzpicture}[scale=1.5]
\filldraw[draw = red, ultra thick, fill=red!60] (.625,.25) -- (1.25,.25) -- (1.25,1) -- (.625,1) -- cycle;
\draw[red, ultra thick] (.625,.5) -- (1.25,.5);
\draw[red, ultra thick] (.916,.25) -- (.916,1);
\draw[red, ultra thick] (.75,.25) -- (.75,1);
\draw[red, ultra thick] (1,.25) -- (1,1);
\draw (0,0) -- (1.5,0) -- (1.5, 1) -- (0,1) -- cycle;
\draw (0,.5) -- (1.5, .5);
\draw (.25,0) -- (.25, 1);
\draw (.5,0) -- (.5, 1);
\draw (.75,0) -- (.75, 1);
\draw (1,0) -- (1, 1);
\draw (1.25,0) -- (1.25, 1);
\draw (0,.75) -- (1,.75);
\draw (.5,.25) -- (1.5,.25);
\draw (.625,1) -- (.625, .25);
\draw (.833,0) -- (.833,.75);
\draw (.916,.25) -- (.916,1);
\end{tikzpicture}
\end{tabular}\vspace{-.15cm}
\end{minipage}\hspace{4.25cm}
}\\
\subfloat[]{
\begin{tikzpicture}[scale=3.6]
\draw (0,0) -- (1.25,0) -- (1.25,.8) -- (0,.8) -- cycle;
\draw (.25,0) -- (.25,.8);
\draw (.5,0) -- (.5,.8);
\draw (.75,0) -- (.75,.8);
\draw (1,0) -- (1,.8);
\draw (0,.6) -- (.75,.6);
\draw (.5,.2) -- (1.25,.2);
\draw (.562, .2) -- (.562,.8);
\draw (.624, .6) -- (.624,0);
\draw (.686,.2) -- (.686,.8);
\end{tikzpicture}
}~
\subfloat[]{
\begin{minipage}{8.5cm}\vspace{-2.85cm}
\begin{tabular}{c}
\begin{tikzpicture}[scale=1.8]
\filldraw[draw=red, ultra thick, fill=red!60] (.25,0) -- (.75,0) -- (.75, .6) -- (.25, .6) -- cycle;
\draw[red, ultra thick]  (.5, 0) -- (.5, .6);
\draw[red, ultra thick] (.624, 0) -- (.624, .6);
\draw (0,0) -- (1.25,0) -- (1.25,.8) -- (0,.8) -- cycle;
\draw (.25,0) -- (.25,.8);
\draw (.5,0) -- (.5,.8);
\draw (.75,0) -- (.75,.8);
\draw (1,0) -- (1,.8);
\draw (0,.6) -- (.75,.6);
\draw (.5,.2) -- (1.25,.2);
\draw (.562, .2) -- (.562,.8);
\draw (.624, .6) -- (.624,0);
\draw (.686,.2) -- (.686,.8);
\end{tikzpicture},
\begin{tikzpicture}[scale=1.8]
\filldraw[draw=red, ultra thick, fill=red!60] (.25,0) -- (1,0) -- (1, .8) -- (.25, .8) -- cycle;
\draw[red, ultra thick]  (.5, 0) -- (.5, .8);
\draw[red, ultra thick] (.75, 0) -- (.75, .8);
\draw (0,0) -- (1.25,0) -- (1.25,.8) -- (0,.8) -- cycle;
\draw (.25,0) -- (.25,.8);
\draw (.5,0) -- (.5,.8);
\draw (.75,0) -- (.75,.8);
\draw (1,0) -- (1,.8);
\draw (0,.6) -- (.75,.6);
\draw (.5,.2) -- (1.25,.2);
\draw (.562, .2) -- (.562,.8);
\draw (.624, .6) -- (.624,0);
\draw (.686,.2) -- (.686,.8);
\end{tikzpicture},
\begin{tikzpicture}[scale=1.8]
\filldraw[draw=red, ultra thick, fill=red!60] (.5,0) -- (1,0) -- (1, .2) -- (.5, .2) -- cycle;
\draw[red, ultra thick]  (.75, 0) -- (.75, .2);
\draw[red, ultra thick] (.624, 0) -- (.624, .2);
\draw (0,0) -- (1.25,0) -- (1.25,.8) -- (0,.8) -- cycle;
\draw (.25,0) -- (.25,.8);
\draw (.5,0) -- (.5,.8);
\draw (.75,0) -- (.75,.8);
\draw (1,0) -- (1,.8);
\draw (0,.6) -- (.75,.6);
\draw (.5,.2) -- (1.25,.2);
\draw (.562, .2) -- (.562,.8);
\draw (.624, .6) -- (.624,0);
\draw (.686,.2) -- (.686,.8);
\end{tikzpicture},
\begin{tikzpicture}[scale=1.8]
\filldraw[draw=red, ultra thick, fill=red!60] (.5,.6) -- (.75,.6) -- (.75, .8) -- (.5, .8) -- cycle;
\draw[red, ultra thick]  (.562, .8) -- (.562, .6);
\draw[red, ultra thick] (.686, .8) -- (.686, .6);
\draw (0,0) -- (1.25,0) -- (1.25,.8) -- (0,.8) -- cycle;
\draw (.25,0) -- (.25,.8);
\draw (.5,0) -- (.5,.8);
\draw (.75,0) -- (.75,.8);
\draw (1,0) -- (1,.8);
\draw (0,.6) -- (.75,.6);
\draw (.5,.2) -- (1.25,.2);
\draw (.562, .2) -- (.562,.8);
\draw (.624, .6) -- (.624,0);
\draw (.686,.2) -- (.686,.8);
\end{tikzpicture}
\\
\begin{tikzpicture}[scale=1.8]
\filldraw[draw=red, ultra thick, fill=red!60] (.25,.6) -- (.686,.6) -- (.686, .8) -- (.25, .8) -- cycle;
\draw[red, ultra thick]  (.562, .8) -- (.562, .6);
\draw[red, ultra thick] (.5, .8) -- (.5, .6);
\draw (0,0) -- (1.25,0) -- (1.25,.8) -- (0,.8) -- cycle;
\draw (.25,0) -- (.25,.8);
\draw (.5,0) -- (.5,.8);
\draw (.75,0) -- (.75,.8);
\draw (1,0) -- (1,.8);
\draw (0,.6) -- (.75,.6);
\draw (.5,.2) -- (1.25,.2);
\draw (.562, .2) -- (.562,.8);
\draw (.624, .6) -- (.624,0);
\draw (.686,.2) -- (.686,.8);
\end{tikzpicture},
\begin{tikzpicture}[scale=1.85]
\filldraw[draw=red, ultra thick, fill=red!60] (.5,.2) -- (.686,.2) -- (.686, .6) -- (.5, .6) -- cycle;
\draw[red, ultra thick]  (.562, .2) -- (.562, .6);
\draw[red, ultra thick] (.624, .2) -- (.624, .6);
\draw (0,0) -- (1.25,0) -- (1.25,.8) -- (0,.8) -- cycle;
\draw (.25,0) -- (.25,.8);
\draw (.5,0) -- (.5,.8);
\draw (.75,0) -- (.75,.8);
\draw (1,0) -- (1,.8);
\draw (0,.6) -- (.75,.6);
\draw (.5,.2) -- (1.25,.2);
\draw (.562, .2) -- (.562,.8);
\draw (.624, .6) -- (.624,0);
\draw (.686,.2) -- (.686,.8);
\end{tikzpicture},
\begin{tikzpicture}[scale=1.8]
\filldraw[draw=red, ultra thick, fill=red!60] (.562,.2) -- (.75,.2) -- (.75, .6) -- (.562, .6) -- cycle;
\draw[red, ultra thick]  (.686, .2) -- (.686, .6);
\draw[red, ultra thick] (.624, .2) -- (.624, .6);
\draw (0,0) -- (1.25,0) -- (1.25,.8) -- (0,.8) -- cycle;
\draw (.25,0) -- (.25,.8);
\draw (.5,0) -- (.5,.8);
\draw (.75,0) -- (.75,.8);
\draw (1,0) -- (1,.8);
\draw (0,.6) -- (.75,.6);
\draw (.5,.2) -- (1.25,.2);
\draw (.562, .2) -- (.562,.8);
\draw (.624, .6) -- (.624,0);
\draw (.686,.2) -- (.686,.8);
\end{tikzpicture},
\begin{tikzpicture}[scale=1.8]
\filldraw[draw=red, ultra thick, fill=red!60] (.562,.2) -- (1,.2) -- (1, .8) -- (.562, .8) -- cycle;
\draw[red, ultra thick]  (.75, .2) -- (.75, .8);
\draw[red, ultra thick] (.686, .2) -- (.686, .8);
\draw (0,0) -- (1.25,0) -- (1.25,.8) -- (0,.8) -- cycle;
\draw (.25,0) -- (.25,.8);
\draw (.5,0) -- (.5,.8);
\draw (.75,0) -- (.75,.8);
\draw (1,0) -- (1,.8);
\draw (0,.6) -- (.75,.6);
\draw (.5,.2) -- (1.25,.2);
\draw (.562, .2) -- (.562,.8);
\draw (.624, .6) -- (.624,0);
\draw (.686,.2) -- (.686,.8);
\end{tikzpicture}
\end{tabular}\vspace{-.15cm}
\end{minipage}\hspace{1.75cm}
}
\caption{In (a) is shown an LR-mesh providing LR B-splines of degree (3,3) in an equivalent arrangement of the LR B-splines of bidegree (2,2) on the mesh in Figure \ref{8LR} (a). In (b) are shown the supports of the eight B-splines in the linear dependence relation. In (c)-(d) are shown the same for bidegree (4,4), in (e)-(f) for bidegree (3,1) and (g)-(h) for bidegree (2,0).}\label{deg34LR}
\end{figure}
\end{ex}
We stress that mesh $\cM + \gamma$ in Figure \ref{8LR} (a) is obtained by refining the mesh in Figure \ref{7} (a) considered in Example \ref{exMS}. What happens is that with the insertion of a new split, the MS B-spline in the center of mesh $\cM$, $B^6$, is split into two MS B-splines that can now be obtained through the knot insertion procedure. 


\section{Improvement of the Peeling Algorithm}\label{peel}
The Peeling Algorithm introduced in \cite{dokken2013polynomial} is a tool to check if the LR B-splines on a given LR-mesh are linearly independent. However it does not handle every possible configuration. In this section, we briefly recall it and we show how it can be improved to sort out more cases. 
\begin{dfn}
An element of the box-partition $\mathcal{E}$ is \textbf{overloaded} if it is covered by more B-splines than necessary for spanning the corresponding polynomial space $\Pi_{\pmb{p}}$. We call a B-spline \textbf{overloaded} if all the elements in its support are overloaded. 
\end{dfn}
\begin{oss}
In 2D, if we consider bidegree $\pmb{p} = (p_1, p_2)$, then the necessary number of B-splines for spanning $\Pi_{\pmb{p}}$ on a given element is $(p_1+1)(p_2+1)$.
\end{oss}
An extra B-spline, in a linear dependence, can be removed without changing spanning properties over the elements of $\mathcal{E}$ it is covering with its support. So, only overloaded B-splines occur in linear dependencies. The Peeling algorithm is based on the fact that a linear dependence relation has to involve at least two overloaded B-splines. Therefore, if on an element there is the support of only one overloaded B-spline, such B-spline cannot be active in a linear dependence. This simple observation is the basis of the Peeling Algorithm.\\
\begin{algorithm}[H]
\vspace{.2cm}
\nonl\nosemic \textbf{Peeling Algorithm}
\dosemic\\
\vspace{.1cm}
From the set of LR splines $\mathcal{B}^{\mathcal{LR}}(\cM)$ create the set $\mathcal{B}^O$ of overloaded B-splines\;
Let $\mathcal{E}^O$ be the elements of $\mathcal{E}$ in the supports of the B-splines in $\mathcal{B}^O$\;
Initialization of a subset $\mathcal{B}^O_1$ of $\mathcal{B}^O$ we are going to define, $\mathcal{B}^O_1 = \emptyset$\;
\For{every element $\beta$ in $\mathcal{E}^O$}{
\If{only one B-spline $B$ of $\mathcal{B}^O$ has $\beta$ in its support}{$\mathcal{B}^O_1 = \mathcal{B}^O_1\cup\{B\}$}}
\eIf{$\mathcal{B}^O\backslash \mathcal{B}^O_1=\emptyset$}{linear independence.}{
\If{$\mathcal{B}^O_1=\emptyset$}{break, but might have linear dependence.}
$\mathcal{B}^O = \mathcal{B}^O\backslash \mathcal{B}^O_1$\;
Go to 2\;
}
\end{algorithm}
The implementation of it is described in \cite{dokken2013polynomial} in terms of matrices.

\begin{figure}[ht!]
\subfloat[]{
\begin{tikzpicture}[scale = 3.8]
\draw[white] (.5,1.05);
\fill[blue!50] (.2,.2) -- (.2,.8) -- (.8,.8) -- (.8, .2) -- cycle;
\draw (0,0) -- (1,0) -- (1,1) -- (0,1) -- cycle;
\draw (0,.2) -- (1,.2);
\draw (0,.4) -- (1,.4);
\draw (0,.6) -- (1,.6);
\draw (0,.8) -- (1,.8);
\draw (.2,0) -- (.2,1);
\draw (.4,0) -- (.4,1);
\draw (.6,0) -- (.6,1);
\draw (.8,0) -- (.8,1);
\draw (0, .7) -- (.7,.7);
\draw (.3,.7) -- (.3,0);
\draw (.3,.3) -- (1,.3);
\draw (.7,.3) -- (.7,1);
\draw[white] (.5,-.1);
\end{tikzpicture}
}
\subfloat[]{
\begin{minipage}{5cm}\vspace{-4.15cm}
\begin{tabular}{c}
\begin{tikzpicture}[scale=2]
\filldraw[draw=red, ultra thick, fill=red!50] (0,0) -- (.66,0) -- (.66,.83) -- (0,.83) -- cycle;
\draw[red, ultra thick] (0,.33) -- (.66,.33);
\draw[red, ultra thick] (0,.66) -- (.66,.66);
\draw[red, ultra thick] (.33,0) -- (.33,.83);
\draw[red, ultra thick] (.165,0) -- (.165,.83);
\draw (0,0) -- (1,0) -- (1,1) -- (0,1) -- cycle;
\draw (0,.33) -- (1,.33);
\draw (0,.66) -- (1,.66);
\draw (.33,0) -- (.33,1);
\draw (.66,0) -- (.66,1);
\draw (0,.83) -- (.83,.83);
\draw (.165,0) -- (.165,.83);
\draw (.165,.165) -- (1,.165);
\draw (.83,.165) -- (.83,1);
\end{tikzpicture},
\begin{tikzpicture}[scale=2]
\filldraw[draw=red, ultra thick, fill=red!50] (.165,0) -- (1,0) -- (1,.66) -- (.165, .66) -- cycle;
\draw[red, ultra thick] (.165,.33) -- (1,.33);
\draw[red, ultra thick] (.165,.165) -- (1,.165);
\draw[red, ultra thick] (.33,0) -- (.33,.66);
\draw[red, ultra thick] (.66,0) -- (.66,.66);
\draw (0,0) -- (1,0) -- (1,1) -- (0,1) -- cycle;
\draw (0,.33) -- (1,.33);
\draw (0,.66) -- (1,.66);
\draw (.33,0) -- (.33,1);
\draw (.66,0) -- (.66,1);
\draw (0,.83) -- (.83,.83);
\draw (.165,0) -- (.165,.83);
\draw (.165,.165) -- (1,.165);
\draw (.83,.165) -- (.83,1);
\end{tikzpicture},
\begin{tikzpicture}[scale=2]
\filldraw[draw=red, ultra thick, fill=red!50] (.33,.165) -- (1,.165) -- (1,1) -- (.33, 1) -- cycle;
\draw[red, ultra thick] (.33,.33) -- (1,.33);
\draw[red, ultra thick] (.33,.66) -- (1,.66);
\draw[red, ultra thick] (.66,.165) -- (.66,1);
\draw[red, ultra thick] (.83,.165) -- (.83,1);
\draw (0,0) -- (1,0) -- (1,1) -- (0,1) -- cycle;
\draw (0,.33) -- (1,.33);
\draw (0,.66) -- (1,.66);
\draw (.33,0) -- (.33,1);
\draw (.66,0) -- (.66,1);
\draw (0,.83) -- (.83,.83);
\draw (.165,0) -- (.165,.83);
\draw (.165,.165) -- (1,.165);
\draw (.83,.165) -- (.83,1);
\end{tikzpicture}\\
\begin{tikzpicture}[scale=2]
\filldraw[draw=red, ultra thick, fill=red!50] (0,.33) -- (.83,.33) -- (.83,1) -- (0, 1) -- cycle;
\draw[red, ultra thick] (0,.83) -- (.83,.83);
\draw[red, ultra thick] (0,.66) -- (.83,.66);
\draw[red, ultra thick] (.33,.33) -- (.33,1);
\draw[red, ultra thick] (.66,.33) -- (.66,1);
\draw (0,0) -- (1,0) -- (1,1) -- (0,1) -- cycle;
\draw (0,.33) -- (1,.33);
\draw (0,.66) -- (1,.66);
\draw (.33,0) -- (.33,1);
\draw (.66,0) -- (.66,1);
\draw (0,.83) -- (.83,.83);
\draw (.165,0) -- (.165,.83);
\draw (.165,.165) -- (1,.165);
\draw (.83,.165) -- (.83,1);
\end{tikzpicture},
\begin{tikzpicture}[scale=2]
\filldraw[draw=red, ultra thick, fill=red!50] (0,0) -- (1,0) -- (1,1) -- (0, 1) -- cycle;
\draw[red, ultra thick] (0,.33) -- (1,.33);
\draw[red, ultra thick] (0,.66) -- (1,.66);
\draw[red, ultra thick] (.33,0) -- (.33,1);
\draw[red, ultra thick] (.66,0) -- (.66,1);
\draw (0,0) -- (1,0) -- (1,1) -- (0,1) -- cycle;
\draw (0,.33) -- (1,.33);
\draw (0,.66) -- (1,.66);
\draw (.33,0) -- (.33,1);
\draw (.66,0) -- (.66,1);
\draw (0,.83) -- (.83,.83);
\draw (.165,0) -- (.165,.83);
\draw (.165,.165) -- (1,.165);
\draw (.83,.165) -- (.83,1);
\end{tikzpicture}
\end{tabular}
\end{minipage}\hspace{2cm}
}
\subfloat[]{
\begin{tikzpicture}[scale=3.75]
\draw (0,0) -- (1,0) -- (1,1) -- (0,1) -- cycle;
\draw (0,.33) -- (1,.33);
\draw (0,.66) -- (1,.66);
\draw (.33,0) -- (.33,1);
\draw (.66,0) -- (.66,1);
\draw (0,.83) -- (.83,.83);
\draw (.165,0) -- (.165,.83);
\draw (.165,.165) -- (1,.165);
\draw (.83,.165) -- (.83,1);
\draw (.082,.165) node{2};
\draw (.83,.082) node{2};
\draw (.165,.905) node{2};
\draw (.905,.83) node{2};
\draw (.247,.082) node{3};
\draw (.247,.247) node{3};
\draw (.5,.247) node{4};
\draw (.5,.082) node{3};
\draw (.752,.247) node{3};
\draw (.905,.247) node{3};
\draw (.082,.5) node{3};
\draw (.247,.5) node{4};
\draw (.5,.5) node{5};
\draw (.752,.5) node{4};
\draw (.905,.5) node{3};
\draw (.082,.752) node{3};
\draw (.247,.752) node{3};
\draw (.5,.752) node{4};
\draw (.752,.752) node{3};
\draw (.5,.905) node{3};
\draw (.752,.905) node{3};
\fill (.165,.165) circle (.02);
\fill (.165,.83) circle (.02);
\fill (.83,.83) circle (.02);
\fill (.83,.165) circle (.02);
\draw[white] (.5,-.1);
\end{tikzpicture}
}
\caption{Consider bidegree (2,2). In the highlighted region in (a) there are the supports of five overloaded LR B-splines, depicted in (b). The numbers in the elements of the region, reported in (c), indicate how many supports of these B-splines are on them. The highlighted vertices are the T-vertices corresponding to pair of knots of the overloaded LR B-splines.}\label{pe}
\end{figure}
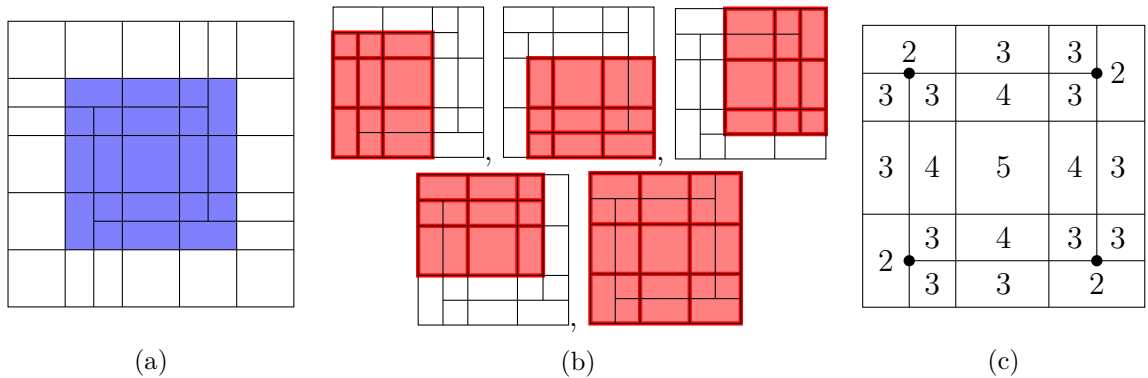

However, it might happen that every element of $\mathcal{E}^O$ is shared but yet the overloaded LR B-splines are linearly independent. An example is reported in Figure \ref{pe}.
We consider bidegree (2,2) and an LR-mesh of multiplicity one. In the highlighted region in (a) there are the supports of five LR B-splines, reported in (b), that form the collection $\mathcal{B}^O$ of the algorithm. Then, for each element of the box-partition in such region we count how many of these supports are on it. If an element is only in one support, the corresponding B-spline is placed in the subcollection $\mathcal{B}^O_1$ of the algorithm. From (c), we see that $\mathcal{B}^O_1 = \emptyset$. Therefore, the algorithm stops without answering whether the LR B-splines on the mesh are linearly independent or not. However, if we now look at the T-vertices in the region, highlighted in (c), we see that none of them is shared, as pair of knots, in two or more B-splines of $\mathcal{B}^O$. Since the necessary condition for linear dependency Corollary \ref{Tver} is not satisfied, we can conclude that the LR B-splines on the mesh are linearly independent.

The Peeling Algorithm can therefore be improved inserting in $\mathcal{B}_1^O$ also the B-splines of $\mathcal{B}^O$ that have an exclusive T-vertex as pair of knots.
\section{Conclusions, conjectures and future work}\label{conc}
In this work we have identified features of the mesh to have a linear dependence relation in the MS and LR B-spline sets for any bidegree $\pmb{p}$. Namely, if the union of the supports of the B-splines involved in the linear dependence relation is called $\mathcal{R}$,\begin{itemize}
\item There are nested B-splines at the corners of $\mathcal{R}$. 
\item Every relevant T-vertex is shared.
\end{itemize}  
These conditions are necessaries for a linear dependence relation on any LR-mesh (with meshlines of arbitrary multiplicities). Moreover, we have proved that the minimal number of MS B-splines needed for a linear dependence relation is 6 while for the LR B-splines is 8. These numbers are sharp for any bidegree $\pmb{p} = (p_1,p_2)$ with $p_k\geq 1$ for the MS B-splines and $p_k \geq 2$ for the LR B-splines for some $k \in \{1,2\}$. 

While, when $(p_1,p_2)=(0,0)$ both the MS and the LR B-splines are linearly independent on any LR-mesh. Finally when $(p_1, p_2) = (0,1), (1,0)$ or $(1,1)$, it seems not possible to have a linear dependence relation in the LR B-spline set.

In our future work, we would like to classify the meshes with a linear dependence relation involving this minimal number of MS B-splines. The number of possible cases would then be dependent on the bidegree chosen. Our conjecture is that every possible configuration of linear dependency is a refinement of one of such cases. Note that this is what happens in the meshes of Example \ref{exLR}.

\bibliographystyle{elsarticle-num-names}
\bibliography{biblio-tor}

\end{document}